\documentclass[11pt]{amsart}
\usepackage[hmargin=3cm,vmargin=3cm]{geometry}

\usepackage[latin1]{inputenc}
\usepackage{xspace,amssymb,amsfonts,euscript}
\usepackage{amsthm,amsmath,amscd,stmaryrd,latexsym}
\usepackage{palatino, euscript}

\usepackage{graphicx}
\usepackage[all]{xy}
\usepackage{tikz}
\usetikzlibrary{calc, matrix, arrows, decorations.pathmorphing}
\usepackage[dvips]{epsfig}
\usepackage{pinlabel}
\usepackage{psfrag}

\newcommand{\ig}[2]{\vcenter{\xy (0,0)*{\includegraphics[scale=#1]{fig/#2}} \endxy}}

\IfFileExists{comments.sty}{\usepackage{comments}}

\RequirePackage{color}
\definecolor{myred}{rgb}{0.75,0,0}
\definecolor{mygreen}{rgb}{0,0.5,0}
\definecolor{myblue}{rgb}{0,0,0.65}

\RequirePackage{ifpdf}
\ifpdf
  \IfFileExists{pdfsync.sty}{\RequirePackage{pdfsync}}{}
  \RequirePackage[pdftex,
   colorlinks = true,
   urlcolor = myblue, 
   citecolor = mygreen, 
   linkcolor = myred, 
   pagebackref,
   bookmarksopen=true]{hyperref}
\else
  \RequirePackage[hypertex]{hyperref}
\fi

\newtheorem{thm}{Theorem}[section]
\newtheorem{lemma}[thm]{Lemma}

\newtheorem{prop}[thm]{Proposition}
\newtheorem{cor}[thm]{Corollary}
\newtheorem{claim}[thm]{Claim}  
\newtheorem{conj}[thm]{Conjecture}

\newtheorem*{prop*}{Proposition}

\theoremstyle{definition}
\newtheorem{defn}[thm]{Definition}

\newtheorem{example}[thm]{Example}

\theoremstyle{remark}
\newtheorem{remark}[thm]{Remark}

\numberwithin{equation}{section}

    \def\CM{{\mathbb{C}}}

  \def\gg{{\mathfrak g}}

    \def\QM{{\mathbb{Q}}}

    \def\ZM{{\mathbb{Z}}}

    \def\CC{{\mathcal{C}}}
    
  \def\eb{{\mathbf e}}  
  \def\fb{{\mathbf f}}  
  \def\gb{{\mathbf g}}  
  \def\hb{{\mathbf h}}

    \def\UC{{\mathcal{U}}}

\def\a{\alpha}
\def\b{\beta}
\def\g{\gamma}
\def\G{\Gamma}
\def\d{\delta}

\def\e{\varepsilon}
\def\k{\kappa}
\def\l{\lambda}
\def\L{\Lambda}

\def\s{\sigma}

\def\t{\tau}
\def\w{\omega}
\def\Om{\Omega}
\def\z{\zeta}

\let\phi=\varphi

\usepackage{bbm}

\def\R{{\mathbbm R}}

\def\1{\mathbbm{1}}

\newcommand{\un}{\underline}

\newcommand{\ot}{\otimes}

\newcommand{\co}{\colon}

\renewcommand{\to}{\rightarrow}

\renewcommand{\sl}{\mathfrak{sl}}
\newcommand{\gl}{\mathfrak{gl}}

\newcommand{\refequal}[1]{\xy {\ar@{=}^{#1}
(-1,0)*{};(1,0)*{}};
\endxy}

\DeclareMathOperator{\Hom}{{\rm Hom}}

\DeclareMathOperator{\End}{{\rm End}}

\DeclareMathOperator{\Ext}{{\rm Ext}}

\DeclareMathOperator{\Lad}{\rm{Lad}}

\DeclareMathOperator{\Fund}{\bf{Fund}}
\DeclareMathOperator{\Rep}{\bf{Rep}}
\DeclareMathOperator{\LL}{\mathbb{LL}}

\DeclareMathOperator{\wt}{\rm{wt}}

\newcommand{\qbinom}[2]{{\genfrac{[}{]}{0pt}{}{#1}{#2}}}

\renewcommand{\top}{+}

\begin{document}

\begin{abstract} Morphisms between tensor products of fundamental representations of the quantum group $U_q(\sl_n)$ are described by the $\sl_n$-webs of Cautis-Kamnitzer-Morrison. Using
these webs, we provide an explicit, root-theoretic formula for the local intersection forms attached to each summand of the tensor product of an irreducible representation with a
fundamental representation. We prove this formula for $\sl_n$ with $n \le 4$, and conjecture that it holds for all $n$.

Given two sequences of fundamental weights which sum to the same dominant weight, the clasp is the morphism between the corresponding tensor products which projects to the top
indecomposable summand. Using our computation of intersection forms, we provide a recursive, ``triple-clasp expansion'' formula for clasps.

In addition, we describe the cellular structure on $\sl_n$-webs, and prove that $\sl_n$-webs are an integral form for tilting modules of quantum groups.
\end{abstract}

\title{Light ladders and clasp conjectures}

\author{Ben Elias} \address{University of Oregon, Eugene}

\maketitle

\tableofcontents

\section{Introduction}
\label{sec-intro}

Consider a complex semisimple Lie algebra $\gg$, its category $\Rep(\gg)$ of finite dimensional representations, and the full subcategory $\Fund(\gg)$ whose objects consist of tensor
products of fundamental representations. In some sense, $\Fund(\gg)$ contains all the information needed to study $\Rep(\gg)$, as every irreducible representation is a summand of some
object in $\Fund(\gg)$. Moreover, $\Fund(\gg)$ has a number of advantages which make it more accessible to algebraic study. The upshot is that, when $\gg = \sl_n$, there is a presentation
of the monoidal category $\Fund(\gg)$ by generators and relations due to Cautis-Kamnitzer-Morrison \cite{CKM}, using the language of planar diagrammatics. A morphism between fundamental
representations is encoded as a certain kind of oriented graph known as an $\sl_n$-web.

The idea to describe morphisms between fundamental representations using webs dates back to the foundational work of Kuperberg \cite{Kupe}. Earlier, morphisms in $\Fund(\sl_2)$ had been
described by the famous Temperley-Lieb algebra \cite{TemLie}. Kuperberg generalized this to rank 2 Lie algebras. Morrison \cite{Morr07} attempted to describe webs for
$\sl_n$ in a similar fashion, giving generators and a correct list of relations (together with some redundant ones), but he was unable to prove at that time that the list of relations was
sufficient. One should also mention Dongseok Kim \cite{Kim03, Kim07} as another early developer of the theory, who also gave a conjectural presentation for $\sl_4$.
Cautis-Kamnitzer-Morrison \cite{CKM} proved Morrison's conjectural results by a clever representation-theoretic transformation, using skew Howe duality to relate $\sl_n$-webs with $m$
inputs and outputs to the universal enveloping algebra of $\sl_m$. To date, nothing is known outside of type $A$ and rank $2$.

A brief motivation of $\sl_n$-webs goes as follows. For $0 \le i \le n$, let $V_{\w_i}$ denote the exterior product $\Lambda^i (\CM^n)$. When $1 \le i \le n-1$, these are the fundamental
representations of $\sl_n$. An $\sl_n$-web (resp. extended $\sl_n$-web) is a kind of oriented graph where the edges are labeled by $i$ for $1 \le i \le n-1$ (resp. for $0 \le i \le n$),
where an edge labeled $i$ represents either the representation $V_{\w_i}$ or its dual, depending on the orientation. The structure of duality is encoded in four maps (trace, cotrace,
evaluation, and coevaluation) which are depicted as oriented cups and caps. When $i + j = k$, there is an obvious map $\Lambda^i(\CM^n) \ot \Lambda^j(\CM^n) \to \Lambda^k(\CM^n)$, given by
the wedge product. This is depicted by a trivalent vertex, with edges labeled $i$, $j$, and $k$. There is also an isomorphism between $\Lambda^i(\CM^n)$ and $(\Lambda^{n-i}(\CM^n))^*$,
which is depicted as an orientation-switching tag. Together these generate all morphisms between tensor products of $V_{\w_i}$ and their duals. The interesting question, answered finally
by Cautis-Kamnitzer-Morrison, is what relations hold between these graphs.

All the results mentioned above apply more generally to the quantum deformation. Let $U_q(\gg)$ denote the $\ZM[q,q^{-1}]$-integral form of the quantum group of $\gg$, and let $U'_q(\gg) =
U_q(\gg) \ot_{\ZM[q,q^{-1}]} \QM(q)$. Let us omit $\gg$ from the notation when it is understood. There are $q$-deformations $\Rep_q$ and $\Rep'_q$, $\Fund_q$ and $\Fund'_q$, which are
categories of representations of $U_q$ and $U'_q$ respectively. One goal is to present $\Fund_q$ as a $\ZM[q,q^{-1}]$-algebra, something which was also achieved by the Temperley-Lieb
algebra for $\sl_2$. The skew Howe duality trick used in Cautis-Kamnitzer-Morrison \cite{CKM} only works over $\QM(q)$, so that they were able to prove that their generators and relations
described $\Fund'_q$ over $\QM(q)$, but could not prove anything about the integral form. They were kind enough in \cite{CKM} to provide relations which, while redundant over $\QM(q)$, are
not redundant over $\ZM[q,q^{-1}]$, and it is safe to say that they conjectured their relations to suffice over $\ZM[q,q^{-1}]$ as well.

Henceforth we usually omit the $q$ from $\Fund_q$ and $\Rep_q$, as we will always be talking about the $q$-deformation.

The first half of this paper will be an in-depth analysis of the category of $\sl_n$-webs described in \cite{CKM}, with the ultimate goal of constructing a basis for morphism spaces over
$\ZM[q,q^{-1}]$. We call it the \emph{double ladders basis}, because elements are built from pairing two \emph{light ladders} together. In particular, we give a purely diagrammatic proof
that double ladders are a basis (without needing to use skew Howe duality), and proving that $\sl_n$-webs describe the integral form $\Fund_q (\sl_n)$. Double ladders are a cellular basis, and are adapted to the monoidal structure of the category in a particular way, which shall be described further in the introduction.

\begin{remark} The elements of the double ladders basis are webs, not linear combinations of webs. It was proven by Khovanov-Kuperberg \cite{KhoKup} that a web basis can not correspond,
under skew Howe duality, to the canonical basis of $U_q(\sl_m)$. \end{remark}

\begin{remark} A basis for $\sl_n$-webs was also constructed by Bruce Fontaine in \cite{Font12}. (Thanks to Joel Kamnitzer for bringing this to my attention.) His proof involves the
geometry of the affine Grassmannian, so that it requires the base ring to be $\CM$, with $q=1$. He identifies all morphism spaces, using adjunction, with webs on a disk, or morphisms from
a tensor product to the trivial representation. Using adjunction in this way is reasonable if one is not interested in certain features of the category, but we will not do it because it
obfuscates the cellular structure. It seems very likely that our bases for this particular morphism space agree (compare his $T_{\w}$ in Theorem 3.1 to our elementary light ladder
\eqref{elementarylightladder}).

Bruce Fontaine's work is related to earlier works of Bruce Westbury, who constructed bases for morphisms between certain tensor products in type $A$ \cite{WestWebs}, and also studied type
$G_2$ \cite{WestG2}. The concept of a cellular category is also due to Westbury \cite{WestCC}. \end{remark}

\begin{remark} That $\sl_n$-webs are a cellular category (at least over $\QM(q)$) was also recently proven by Andersen-Stroppel-Tubbenhauer \cite{AST}, which appeared as this paper was in
preparation. We will discuss their work in greater detail below. \end{remark}

It is all well and good to say that $\Fund(\gg)$ contains the information needed to study $\Rep(\gg)$. After all, the categories are Morita-equivalent (in the semisimple case). It is
another thing altogether to actually try to recover $\Rep(\gg)$ from $\Fund(\gg)$ in an explicit fashion. One needs to compute the idempotent inside a tensor product of fundamental
representations which projects to each of its irreducible summands, which is an extremely difficult task.

Each tensor product of fundamental representations has a unique irreducible summand which does not appear in a shorter tensor product (i.e. the one with the highest highest weight), which
we call the \emph{top} summand. In the literature, the idempotent projecting to the top summand is called a \emph{clasp} or a \emph{generalized Jones-Wenzl projector}, after Jones
\cite{Jon3} and Wenzl \cite{Wenzl} who independently gave formulae for this idempotent in the $\sl_2$ case. The word clasp comes from Kuperburg \cite{Kupe} (who called them \emph{internal
clasps}). There are a number of recursive formulas for clasps in the $\sl_2$ case, and these formulas are called \emph{single-clasp formulas}, \emph{double-clasp formulas}, etcetera, based
on the precise form of the recursion (this will be discussed below). Dongseok Kim \cite{Kim03, Kim07} has investigated the problem in rank 2, and derived a number of single-clasp
formulas. Beyond that, before this paper, essentially nothing was known.

\begin{remark} Sabin Cautis \cite{CautisClasp} has an entirely different approach to clasps, expressing them as an infinite power of the full twist (a convergent limit) coming from a braid
group action. This lovely paper does much more, using clasps to study categorified skew Howe duality. In theory, this description of the clasp should also allow one to compute it, but more
work needs to be done in this direction. \end{remark}

The same plethyism principle which makes the double ladders basis adapted to the monoidal structure also provides a tautological form for the \emph{triple-clasp recursive formula}, which
is the one we shall study in this paper. What remains is to compute certain scalars known as \emph{local intersection forms} or \emph{cellular forms} in different contexts, which appear as
coefficients in this recursion. We conjecture an explicit formula for the local intersection forms involved in the triple clasp formula, which describes them as a product of ratios of
quantum numbers. The form of this conjecture is superficially similar to the quantum Weyl dimension formula, but different enough to be entirely new and unfamiliar. By computing these
numbers, one can begin to analyze precisely what goes wrong when $q$ is specialized to a root of unity.

We develop some preliminary techniques to study this conjecture, and prove it for $\sl_n$ with $n \le 4$. These results are new even for $n=3$. We feel that the techniques of this paper,
applied with even more fervor and elbow grease, should be sufficient to prove this conjecture in general. However, a purely representation-theoretic explanation would be a desirable
alternative.

\subsection{Monoidal cellular categories}
\label{subsec-monoidalcellular}

The classic Artin-Wedderburn theorem states that a semisimple algebra is isomorphic to a product of matrix algebras. In categorical language, the implication is that, in a semisimple
category, every morphism is a linear combination of morphisms which factor as a projection to an irreducible summand, followed by an inclusion. This \emph{factorization principle} is what
underlies the notion of an (object-adapted) cellular category.

The category $\Rep'$ is semisimple, with irreducibles $V_\l$ parametrized by the set $\L$ of dominant weights. Thus every morphism between representations, and in particular every
morphism in $\Fund'$, should (be a linear combination of morphisms which) factor as a projection to some $V_\l$ followed by an inclusion. Morphisms which factor through $V_\l$ are
orthogonal to those which factor through $V_\mu$, for $\l \ne \mu$. However, as mentioned above, the projection maps from a tensor product of fundamental representations to an irreducible
module are hard to describe. (Rather, we do not even have a language to describe morphisms between arbitrary representations. But we do have a language to describe morphisms between tensor
products in $\Fund'$; it is the composition of projections with inclusions which live in $\Fund'$, and are hard to describe.) More signficantly, such projection maps can not be described
over the integral form $\ZM[q,q^{-1}]$, because $\Rep$ is not semisimple, and the projection maps may not exist.

However, it turns out that what $\sl_n$-webs can encode easily are morphisms which factor through $V_\l$, plus terms which factor through $V_{\mu}$ for $\mu < \l$ in the dominance order.
These morphisms will be our double ladders basis. They are defined in the integral form $\Fund$, and (after base change to $\QM(q)$) are related to the semisimple factorized maps by an
upper triangular change of basis. Note that, again, we can not visualize what it means for a morphism to factor through $V_\l$ because this irreducible object is not inside $\Fund$.
However, we can assert that our morphism factors through a tensor product of fundamental representations which has $V_\l$ as its top summands; these two objects agree modulo summands which
are lower in the dominance order. This generalization of a product of matrix algebras, where morphisms ``factor through irreducibles modulo lower terms," is roughly the idea of a cellular
category.

For the rest of this introduction, ``lower terms'' will refer to morphisms which factor through summands which are lower in the dominance order than the $\l$ we are interested in. Note
however, that the summand $V_\l$ need not exist in $\Rep$ or in $\sl_n$-webs, but only in $\Rep'$. We will talk about morphisms which project to summands modulo lower terms; these are
morphisms in the integral form of $\sl_n$-webs which, after base change, have the desired properties in $\Rep'$.

\begin{remark} Many things are more simply described if one is willing to allow lower terms. For example, $\Fund$ is a symmetric monoidal category, and the braiding is defined by a nasty
linear combination of webs \cite[Corollary 6.2.3]{CKM}. However, certain simple webs, which we call \emph{neutral ladders} in this paper, provide morphisms which agree with the braiding on
the top summand, and thus only disagree modulo lower terms. \end{remark}

Strictly speaking, the original notions of a cellular algebra \cite{GraLeh} and a cellular category \cite{WestCC} do not involve the idea of factorization. There is a special basis
$\{c^\l_{S,T}\}$ parametrized by $\l \in \L$, a poset of \emph{cells}, and by $S$ and $T$ in some sets attached to each cell. One can think of a matrix algebra, which has a single cell,
and where the basis is parametrized by a choice $T$ of column and a choice $S$ of row. What we have described above is a cellular category where the cellular basis factors: $c^\l_{S,T} =
c^\l_S \circ c^\l_T$ into two ``half-basis'' maps, where $c^\l_S$ (resp. $c^\l_T$) is a morphism from (resp. to) some object attached to the cell $\l$. We call this a \emph{strictly
object-adapted cellular category} or \emph{SOACC}. The axiomatics of SOACCs have appeared in joint work with Lauda \cite{ELauda}. There is a very simple criterion \cite[Lemma 2.8]{ELauda}
by which one can verify that a category which looks like an SOACC (i.e. has a factorized basis) is in fact cellular, a criterion which could not be as simple for general cellular
categories.

\begin{remark} Most cellular algebras in the literature are endomorphism rings inside cellular categories (which is true tautologically, but we refer to more interesting examples). For
example, the Temperley-Lieb algebras are endomorphism rings inside the Temperley-Lieb category. Most cellular categories in the literature are also strictly object-adapted cellular
categories. Some examples we have in mind are: Temperley-Lieb algebras and similar algebras (like blob algebras), Soergel diagrammatics \cite{EWGr4sb}, and $q$-Schur algebras. Ironically,
the most famous example of a cellular algebra, the Hecke algebra of the symmetric group, is not known to come from an SOACC. \end{remark}

\begin{remark} In \cite{ELauda} we generalize SOACCs to define \emph{fibered cellular categories}. The analogy here is that the endomorphism ring of an ``irreducible'' $V_\l$ (let us now
say ``indecomposable'') is not just the scalars, but some (graded) local ring $A_\l$. Thus morphisms should factor as a composition $c^\l_T \circ a \circ c^\l_S$, projecting to $V_\l$,
applying an endomorphism, and then including from $V_\l$. Most of the generalities discussed here and below should apply to monoidal fibered cellular categories. For examples, see
\cite{ELauda}. \end{remark}

\begin{remark} There is one more structure attached to cellular categories and SOACCs which we have ignored so far, a duality involution. For $\sl_n$-webs, this comes from flipping a
diagram upside-down and reversing the orientation. \end{remark}

We identify an object in $\Fund$ (or the category of $\sl_n$-webs) with a sequence $\un{w} = (w_1, w_2, \ldots, w_d)$ of fundamental weights, with $w_k = \w_{i_k}$. We let $V_{\un{w}}$
denote the tensor product $V_{w_1} \ot V_{w_2} \ot \cdots \ot V_{w_d}$. In $\Rep'$ (that is, after base change to $\QM(q)$) the indecomposable summands of $V_{\un{w}}$ of the form $V_\l$
are parametrized by the set of miniscule Littelmann paths $E(\un{w}, \l)$ (more on this shortly). In summary, one of the things we accomplish in this paper is an explicit description of a
web, which we call a \emph{light ladder}, one for each miniscule Littelmann path in $E(\un{w}, \l)$, which gives a map from $V_{\un{w}}$ to $V_{\un{x}_\l}$ for some sequence $\un{x}_\l$
where the weights add up to $\l$, and which agrees with projection to the corresponding summand modulo lower terms.

\begin{remark} In fact, an object in $\Fund$ or $\sl_n$-webs is a sequence consisting of fundamental representations \emph{or their duals}. However, each dual of a fundamental
representation is isomorphic to some other fundamental representation (by the orientation-switching tag morphism). Tensor products of fundamentals, without their duals, form an essentially
surjective subcategory $\Fund^+$, which is actually our main object of study. We will ignore the distinction between $\Fund$ and $\Fund^+$ in the rest of this introduction. \end{remark}

Before moving on to a discussion of the monoidal structure, let us comment briefly on the recent work of Andersen-Stroppel-Tubbenhauer \cite{AST}. This elegant paper gives an explicit
representation-theoretic construction of the cellular algebra structure on the endomorphism ring of any \emph{tilting module}. For example, the tilting modules of the integral form $U_q$
are precisely the summands of tensor products of fundamental representations, so that the Karoubi envelope of $\Fund$ is naturally equivalent to the category of tilting modules. Tilting
modules have a filtration by standard modules and a filtration by costandard modules, and the cellular basis constructed in \cite{AST} also factors, via morphisms to and from a standard
module. However, the standard modules being factored through are not themselves tilting, so that their construction does not immediately lend itself to the formalism of object-adapted
cellular categories. It seems very likely that one could adjust their definition and proofs to define a cellular basis which factors through tilting modules instead, thus producing an
SOACC.

Their paper \cite{AST} achieves its goal: a sweeping result about tilting modules in generality, and an abstract construction which, with additional work, can be made explicit in examples
(as they do for the Temperley-Lieb algebra). The first half of our paper has a different goal: the explicit construction for the special case of $\sl_n$-webs. For what it is worth, note
that the results of \cite{AST} do not apply directly to $\sl_n$-webs, only to the representation category $\Fund$, and they are appropriately careful in \cite[\S 6.1.4]{AST} when they
state the connection with $\sl_n$-webs vaguely. Our results here prove that these integral forms are equivalent categories, but this was not known previously.

Most of the cellular categories in the literature have an additional feature: they are monoidal. Decomposing a tensor product of irreducible modules into irreducibles in a semisimple
category is known as the art of \emph{plethyism}. Plethyism for fundamental representations in semisimple $\sl_n$-representation theory is particularly easy. Let $\w_a$ be the $a$-th
fundamental weight, for $1 \le a \le n-1$, and $\Om(a)$ be the set of weights with nonzero weight spaces in $V_{\w_a}$. Then for an arbitrary dominant weight $\l$, $V_\l \ot V_{\w_a} \cong
\oplus V_{\l + \mu}$ where the direct sum is over $\mu \in \Om(a)$ such that $\l + \mu$ is dominant. A sequence of dominant weights, each differing from its neighbors by an element of
$\Om(a)$ for some $a$, is called a \emph{miniscule Littelman path}; if the sequence of $a$'s is named $\un{w}$, and the final weight in the sequence is $\l$, then this is an element of the
set we call $E(\un{w}, \l)$. In particular, the tensor product $V_{\un{w}}$ decomposes into simples, with one copy of $V_\l$ for each element of $E(\un{w}, \l)$.

A miniscule Littelman path is an inductive decomposition of the tensor product $V_{\un{w}}$: the $d$-th weight in the path is the summand chosen in the tensor product of the first $d$
fundamental representations in $\un{w}$. In particular, in a monoidal category, one can always perform plethyism inductively. This has implications for the form of the morphisms which we
call light ladders. They should have an inductive, tiered construction, as in the following schematic diagram.

\begin{equation} \label{LLschematic} {
\labellist
\small\hair 2pt
 \pinlabel {$L_d$} [ ] at 43 18
 \pinlabel {$E_{\mu}$} [ ] at 100 73
 \pinlabel {$L_{d+1}$} [ ] at 184 49
 \pinlabel {$=$} [ ] at 116 51
 \pinlabel {$a$} [ ] at 104 8
 \pinlabel {$\un{w}$} [ ] at -10 3
 \pinlabel {$\un{x}_{\l}$} [ ] at 12 34
 \pinlabel {$\un{x}_{\l + \mu}$} [ ] at 12 114
 \pinlabel {$\un{w}$} [ ] at 180 23
 \pinlabel {$\un{x}_{\l + \mu}$} [ ] at 178 76
\endlabellist
\centering
\ig{1}{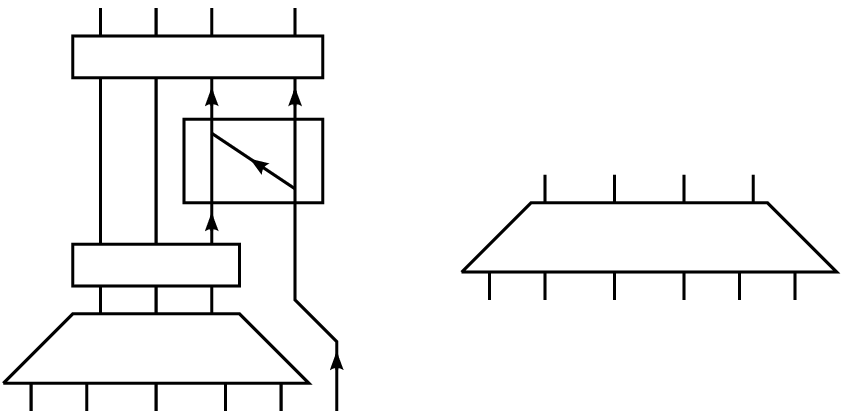}
} \end{equation}

This diagram should be read from bottom to top, as for all diagrams in this paper. A trapezoid represents a light ladder map, which is a cellular ``half-basis'' element. Let us arbitrarily
choose a sequence of fundamental weights $\un{x}_\l$ for each dominant weight $\l$, such that the sum of the weights is $\l$. Suppose that $\un{w}$ is a sequence of length $d$, and we have
already constructed the light ladder web $L_d$ from $\un{w}$ to $\un{x}_\l$ corresponding to a miniscule Littelman path in $E(\un{w}, \l)$. Now we want to construct the light ladder
$L_{d+1}$ for the sequence $\un{w}a$, associated to a choice of weight $\mu \in \Om(a)$ with $\l + \mu$ dominant. There is a particular web $E_{\mu}$, the \emph{elementary light ladder},
which when tensored with an identity map, gives the projection from $V_\l \ot V_{\w_a}$ to $V_{\l + \mu}$ modulo lower terms. The web $E_{\mu}$ has inputs and outputs that are actually
independent of $\l$, only depending on $\mu$. In order to put those inputs in the correct place, we first apply a \emph{neutral ladder}, symbolized above by a rectangle, which is equal
(modulo lower terms) to the braiding isomorphism that reorders the tensor factors. If it is not possible to rearrange the sequence $\un{x}_\l$ so that it ends with the inputs of $E_\mu$,
then in fact $\l + \mu$ was not dominant to begin with. Finally, after applying $E_{\mu}$, we rearrange the outputs again to obtain $\un{x}_{\l + \mu}$. Thus we have constructed $L_{d+1}$,
the next light ladder in the inductive algorithm.

Our actual cellular basis, the \emph{double ladders basis}, just puts two halves with the same middle together.
\begin{equation} \label{doubleintro} {
\labellist
\small\hair 2pt
 \pinlabel {$LL_{\eb}$} [ ] at 44 17
 \pinlabel {$\overline{LL}_{\fb}$} [ ] at 44 46
 \pinlabel {$\un{w}$} [ ] at 91 0
 \pinlabel {$\un{x}_\l$} [ ] at 75 32
 \pinlabel {$\un{y}$} [ ] at 91 67
\endlabellist
\centering
\ig{1}{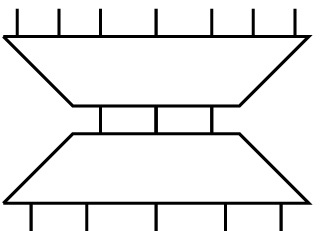}
} \end{equation}
Here, $\eb \in E(\un{w}, \l)$ and $\fb \in E(\un{y}, \l)$. The overline indicates the duality involution.

\begin{remark} The reader familiar with Soergel diagrammatics will note that this algorithm for producing a basis is identical in almost all regards to the construction of the double
leaves basis from light leaves in \cite{EWGr4sb}. This was, of course, the inspiration for light ladders and double ladders. The light leaves basis itself dates back to work of Libedinsky
\cite{LibLL}. More will be said about the connections between Soergel diagrammatics and $\sl_n$-webs below. \end{remark}

That a semisimple monoidal category should have a basis which is constructible by such an algorithm is almost a tautological unwinding of plethyism and the observations about factoring of
moprhisms through irreducible summands. We expect many other settings (such as Khovanov and Lauda's category $\UC$ \cite{LauSL2, KhoLau10}) to have bases constructed by similar
algorithms. The interesting part of the story is the explicit computation of the elementary light ladders $E_{\mu}$ (and the fact that they are independent of $\l$). We define them in
\eqref{elementarylightladder}.

That an integral form (or a deformation) of a semisimple monoidal category should have such a basis is also not just a happy coincidence. We suspect that the arguments of \cite{AST} can be
adapted to prove that such bases are quite general. Nonetheless, it is much more difficult to prove (Theorem \ref{thm:doubleladdersspan}) that the webs produced by this algorithm are, in
fact, a basis. The proof of linear independence in \S\ref{subsec-lightladdersindep} involves evaluating these webs after applying a functor, in this case the functor $\G$ from $\sl_n$-webs
to $U_q$-modules, defined in Cautis-Kamnitzer-Morrison \cite{CKM}. The proof that they span involves three components. \begin{enumerate} \item A topological, Morse-theory like argument
which rigidifies the kind of diagrams one needs to analyze. This was accomplished already in \cite{CKM} in their proof that ladders, certain special $\sl_n$-webs, span the category. \item
An argument specific to the setting, which states that elementary light ladders span a certain class of morphisms modulo certain transformations (see \S\ref{subsec-lightladdersspan1}).
\item A general argument which applies to any basis constructed by an analogous tiered algorithm, which bootstraps the previous result into a proof of spanning (see
\S\ref{subsec-lightladdersspan2}). \end{enumerate} We have tried to phrase our arguments in a language which makes it clear how they adapt to other settings.

\begin{remark} The proof that the diagrammatic Soergel category is spanned by double leaves, found in the last chapter of \cite{EWGr4sb}, follows this same rubrick. At the time, we (myself
and Williamson) were somewhat dissatisfied with the proof, as it seemed overly complicated and ad hoc. Having had to repeat the proof in a different context, it now seems far more pleasing
and philosophically correct (at least to me). \end{remark}

Once one knows that double ladders form a basis for $\sl_n$-webs over $\ZM[q,q^{-1}]$, it is not difficult to prove that the functor $\G$ is an isomorphism from this integral form of
$\sl_n$-webs to the category $\Fund$. This is done in Theorem \ref{EquivThm}, following an analogous proof (due to Ben Webster) from the appendix to \cite{ELib}.

\subsection{Intersection forms and clasps}
\label{subsec-ifclasp}

For a sequence $\un{w}$ of fundamental weights, if the sum of the fundamental weights is the dominant integral weight $\l$, then we say that $\un{w}$ \emph{expresses} $\l$, and we write
$\un{w} \in P(\l)$. The objects in $P(\l)$ are all isomorphic, as $\sl_n$-webs are a braided monoidal category. After base change to $\QM(q)$, they have an irreducible top summand $V_\l$.

For $\un{w}, \un{x} \in P(\l)$, the projection map $\un{w} \to V_\l$ (resp. the inclusion map $V_\l \to \un{x}$) can not be described using $\sl_n$-webs, because $V_\l$ is not an object in
$\Fund$. However, the composition $\un{w} \to V_\l \to \un{x}$, which projects to this common summand, can be described using $\sl_n$-webs. These projections will be packaged as a clasp.

\begin{defn} Let $\Bbbk$ be an extension of $\ZM[q,q^{-1}]$. The \emph{$\l$-clasp} $\phi^\l$ associated to a dominant weight $\l$ is a family of maps $\{ \phi_{\un{w}, \un{x}} \}$
associated to pairs $\un{w}, \un{x} \in P(\l)$, which are $\Bbbk$-linear combinations of webs. They should satisfy the following properties: \begin{itemize} \item \emph{Compatibility:}
$\phi_{\un{x}, \un{y}} \circ \phi_{\un{w}, \un{x}} = \phi_{\un{w}, \un{y}}$ for any three elements of $P(\l)$. \item \emph{Orthogonality:} $\phi_{\un{w}, \un{x}} \circ a = 0$ and $b \circ
\phi_{\un{w}, \un{x}} = 0$ for any $a, b \in I_{< \l}$. Here, the \emph{ideal of lower terms} $I_{< \l}$ is the ideal of morphisms which factor through shorter tensor products. \item
\emph{Unitality:} $\phi_{\un{w}, \un{w}} \equiv 1$ modulo $I_{< \l}$. \end{itemize} \end{defn}

The uniqueness of the clasp is a straightforward argument following from the cellular structure. It is guaranteed to exist when $\Bbbk = \QM(q)$, because the projection to $V_\l$ satisfies
these three properties. One of the major questions is: what is the smallest extension $\Bbbk$ for which the $\l$-clasp is defined? In other words, what are the denominators in the formula
for the clasp?

Compatibility implies that, for each $\un{w} \in P(\l)$, the map $\phi_{\un{w}, \un{w}}$ is an idempotent, which is what we had called the clasp earlier in this introduction, and which the
literature typically refers to as a clasp. Instead, we use the word clasp to denote the entire family $\phi$ or any morphisms therein, not necessarily with the same source and target. We
would like to advocate the entire family $\phi$ as a more natural object of study. An idempotent is good enough to pin down a summand inside an object, while a family $\phi$ as above is
required to pin down a common summand inside a family of objects. As shall be seen, it is also a more practical object to compute, because of the nature of recursive formulas.

Finding a closed formula for the clasp seems currently out of reach. Even for the Temperley-Lieb algebra, this was quite a chore (see \cite{Morr02}). However, we provide a conjecture
which gives a recursive formula to compute all clasps.

\begin{remark} For any $\Bbbk$, and any $\un{w} \in P(\l)$, there is a unique summand $T_\l$ of $\un{w}$ which does not occur inside any shorter tensor products (and this summand is
common to all $\un{w} \in P(\l)$). This can be shown in general for an SOACC. Speaking from a representation theory perspective, $T_\l$ is the indecomposable tilting module in $\Rep \ot
\Bbbk$. However, $T_\l$ will typically be much bigger than $V_\l$, and the family of projection maps $\un{w} \to T_\l \to \un{x}$ will not satisfy orthogonality. We might call this
family a \emph{$\Bbbk$-clasp}. Its computation is beyond the scope of this paper, and will not follow from our conjecture, but will follow from more general computations along the same
lines. \end{remark}

\begin{remark} The geometric Satake equivalence, or a ($q$-deformed) algebraic version developed in \cite{EQuantumI}, will take a tensor product of fundamental representations and return a
singular Soergel bimodule associated to the affine Weyl group, as in \cite{WillSingular}. The $\sl_n$-webs are transformed into morphisms between bimodules, and clasps become projections to
indecomposable Soergel bimodules. Thus, finding clasps (and $\Bbbk$-clasps) will also solve problems related to Soergel bimodules for affine Weyl groups, which was the author's original
motivation.

Along these lines, let us mention that the double ladders basis, after its transformation into bimodules, is part of a larger ``double singular light leaves'' basis, which makes singular
Soergel bimodules in (finite and affine) type $A$ into a fibered cellular category. This will be explored in forthcoming work with Williamson. \end{remark}

In the literature, recursive formulas are described as \emph{single
clasp expansions}, \emph{double clasp expansions}, etcetera, based on the number of clasps which appear in the most complex diagram. For example, here is the familiar single clasp
expansion for $\sl_2$.
\begin{equation} \label{TL1clasp} {
\labellist
\small\hair 2pt
 \pinlabel {$n$} [ ] at 33 39
 \pinlabel {$n+1$} [ ] at 139 39
 \pinlabel {$\displaystyle \sum_{i=1}^{n} \frac{-[i]}{[n+1]}$} at 237 41
 \pinlabel {$i$} [ ] at 332 76
 \pinlabel {$n$} [ ] at 324 23
\endlabellist
\centering
\ig{1}{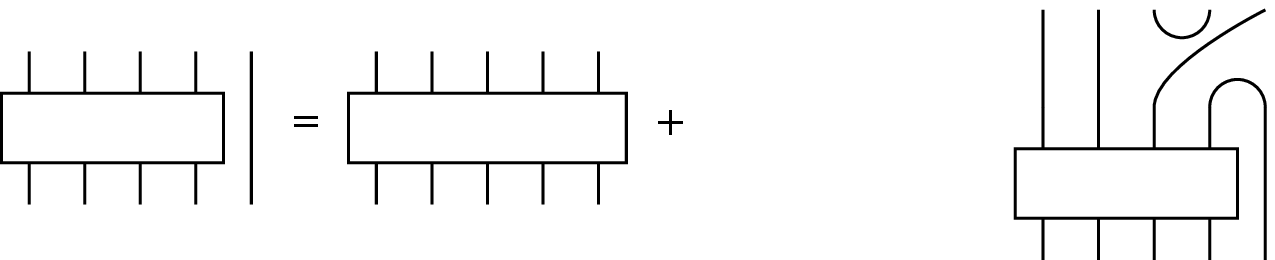}
} \end{equation}
(We use the diagrammatics for the Temperley-Lieb algebra developed by Kauffman \cite{Kauf}. Henceforth, we denote a clasp with a rectangle labeled by the
weight $\l$ it represents, which in this case is the number of strands which enter it; pictured is the case $n=4$. We write $[n]$ for the $n$-th quantum number, $[n] = q^{n-1} + q^{n-3} + \ldots + q^{3-n} + q^{1-n}$. In this formula, $i$ is the first strand which takes part in the cup.) Here is the double clasp expansion.
\begin{equation} \label{TL2clasp} {
\labellist
\small\hair 2pt
 \pinlabel {$n$} [ ] at 33 50
 \pinlabel {$n+1$} [ ] at 140 50
 \pinlabel {$\displaystyle \frac{-[n]}{[n+1]}$} [ ] at 217 50
 \pinlabel {$n$} [ ] at 273 86
 \pinlabel {$n$} [ ] at 273 22
\endlabellist
\centering
\ig{1}{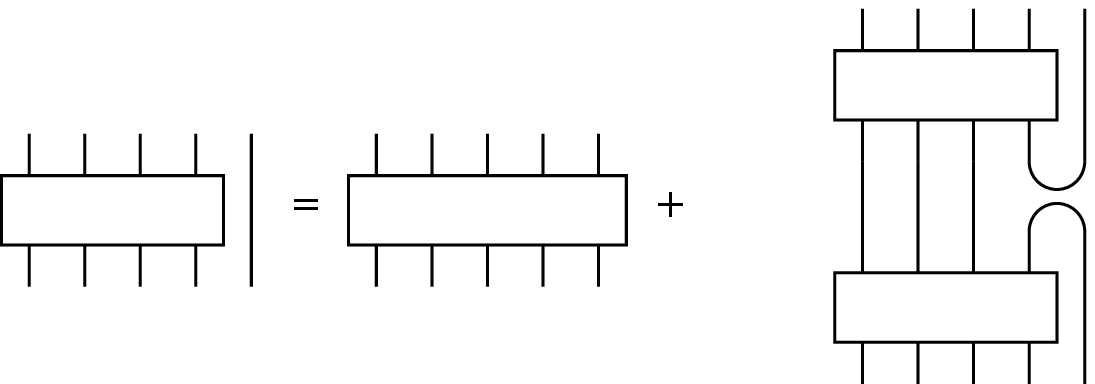}
} \end{equation}
However, we prefer the \emph{triple clasp expansion} below.
\begin{equation} \label{TL3clasp} {
\labellist
\small\hair 2pt
 \pinlabel {$n$} [ ] at 33 50
 \pinlabel {$n+1$} [ ] at 140 50
 \pinlabel {$\displaystyle \frac{-[n]}{[n+1]}$} [ ] at 217 50
 \pinlabel {$n$} [ ] at 273 86
 \pinlabel {$n$} [ ] at 273 22
 \pinlabel {$n-1$} [ ] at 263 54
\endlabellist
\centering
\ig{1}{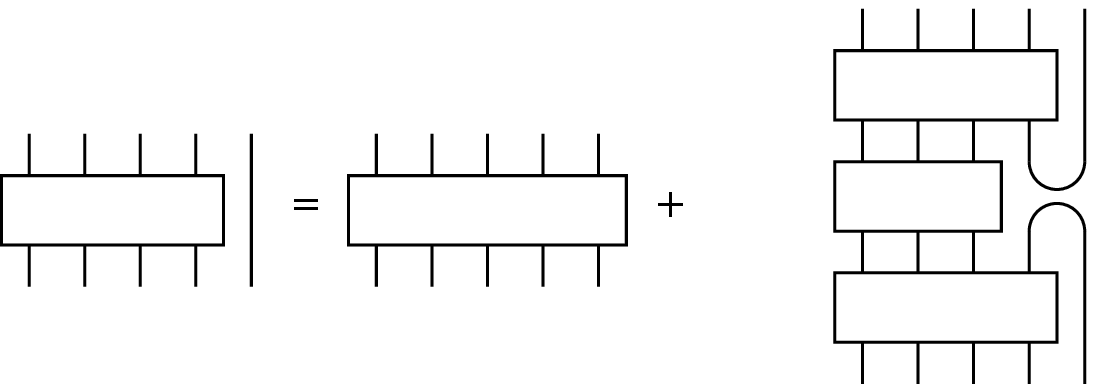}
} \end{equation}

Experts will already recognize that the \emph{central clasp}, which distinguishes \eqref{TL3clasp} from \eqref{TL2clasp}, is redundant. This clasp is a linear combination of diagrams, all of which vanish (by orthogonality) when plugged in, except the identity diagram with coefficient $1$ (by unitality).  However, an analogous central clasp is not always redundant in the expansion formulas for $\sl_n$, $n > 2$.

The central clasp plays an important philosophical role. We know in the representation theory of $\sl_2$ that $V_n \ot V_1 \cong V_{n+1} \oplus V_{n-1}$. The formula \eqref{TL3clasp}
explicitly takes the identity morphism of $V_n \ot V_1$, decomposes it as a sum of two orthogonal idempotents, and further still, it factors the second idempotent as projection $V_n \ot V_1
\to V_{n-1}$ followed by an inclusion $V_{n-1} \to V_n \ot V_1$. The triple clasp expansion is the morphism-theoretic equivalent to the object-theoretic statement $V_n \ot V_1 \cong V_{n+1}
\oplus V_{n-1}$. A priori, the expansion \eqref{TL2clasp} only factors the second idempotent through $V_1^{\ot (n-1)}$, but it is not immediately clear that it factors further through
$V_{n-1}$. In the single clasp expansion, it is even less clear what the individual terms of the sum mean. The single clasp expansion may be most useful in some computer applications, however.

As mentioned before, Kim \cite{Kim07} has computed some recursion formulas for clasps in rank 2, such as a single and a double clasp expansion in type $A_2$ and $B_2$, and a quadruple
clasp expansion in type $A_2$. Outside of type $A_1$ and the work of Kim, it appears that no other recursive formulas are to be found in the literature. The work of Cautis \cite{CautisClasp} gives another possible avenue to computing clasps, although its implications for the Grothendieck group are not clear.

This factorization of idempotents gives us an immediate interpretation of the coefficient $\frac{-[n]}{[n+1]}$ which appears in \eqref{TL3clasp}; it is the reciprocal of (the value of) a
\emph{local intersection form}. There is a pairing $\Hom(V_n \ot V_1, V_{n-1}) \times \Hom(V_{n-1}, V_n \ot V_1) \to \End(V_{n-1}) = \QM(q)$ induced by composition; after using duality to
identify these two $\Hom$ spaces, this becomes a form on $\Hom(V_{n-1}, V_n \ot V_1)$ called the \emph{local intersection form} or \emph{LIF}. The space $\Hom(V_{n-1}, V_n \ot V_1)$ is one-dimensional; this
crucial fact comes from the fact that the weight space $-1$ in $V_1$ is one-dimensional. The rank of the LIF determines how many copies of $V_{n-1}$ appear as summands inside $V_n \ot
V_1$. The following computation shows that the LIF corresponds to the $1 \times 1$ matrix with entry $\frac{-[n+1]}{[n]}$.
\begin{equation} \label{TLlif} {
\labellist
\small\hair 2pt
 \pinlabel {$n-1$} [ ] at 23 22
 \pinlabel {$n-1$} [ ] at 24 86
 \pinlabel {$n-1$} [ ] at 166 54
 \pinlabel {$n$} [ ] at 31 54
 \pinlabel {$\displaystyle \frac{-[n+1]}{[n]}$} [ ] at 114 54
\endlabellist
\centering
\ig{1}{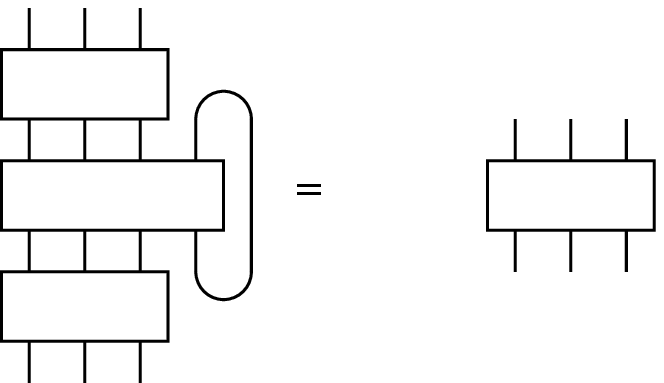}
} \end{equation}
The LHS of \eqref{TLlif} can be thought of as the idempotent ``turned inside-out." This computation can be performed inductively: one takes the middle clasp on the LHS of \eqref{TLlif}, and replaces it using
\eqref{TL3clasp} (but for $n-1$ instead of $n$). Thus there are two formulas at work: the formula for a clasp \eqref{TL3clasp} which involves the reciprocal of a LIF, and the formula for the
LIF \eqref{TLlif} which involves a shorter clasp. These can be solved together by recursion; in this example, it is a good exercise to deduce that the LIF is $\frac{-[n+1]}{[n]}$ (assuming
the case $n=1$, which is one of the relations in the Temperley-Lieb algebra).

\begin{remark} In a specialization $\Bbbk$ where $[n+1]=0$ (and no smaller quantum numbers vanish), the LIF has rank $0$, and $V_{n-1}$ is not a summand of $V_n \ot V_1$. Instead, this
tensor product will be indecomposable, and will be the top summand $T_{n+1}$. This will ruin the inductive procedure, which is why computing $\Bbbk$-clasps is difficult. \end{remark}

To give another example, let us demonstrate the triple clasp expansions for $\sl_3$ (new in this paper), which we draw using the $\sl_3$-webs of Kuperberg \cite{Kupe}. We encourage the reader to contrast these formulas with the single-clasp expansions from \cite{Kim07}.
\begin{equation} \label{sl3claspU} {
\labellist
\small\hair 2pt
 \pinlabel {$m,n$} [ ] at 24 70
 \pinlabel {$m+1,n$} [ ] at 121 70
 \pinlabel {$\displaystyle \frac{[m]}{[m+1]}$} [ ] at 189 70
 \pinlabel {$m,n$} [ ] at 236 117
 \pinlabel {$m,n$} [ ] at 236 22
 \pinlabel {\tiny{$m-1, n+1$}} [ ] at 241 70
 \pinlabel {$\frac{[n][m+n+1]}{[n+1][m+n+2]}$} [ ] at 317 69
 \pinlabel {$m,n$} [ ] at 379 106
 \pinlabel {$m,n$} [ ] at 380 42
 \pinlabel {\tiny{$m,n-1$}} [ ] at 372 73
\endlabellist
\centering
\ig{1}{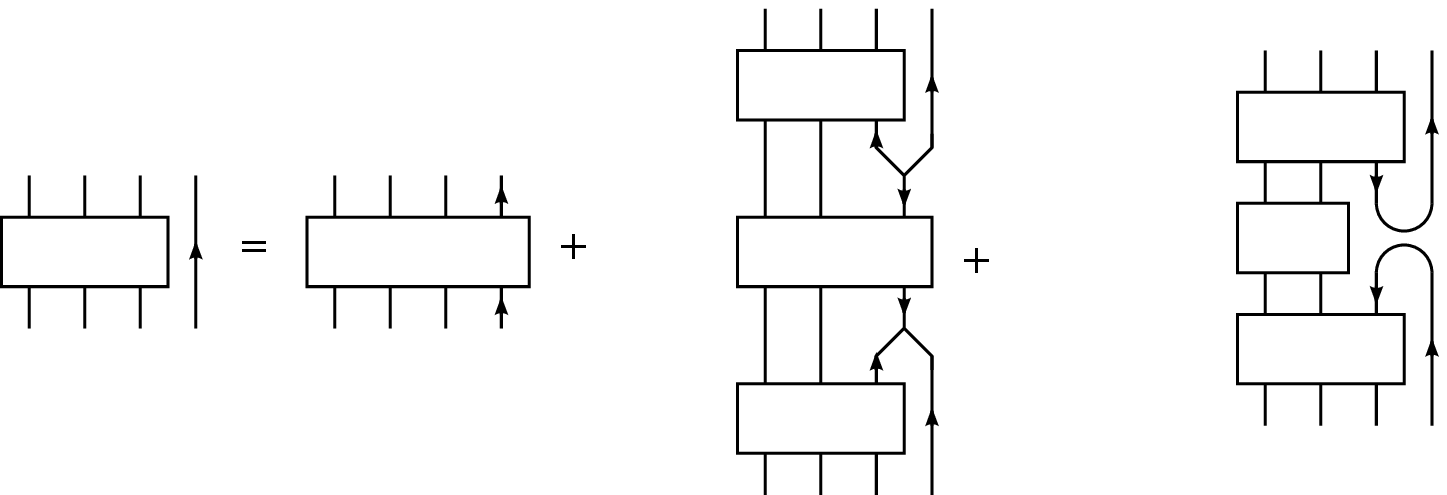}
} \end{equation}
\begin{equation} \label{sl3claspD} {
\labellist
\small\hair 2pt
 \pinlabel {$m,n$} [ ] at 24 70
 \pinlabel {$m,n+1$} [ ] at 121 70
 \pinlabel {$\displaystyle \frac{[n]}{[n+1]}$} [ ] at 189 70
 \pinlabel {$m,n$} [ ] at 236 117
 \pinlabel {$m,n$} [ ] at 236 22
 \pinlabel {\tiny{$m+1, n-1$}} [ ] at 241 70
 \pinlabel {$\frac{[m][m+n+1]}{[m+1][m+n+2]}$} [ ] at 317 69
 \pinlabel {$m,n$} [ ] at 379 106
 \pinlabel {$m,n$} [ ] at 380 42
 \pinlabel {\tiny{$m-1,n$}} [ ] at 372 73
\endlabellist
\centering
\ig{1}{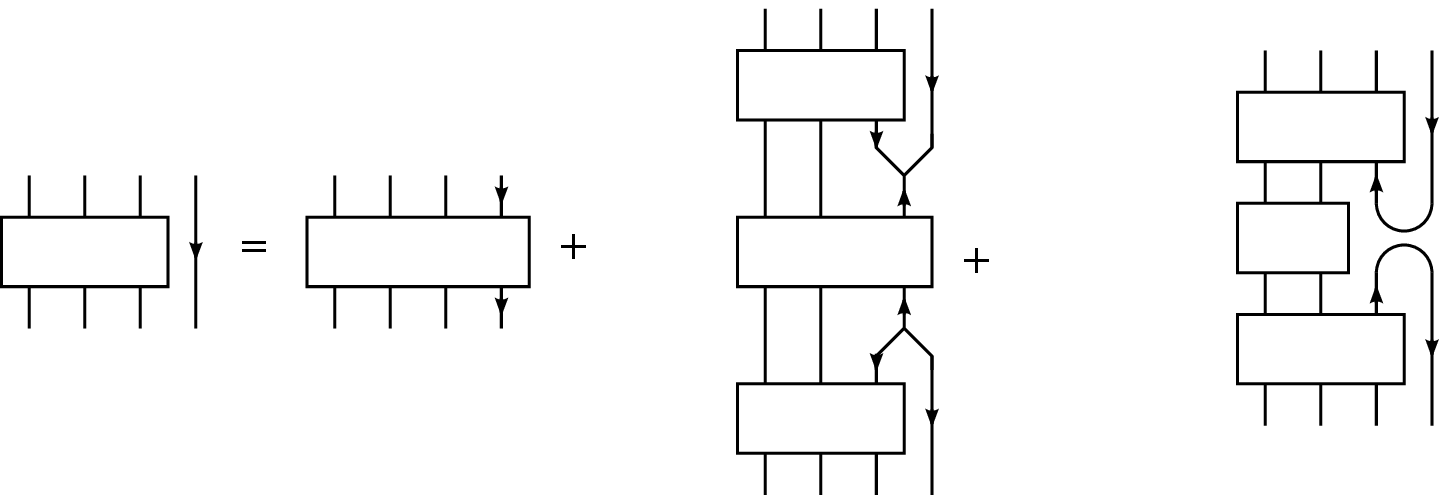}
} \end{equation}
For Kuperberg's $\sl_3$-webs, all strands should be oriented: an upward orientation represents $\w_1$, and a downward orientation $\w_2$. The rectangle labeled $(m,n)$ has $m$ upward and $n$ downward inputs (i.e. strands coming in from below), and the same number of outputs (although possibly in a different order). We have not placed orientations on certain strands because those orientations are irrelevant, so long as the appropriate number enter and leave each clasp. Note, for instance, that the last diagram in \eqref{sl3claspU} is not consistent when $n=0$ because no downward strand enters the $(m,n)$ clasp. However, the coefficient of this diagram is zero when $n=0$, so there is no problem. In fact, this happens precisely when the object being factored through, $(m, n-1)$, is not actually a dominant weight.

The triple clasp expansion \eqref{sl3claspU} is really the decomposition of the identity map of $$V_{m \w_1 + n \w_2} \ot V_{\w_1}$$ into orthogonal idempotents which factor, respectively,
though $$V_{(m+1) \w_1 + n \w_2},\;\; V_{(m-1) \w_1 + (n+1) \w_2},\;\; \textrm{ and }\;\; V_{m \w_1 + (n-1) \w_2}.$$ Once again, the morphism spaces from each summand into the tensor
product are one-dimensional; this is because the fundamental representation $V_{\w_1}$ is miniscule. The coefficients appearing in \eqref{sl3claspU} are the reciprocals of the values of $1
\times 1$ local intersection forms. To compute these local intersection forms, we turn each idempotent inside-out, and resolve the middle clasp using the same triple clasp expansions, but
for a smaller weight. This methodology produces a recursive formula for the coefficients. Finding this recursive formula is an exercise in $\sl_3$-webs, and solving it is an exercise in
quantum numbers.

As the reader may have noticed, the same general plethyism formalities which led to the double ladder algorithm also lead to the triple clasp expansion. Thus we claim that, for any $\sl_n$, clasps (which are drawn as ovals) satisfy the following recursion.

\begin{equation} \label{tripleintro}{
\labellist
\small\hair 2pt
 \pinlabel {$\l+\w_a$} [ ] at 66 118
 \pinlabel {$\l$} [ ] at 260 118
 \pinlabel {$\l$} [ ] at 466 25
 \pinlabel {$\l+\mu$} [ ] at 482 118
 \pinlabel {$\l$} [ ] at 466 210
 \pinlabel {$\displaystyle - \sum_{\mu} \frac{1}{\k_{\l, \mu}}$} [ ] at 370 118
 \pinlabel {\tiny $a$} [ ] at 138 101
 \pinlabel {\tiny $a$} [ ] at 138 138
 \pinlabel {\tiny $a$} [ ] at 337 113
 \pinlabel {\tiny $a$} [ ] at 547 9
 \pinlabel {\tiny $a$} [ ] at 547 214
 \pinlabel {$E_{\mu}$} [ ] at 560 75
 \pinlabel {$\overline{E}_{\mu}$} [ ] at 560 167
\endlabellist
\centering
\ig{.73}{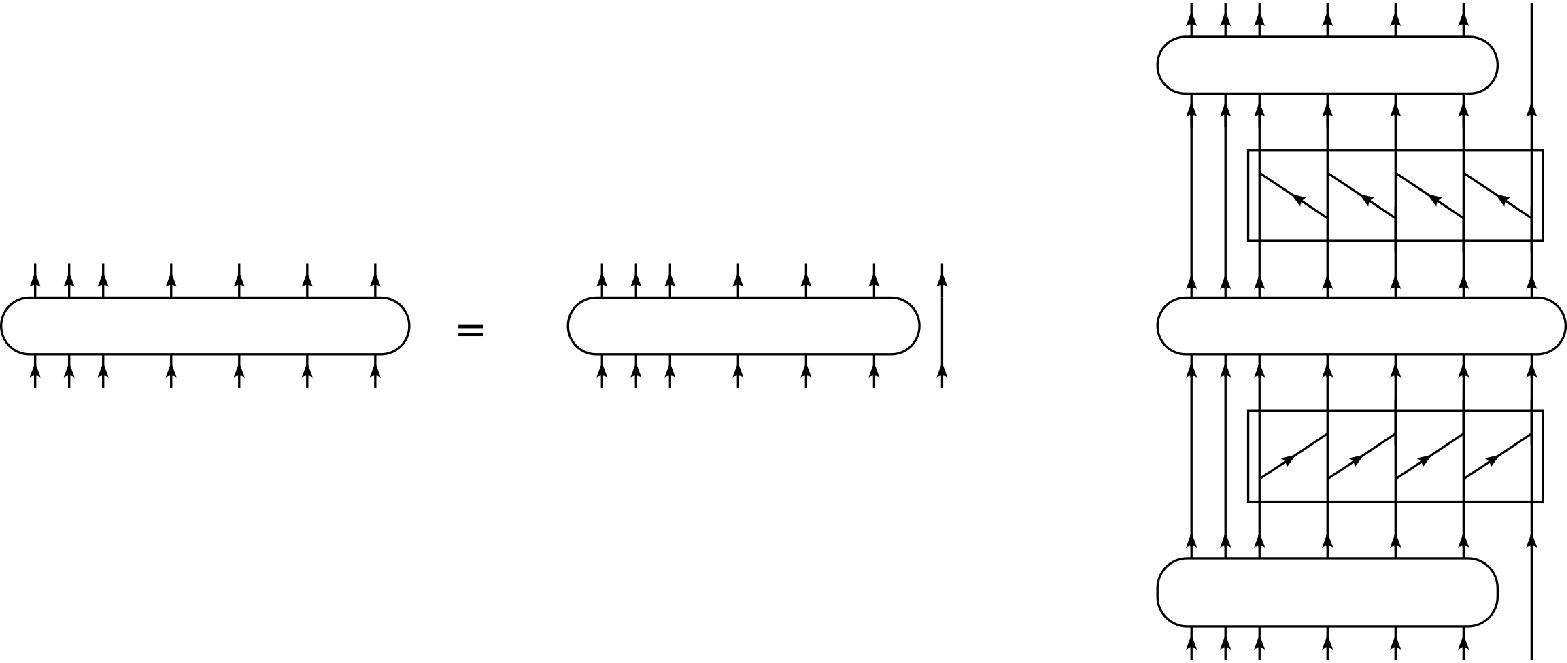}
} \end{equation} 
The coefficients $\k_{\l, \mu}$, which are the \emph{local intersection forms of $V_\l \ot V_{\w_a}$ at $V_{\l + \mu}$} for $\mu \in \Om(a)$ with $\l + \mu$ dominant, satisfy
\begin{equation} \label{intersectionformintro}{
\labellist
\tiny\hair 2pt
 \pinlabel {\small $\l + \mu$} [ ] at 71 22
 \pinlabel {\small $\l + \mu$} [ ] at 71 212
 \pinlabel {\small $\l$} [ ] at 57 118
 \pinlabel {\small $\l + \mu$} [ ] at 310 118
 \pinlabel {\small $\k_{\l,\mu}$} [ ] at 208 118
 \pinlabel {\small $\overline{E}_{\mu}$} [ ] at 150 69
 \pinlabel {\small $E_{\mu}$} [ ] at 150 166
\endlabellist
\centering
\ig{.8}{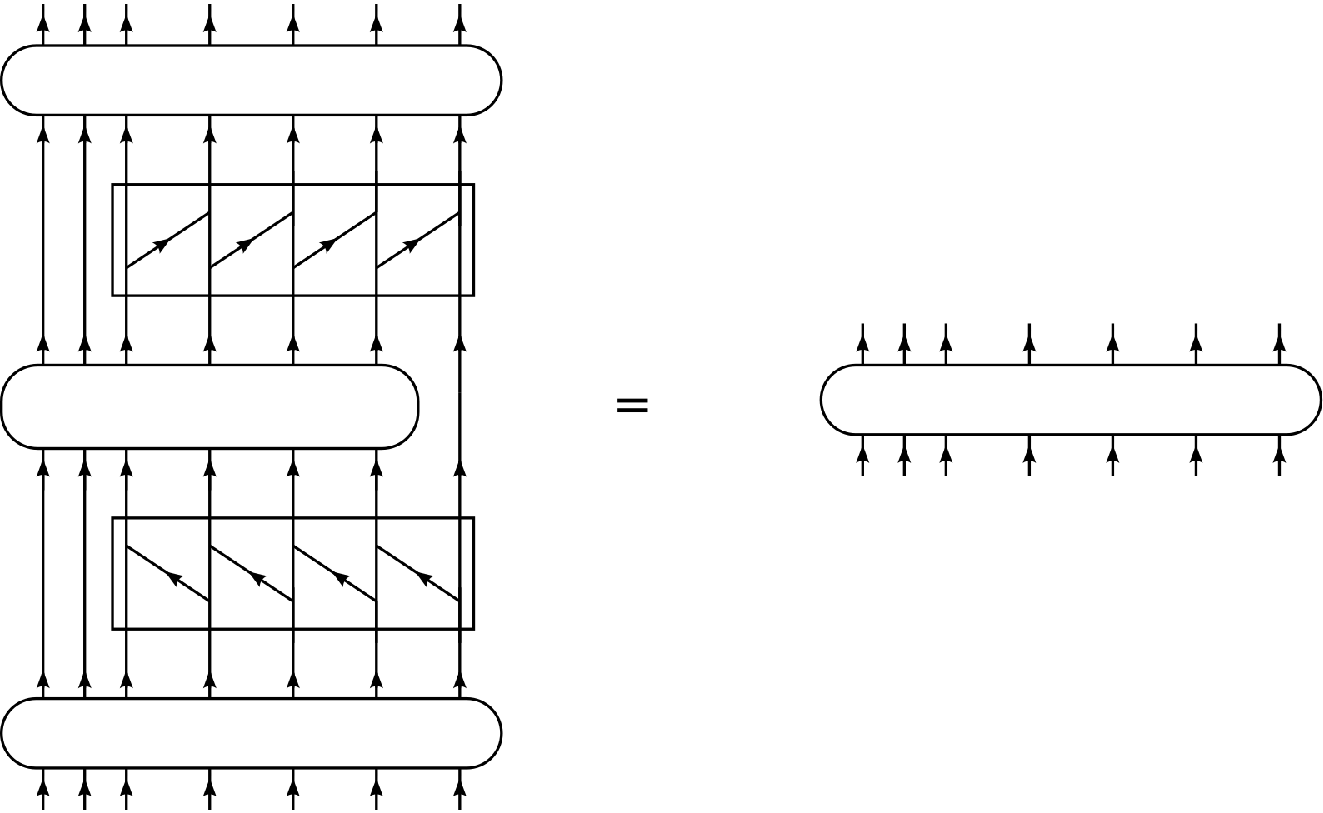}
}. \end{equation}
We prove this recursion in \S\ref{subsec-tripleclasp}, but it is mostly a tautology. The interesting question is the computation of $\k_{\l, \mu}$, which is done by using \eqref{tripleintro} to replace the $\l$-clasp in \eqref{intersectionformintro} and resolving. But resolving is difficult, and therein lies the rub.

\begin{remark} The unlabeled strands in these diagrams are arbitrary, so long as the fundamental weights add up to $\l$ or the appropriate weight. Thus this formula describes the clasp for
$\l+\w_a$ as a morphism between any two sequences which both end in $a$. To describe the clasp as a morphism between sequences in $P(\l + \w_a)$ which do not end in $a$, one can use
neutral ladders to reorder the sequence. One can prove that a neutral ladder pre- or post-composed with a clasp is also a clasp. (Similarly, to give a clasp between objects that involve
the duals of fundamental representations, use the tag isomorphism to swap the orientation first. Tag morphisms send a clasp to a clasp.) \end{remark}

In this paper we compute $\k_{\l, \mu}$ for $\sl_3$ and $\sl_4$ (and the same techniques compute many other examples in $\sl_n$ for $n > 4$). They can be found in \S\ref{subsec-conjecture}. Having computed them, a pattern arose, which led to the following conjecture.

\begin{defn} Fix $n \ge 2$. Let $\rho$ denote the half-sum of the positive roots, which is also the sum of all fundamental weights. For any dominant weight $\l$ and
positive root $\a$ let $A(\l, \a) = \langle \l + \rho, \a \rangle$. \end{defn}

\begin{defn} For any weight $\mu \in \Om(a)$ in a fundamental representation $V_{\w_a}$, let $w_{\mu}$ be the unique element of $S_n$ which sends $\w_a$ to $\mu$, and has minimal length
(i.e. is a minimal coset representative for $S_a \times S_{n-a}$, the stabilizer of $\w_a$). Let $\Phi(\mu)$ denote the subset of positive roots $\a$ for which $w_{\mu}^{-1}(\a)$ is a
negative root. \end{defn}

\begin{claim} For any dominant weight $\l$, $A(\l+\mu, \a) = A(\l,\a)-1$ if and only if $\a \in \Phi(\mu)$, and $A(\l,\a)-1 = 0$ for some $\a \in \Phi(\mu)$ if and only if $\l + \mu$ is
not dominant. \end{claim}

\begin{conj} \label{MainConjIntro} Let $\l$ be dominant, and $\mu \in \Om(a)$. We conjecture that, whenever $\l + \mu$ is dominant, one has
\begin{equation} \k_{\l, \mu} = \prod_{\a \in \Phi(\mu)} \frac{[A(\l, \a)]}{[A(\l, \a)-1]}, \end{equation}
or equivalently,
\begin{equation} \k_{\l, \mu} = \prod_{\a \in \Phi(\mu)} \frac{[A(\l, \a)]}{[A(\l+\mu, \a)]}. \end{equation}
These products of ratios of quantum numbers are well-defined and equal by the above claim. \end{conj}

We have outlined a plan of attack for a direct, computational proof in \S\ref{subsec-recursive1} and \S\ref{subsec-recursive2}. This outline is sufficient for $n=3, 4$. At the moment, we
do not have a good philosophical explanation for this conjecture. One also hopes that more sophisticated techniques, such as skew Howe duality, will provide an easier proof.

\subsection{Complements}
\label{subsec-complements}

We record here some additional remarks.

Let $X$ be an indecomposable module in a $\Bbbk$-linear additive category, whose endomorphism ring is $\Bbbk$. Let $B$ be an arbitrary object. Then the composition map \[ \Hom(B,X) \times
\Hom(X,B) \to \End(X) \cong \Bbbk \] gives a pairing called the \emph{local intersection pairing}. In the presence of a duality involution which identifies these two $\Hom$ spaces, one
obtains a form on $\Hom(B,X)$ called the \emph{local intersection form} of $B$ at $X$.

\begin{remark} The term ``local intersection form" is used in analogy to certain forms in the geometry of perverse sheaves. See \cite{EWShadows} or \cite{dCaMig05} for more details. For
SOACCs, this agrees with the familiar cellular pairing. \end{remark}

Let us temporarily assume that scalars are $\QM(q)$, so that we are working in a semisimple category $\Rep'$. Ultimately, one is interested in the strict monoidal category consisting of
all tensor products of simples. The reason that one restricts to tensor products of fundamental representations is to be able to describe the category by generators and relations. However,
if one could compute all the clasps, then one could formally extend the $\sl_n$-web calculus to a calculus for its Karoubi envelope, which would include all tensor products of simples.

Nonetheless, for this description to be interesting, one would need to compute even more relations, such as the quintuple clasp expansion. This should be a formula along the lines of
\begin{equation} \label{quintuple} {
\labellist
\small\hair 2pt
 \pinlabel {$\l$} [ ] at 24 74
 \pinlabel {$\mu$} [ ] at 73 73
 \pinlabel {$\l + \mu$} [ ] at 153 74
 \pinlabel {$\l$} [ ] at 289 22
 \pinlabel {$\mu$} [ ] at 338 21
 \pinlabel {$\l$} [ ] at 290 134
 \pinlabel {$\mu$} [ ] at 336 135
 \pinlabel {$\nu$} [ ] at 309 78
 \pinlabel {$O$} [ ] at 308 48
 \pinlabel {$I$} [ ] at 308 105
 \pinlabel {$\displaystyle \sum \k_{U,D,\nu}$} [ ] at 237 77
\endlabellist
\centering
\ig{1}{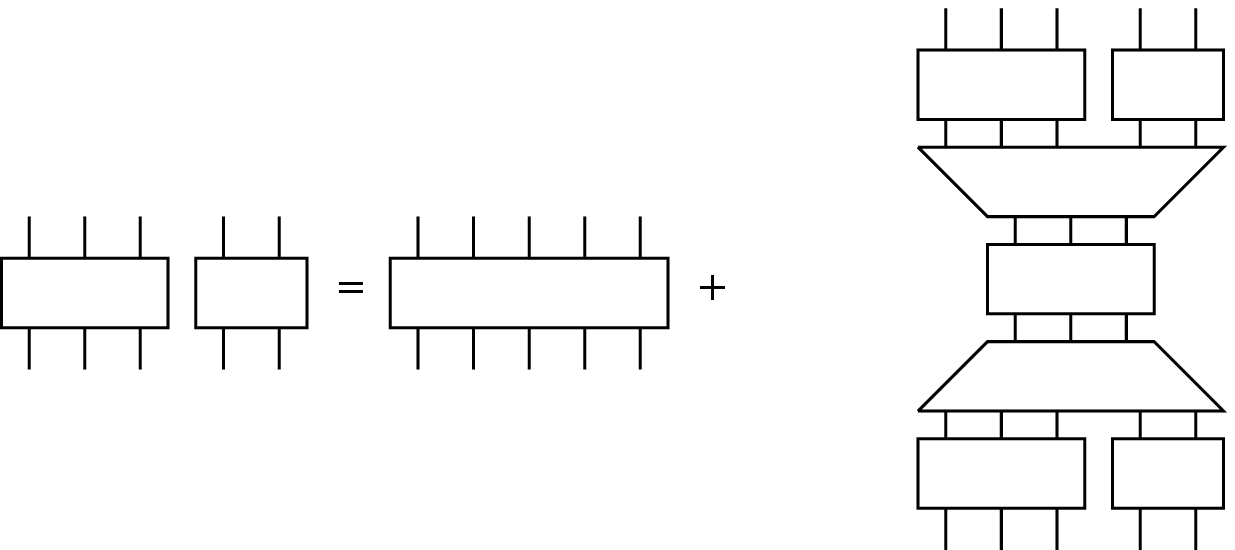}
} \end{equation}
which takes the identity of $V_{\l} \ot V_{\mu}$ and decomposes it into orthogonal idempotents, according to the plethyism of $\sl_n$ representations. The trapezoidal maps $O$ and $I$ are maps similar to light ladders.

A quintuple clasp expansion requires the computation of the local intersection form for $V_{\nu}$ inside $V_{\l} \ot V_{\mu}$, which no longer comes from a one-dimensional morphism space
(the entries of its matrix are $\k_{U,D,\nu}$). One can also seek recursive formulas which add multiple strands at once (like the quintuple clasp expansion, except without applying the
clasp for $V_{\mu}$). Ultimately, these formulas compute local intersection forms which are important in the eventual task of finding $\Bbbk$-clasps for arbitrary specializations.

The form of Conjecture \ref{MainConjIntro} gives hope that there might be formulas for arbitrary intersection forms involving roots, weights, the Weyl group, and other combinatorial data.

Outside of type $A$, plethyism is slightly more complicated, but exactly the same principles apply. The results of this paper should be generalizable to other semisimple Lie algebras.
Conjecture \ref{MainConjIntro} needs modification, as not every weight in a fundamental representation is in the same Weyl group orbit, but it is conceivable that something very similar
should hold.

One significant consequence of Conjecture \ref{MainConjIntro} is that the local intersection forms are positive. That is, they are mutliples of ratios of quantum numbers, with no signs
involved, and they become positive real numbers when $q=1$. Let us make some observations about positivity, for which we assume the reader is somewhat familiar with \cite{EWHodge}.

In \cite{EWHodge}, the author and Williamson prove some Hodge-theoretic statements about local intersection forms in the category of Soergel bimodules for any Coxeter group, which were
used to prove the Soergel conjecture. In particular, we prove that local intersection forms are $\pm$-definite, with an easily computible sign. This positivity result can also be proven
geometrically for Weyl groups and affine Weyl groups. The same geometric proof should suffice for local intersection forms in the category of singular Soergel bimodules \cite{WillSingular}
for Weyl groups and affine Weyl groups. Meanwhile, the geometric Satake equivalence was reinterpreted by the author in \cite{EQuantumI} to give an equivalence between $\sl_n$-webs at $q=1$
and (a subcategory of) singular Soergel bimodules for the affine Weyl group in type $A$. Parity considerations force all the intersection forms to be positive definite, which gives an
explanation for the positive signs at $q=1$.

There is no geometric understanding of what happens when $q \ne 1$, but the algebraic version does deform: the results of \cite{EQuantumI} continue to give an equivalence between
$\sl_n$-webs over $\ZM[q,q^{-1}]$ and (a subcategory of) a $q$-deformation of singular Soergel bimodules in type $A$. We have called this the \emph{quantum algebraic Satake equivalence}.
More precisely, there is a diagrammatic version of singular Soergel bimodules, which agrees with the algebraic version when $q$ is generic, and (has a subcategory which) is equivalent to
$\sl_n$-webs. When $q$ is specialized to a root of unity, additional morphisms appear between Soergel bimodules which are not encoded in the diagrammatic category, and for which a new
diagrammatic calculus is needed (a current research project of the author).

What happens when $q = \z_{2m}$ is a primitive $2m$-th root of unity is extremely interesting. The $q$-deformed reflection representation of the affine Weyl group used to define Soergel
bimodules now factors through the finite complex reflection group $G(m,m,n)$. We think of these Soergel bimodules as ``Soergel bimodules for $G(m,m,n)$," even though Soergel bimodules are
not defined for complex reflection groups in general. Now, the quantum number $[k]$ for $k < m$ is actually a positive real number! Moreover, for any $k$, $[k]$ and $[k+1]$ have the same
sign, or one of them is zero. Thus Conjecture \ref{MainConjIntro} can be used to deduce that the local intersection forms that are still defined are in fact still positive! This should
correspond to an analog of the Soergel conjecture for $G(m,m,n)$. The author believes that Soergel bimodules for $G(m,m,n)$ form yet another example of a beautiful non-geometric category
with extremely geometric properties.

\section{Web calculus}
\label{sec-webs}

\subsection{$\sl_n$-webs}
\label{subsec-webs}

We recall the $\sl_n$ web calculus from \cite{CKM}.

\begin{defn} A \emph{decorated planar graph with boundary} is a kind of graph embedded in the planar strip $\R \times [0,1]$. Edges are labeled by some index set $I$. Edges in the graph
may ``run to the boundary," that is, they may terminate in invisible univalent vertices on the boundary $\R \times \{0,1\}$ of the strip. Such edges are called \emph{boundary edges}. By
reading the indices of the boundary edges on the top boundary $\R \times \{1\}$, one obtains a sequence in $I$ called the \emph{top} or \emph{top sequence} or \emph{output} of the
graph. Similarly, the \emph{bottom} or \emph{input} of the graph encodes the intersection of the graph with the bottom boundary $\R \times \{0\}$. As indicated by the words ``input'' and ``output,'' our diagrams are read from bottom to top.

An \emph{oriented decorated planar graph with boundary} is as above but the edges are oriented. The top and bottom sequences now keep track of the orientations on the boundary edges.
Formally, the top and bottom are sequences in $I^+ \coprod I^-$, where $I^+$ and $I^-$ are copies of $I$. \end{defn}

\begin{defn} \label{websdefn} Fix an integer $n \ge 2$. An \emph{extended $\sl_n$-web} is a kind of oriented decorated planar graph with boundary. The edges are labeled by indices in $\ZM$. The allowed vertices are \emph{trivalent vertices} and \emph{tags}.
	\vskip .15cm
\begin{equation*} {
\labellist
\small\hair 2pt
 \pinlabel {$k$} [ ] at 1 -5
 \pinlabel {$l$} [ ] at 18 -5
 \pinlabel {$k+l$} [ ] at 9 33
 \pinlabel {$k$} [ ] at 45 -5
 \pinlabel {$l$} [ ] at 63 -5
 \pinlabel {$k+l$} [ ] at 54 33
 \pinlabel {$k$} [ ] at 86 -5
 \pinlabel {$n-k$} [ ] at 86 33
 \pinlabel {$k$} [ ] at 118 -5
 \pinlabel {$n-k$} [ ] at 118 33
\endlabellist
\centering
\ig{1.1}{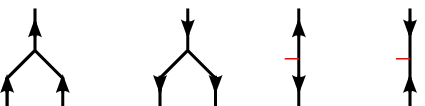}
} \end{equation*}
\vskip .25cm
A \emph{(non-extended) $\sl_n$-web} is an extended $\sl_n$-web where the edges are labeled by indices within the subset $I_n = \{1, \ldots, n-1\}$.
\end{defn}

We will soon define a category $\Fund$ whose morphism spaces are spanned by $\sl_n$-webs. However, relations among $\sl_n$-webs are most easily stated using extended $\sl_n$-webs. We view an
extended $\sl_n$-web as a linear combination of $\sl_n$-webs, following the conventions below (also found in \cite[p6]{CKM}). If an extended $\sl_n$-web has an edge labeled outside of
the set $\hat{I_n} = \{0, 1, \ldots, n-1, n\}$, then it equals zero. For an extended $\sl_n$-web labeled within $\hat{I_n}$, we remove all edges labeled either $0$ or $n$, and replace
trivalent vertices with tags as below.

\vskip .15cm

\begin{equation} \label{remove0n} {
\labellist
\small\hair 2pt
 \pinlabel {$k$} [ ] at 2 52
 \pinlabel {$0$} [ ] at 18 52
 \pinlabel {$k$} [ ] at 10 90
 \pinlabel {$k$} [ ] at 46 52
 \pinlabel {$k$} [ ] at 46 90
 \pinlabel {$k$} [ ] at 78 52
 \pinlabel {$n-k$} [ ] at 95 52
 \pinlabel {$n$} [ ] at 86 90
 \pinlabel {$k$} [ ] at 117 52
 \pinlabel {$n-k$} [ ] at 134 52
 \pinlabel {$k$} [ ] at 2 -5
 \pinlabel {$0$} [ ] at 18 -5
 \pinlabel {$k$} [ ] at 10 35
 \pinlabel {$k$} [ ] at 46 -5
 \pinlabel {$k$} [ ] at 46 35
 \pinlabel {$k$} [ ] at 78 -5
 \pinlabel {$n-k$} [ ] at 95 -5
 \pinlabel {$n$} [ ] at 86 35
 \pinlabel {$k$} [ ] at 117 -5
 \pinlabel {$n-k$} [ ] at 134 -5
\endlabellist
\centering
\ig{1.3}{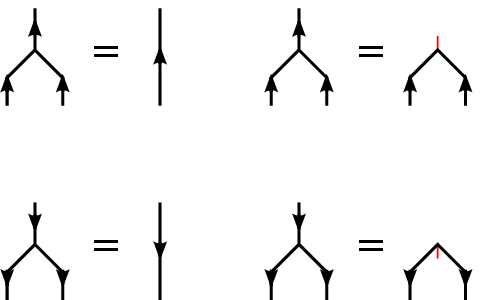}
} \end{equation}

\vskip .25cm

Now, let $\d$ be an indeterminate, and let $\Bbbk$ be an algebra over $\ZM[\d]$. The typical choice is that $\Bbbk = \ZM[q,q^{-1}]$ or $\Bbbk = \QM(q)$, and $\d = q+q^{-1}$. When $\Bbbk
= \ZM[q,q^{-1}]$ it contains elements known as quantum numbers and quantum binomial coefficients, which all happen to be polynomials in $\d = q+q^{-1}$. We now define quantum
numbers and quantum binomial coefficients in an arbitrary $\ZM[\d]$-algebra to be the images of the corresponding polynomials in $\d$.

\begin{defn} Let $\Bbbk$ be an algebra over $\ZM[\d]$. Then set $[0]=0$, $[1]=1$, and define $[n]$ inductively for any $n \in \ZM$ by \[ [n] \d = [n+1] + [n-1]. \] This implies that
$[-n] = -[n]$. For $n \in \ZM$ and $k \in \ZM_{\ge 0}$ let $\qbinom{n}{k}$ be defined as \[ \qbinom{n}{k} = \frac{[n] \cdot [n-1] \cdots [n-k+1]}{[k] \cdot [k-1] \cdots [1]}. \] This is actually a polynomial in $\d$. By convention, $\qbinom{n}{k}=0$ for $k < 0$. \end{defn}

\begin{remark} \label{negativechoose} It is very important to observe that $\qbinom{n}{k}$ is well-defined and non-zero for $n < 0$ and $k>0$. Be careful: the familiar relation $\qbinom{n}{k} = \qbinom{n}{n-k}$ is false for $n<0$, except when both sides are zero. \end{remark}

\begin{defn} Fix a $\ZM[\d]$-algebra $\Bbbk$ and an integer $n \ge 2$. Let $\Fund_{\Bbbk,n}$, or $\Fund$ when $n$ and $\Bbbk$ are understood, denote the monoidal category defined as
below. The objects are sequences of elements of $I_n^+ \coprod I_n^-$, with monoidal structure given by concatenation. The morphisms from a sequence $\un{w}$ to a sequence $\un{x}$ are
the $\Bbbk$-linear span of $\sl_n$-webs with bottom $\un{w}$ and top $\un{x}$, modulo the relations below (and their reflections and orientation reversals). As discussed in
\cite[p8]{CKM}, some of these relations are redundant. \end{defn}

\begin{subequations} \label{sln-web-relations}

\noindent \textbf{Tag flip:}
\begin{equation} \label{tagflip} {
\labellist
\small\hair 2pt
 \pinlabel {$(-1)^{k(n-k)}$} [ ] at 43 16
 \pinlabel {$k$} [ ] at 77 -5
 \pinlabel {$n-k$} [ ] at 77 32
 \pinlabel {$k$} [ ] at 4 -5
 \pinlabel {$n-k$} [ ] at 4 32
\endlabellist
\centering
\ig{1.2}{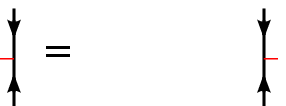}
} \end{equation}
\vskip .25cm

\noindent \textbf{Tag slide:}
\begin{equation} \label{tagslide} \ig{1.2}{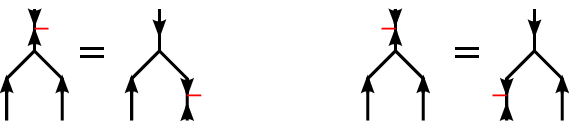} \end{equation}

\noindent \textbf{Tag removal:}
\begin{equation} \label{tagremoval} \ig{1.2}{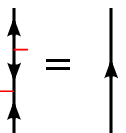} \end{equation}

\noindent \textbf{Bigon removal:}
\begin{equation} \label{bigonremoval} {
\labellist
\tiny\hair 2pt
 \pinlabel {$k$} [ ] at 219 4
 \pinlabel {$k$} [ ] at 164 4
 \pinlabel {$k+l$} [ ] at 77 4
 \pinlabel {$k+l$} [ ] at 21 4
 \pinlabel {$k$} [ ] at 164 51
 \pinlabel {$k+l$} [ ] at 21 51
 \pinlabel {$k+l$} [ ] at 174 18
 \pinlabel {$l$} [ ] at 146 18
 \pinlabel {$l$} [ ] at 22 18
 \pinlabel {$k$} [ ] at -2 18
 \pinlabel {$\displaystyle \qbinom{n-k}{l}$} [ ] at 196 29
 \pinlabel {$\displaystyle \qbinom{k+l}{l}$} [ ] at 49 29
\endlabellist
\centering
\ig{1}{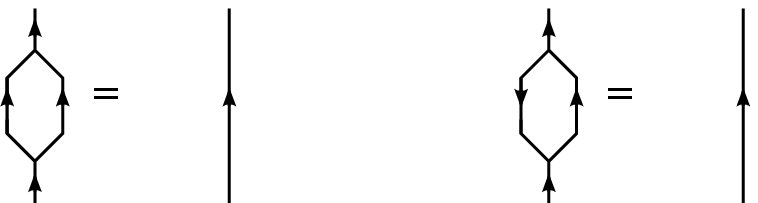}
} \end{equation}

\noindent \textbf{Circle removal:}
\begin{equation} \label{circleremoval} {
\labellist
\small\hair 2pt
 \pinlabel {$k$} [ ] at 23 20
 \pinlabel {$\displaystyle \qbinom{n}{k}$} [ ] at 45 13
\endlabellist
\centering
\ig{1}{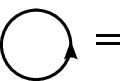}
} \end{equation}

\noindent \textbf{Associativity:}
\begin{equation} \label{associativity} {
\labellist
\tiny\hair 2pt
 \pinlabel {$k$} [ ] at 1 -5
 \pinlabel {$l$} [ ] at 18 -5
 \pinlabel {$m$} [ ] at 33 -5
 \pinlabel {$k$} [ ] at 58 -5
 \pinlabel {$l$} [ ] at 74 -5
 \pinlabel {$m$} [ ] at 90 -5
 \pinlabel {$k+l$} [ ] at -2 29
 \pinlabel {$l+m$} [ ] at 94 29
 \pinlabel {$k+l+m$} [ ] at 20 53
 \pinlabel {$k+l+m$} [ ] at 77 53
\endlabellist
\centering
\ig{1.2}{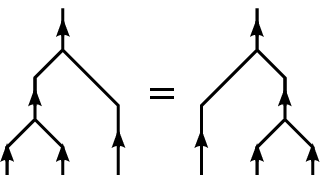}
} \end{equation}
\vskip .1cm
	
\noindent \textbf{Rung squash:}
\begin{equation} \label{rungsquash} {
\labellist
\tiny\hair 2pt
 \pinlabel {$l$} [ ] at 94 -5
 \pinlabel {$k$} [ ] at 70 -5
 \pinlabel {$l$} [ ] at 26 -5
 \pinlabel {$k$} [ ] at 2 -5
 \pinlabel {$l+s+r$} [ ] at 94 60
 \pinlabel {$k-s-r$} [ ] at 70 60
 \pinlabel {$l+s+r$} [ ] at 26 60
 \pinlabel {$k-s-r$} [ ] at 2 60
 \pinlabel {$r+s$} [ ] at 82 19
 \pinlabel {$s$} [ ] at 15 14
 \pinlabel {$r$} [ ] at 11 40
 \pinlabel {$l+s$} [ ] at 35 39
 \pinlabel {$k-s$} [ ] at -6 21
 \pinlabel {\normalsize $\displaystyle \qbinom{r+s}{r}$} [ ] at 55 28
\endlabellist
\centering
\ig{1.5}{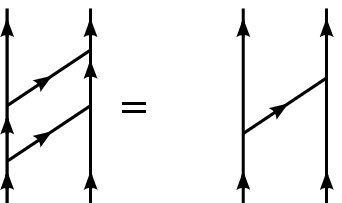}
} \end{equation}
\vskip .1cm

\noindent \textbf{Rung swap:}
\begin{equation} \label{rungswap} {
\labellist
\tiny\hair 2pt
 \pinlabel {$l$} [ ] at 130 -5
 \pinlabel {$k$} [ ] at 105 -5
 \pinlabel {$l$} [ ] at 26 -5
 \pinlabel {$k$} [ ] at 1 -5
 \pinlabel {$l+s-r$} [ ] at 130 60
 \pinlabel {$k+r-s$} [ ] at 106 60
 \pinlabel {$l+s-r$} [ ] at 26 60
 \pinlabel {$k+r-s$} [ ] at 1 60
 \pinlabel {$r-t$} [ ] at 115 10
 \pinlabel {$s-t$} [ ] at 113 43
 \pinlabel {$s$} [ ] at 17 12
 \pinlabel {$r$} [ ] at 16 44
 \pinlabel {$l+t-r$} [ ] at 142 28
 \pinlabel {$k+r-t$} [ ] at 116 28
 \pinlabel {$l+s$} [ ] at 20 28
 \pinlabel {$k-s$} [ ] at -7 28
 \pinlabel {\normalsize $\displaystyle \sum_t \qbinom{k-l+r-s}{t}$} [ ] at 72 26
\endlabellist
\centering
\ig{1.6}{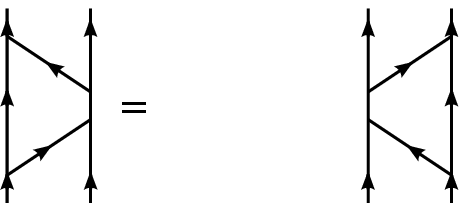}
} \end{equation}
\vskip .25cm

\end{subequations}

This ends the list of relations.

For future use, we record a special case of the rung swap \eqref{rungswap}. We assume that $a>b$.
\vskip .2cm
\begin{equation} \label{neutralundo} {
\labellist
\tiny\hair 2pt
 \pinlabel {$b$} [ ] at 130 -5
 \pinlabel {$a$} [ ] at 105 -5
 \pinlabel {$b$} [ ] at 26 -5
 \pinlabel {$a$} [ ] at 1 -5
 \pinlabel {$b$} [ ] at 130 60
 \pinlabel {$a$} [ ] at 106 60
 \pinlabel {$b$} [ ] at 26 60
 \pinlabel {$a$} [ ] at 1 60
 \pinlabel {$x$} [ ] at 115 12
 \pinlabel {$x$} [ ] at 113 41
 \pinlabel {$a-b$} [ ] at 17 10
 \pinlabel {$a-b$} [ ] at 17 45
 \pinlabel {$b-x$} [ ] at 140 28
 \pinlabel {$a+x$} [ ] at 114 28
 \pinlabel {$a$} [ ] at 22 28
 \pinlabel {$b$} [ ] at -4 28
 \pinlabel {\normalsize $\displaystyle \sum_x \qbinom{a-b}{x}$} [ ] at 72 26
\endlabellist
\centering
\ig{1.6}{rungswap}
} \end{equation}
\vskip .2cm
Note that the $x=0$ term on the RHS is actually just the identity map with coefficient $1$. Thus we can use this relation to replace the identity of the object $(a,b)$ with a sum of more complicated diagrams.

Finally, for an $\sl_n$-web $\phi$, we let $\overline{\phi}$ denote the same diagram flipped upside down, with orientation reversed. For example, the following two diagrams are interchanged by this \emph{duality involution}.
\begin{equation} \label{dualityexample} {
\labellist
\small\hair 2pt
 \pinlabel {$2$} [ ] at 11 8
 \pinlabel {$4$} [ ] at 35 8
 \pinlabel {$4$} [ ] at 61 8
 \pinlabel {$1$} [ ] at 11 76
 \pinlabel {$3$} [ ] at 35 76
 \pinlabel {$6$} [ ] at 61 76
 \pinlabel {$1$} [ ] at 124 8
 \pinlabel {$3$} [ ] at 148 8
 \pinlabel {$6$} [ ] at 172 8
 \pinlabel {$2$} [ ] at 124 76
 \pinlabel {$4$} [ ] at 148 76
 \pinlabel {$4$} [ ] at 172 76
\endlabellist
\centering
\ig{1}{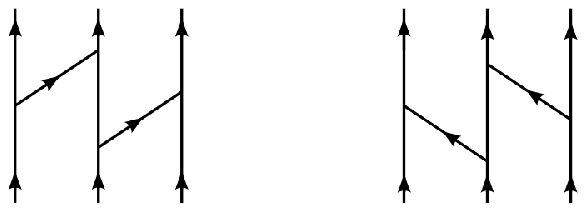}
} \end{equation}
The duality involution gives an isomorphism between $\Hom(\un{w}, \un{x})$ and $\Hom(\un{x}, \un{w})$ for any two sequences $\un{w}$ and $\un{x}$.

\subsection{Ladders}
\label{subsec-ladders}

In the category $\Fund$, a downward-oriented object labeled $k$ is isomorphic an upward-oriented object labeled $n-k$, via the tag morphism. Let $\Fund^+$ denote the full subcategory
whose objects are upward-oriented. We now begin a long investigation of morphisms in $\Fund^+$.

Consider an extended $\sl_n$-web of the following form. We assume that all indices live in $\hat{I}_n = \{0, \ldots, n\}$.
\begin{equation} \label{ladderexample} {
\labellist
\small\hair 2pt
 \pinlabel {$b$} [ ] at 98 -5
 \pinlabel {$a$} [ ] at 74 -5
 \pinlabel {$b$} [ ] at 25 -5
 \pinlabel {$a$} [ ] at 2 -5
 \pinlabel {$d$} [ ] at 98 45
 \pinlabel {$c$} [ ] at 74 45
 \pinlabel {$d$} [ ] at 25 45
 \pinlabel {$c$} [ ] at 2 45
 \pinlabel {$s$} [ ] at 86 25
 \pinlabel {$s$} [ ] at 13 25
 \pinlabel {or} [ ] at 48 19
\endlabellist
\centering
\ig{1}{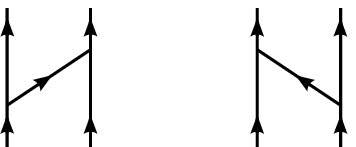}
} \end{equation}
\vskip .15cm

\noindent We refer to such a web as a \emph{rung}, and to compositions of such webs (tensored with identity maps) as \emph{ladders} (as in \cite[\S 5]{CKM}). For example, the morphisms
in \eqref{dualityexample} are ladders. The first rung pictured \emph{tilts northeast}, and the second \emph{tilts northwest}. For any ladder, the sum of the indices in the input is
equal to the sum of the indices in the output (in \eqref{dualityexample}, this sum is $10$).

Ladders are always upward-oriented, and consequently we may neglect to draw the orientations. Every $\sl_n$-web in the remainder of this paper will be upward-oriented (with one exception
that will be mentioned).

We think of a rung visually as a \emph{crossbar} connecting two \emph{uprights}. The edges in a ladder are either crossbars or \emph{segments of uprights}. When the crossbar has edge
label $s$ which is zero, the rung can effectively be ignored (it is the identity map). However, segments of uprights can be labeled $0$ and $n$ in interesting ways: thus, trivalent
vertices and cups and caps are all special examples of rungs, using \eqref{remove0n}. In fact, these morphisms in $\Fund^+$ can be viewed as rungs in multiple ways, so it behooves us to
distinguish between a ladder and the morphism in $\Fund^+$ that it represents.

Following \cite[\S 5]{CKM}, we let $\Lad^n_m$ denote the category of ladders with $m$ uprights built from extended $\sl_n$-webs. That is, the objects of $\Lad^n_m$ are sequences in
$\hat{I}_n$ of length $m$, and the morphisms are the formal span of ladders, being monoidally generated by rungs. No relations are imposed, except that rungs where the crossbar has
label $0$ are counted as the identity map and ignored. Together, $\Lad^n = \oplus_m \Lad^n_m$ is a monoidal category. There is a functor $\Phi \co \Lad^n \to \Fund^+$, which identifies
many objects, sending both $0$ and $n$ to the monoidal identity.

The set of non-extended $\sl_n$-webs which are the image of ladders under $\Phi$ form a suprisingly vast set of diagrams, and this flexibility comes from the fact that $0$ and $n$ are
allowed as segment labels. There are many complicated ladders expressing morphisms from the empty sequence to itself. Ladders are essentially just rigid versions of arbitrary webs. In
\cite[Theorem 5.3.1]{CKM}, it was proven that ladders (i.e. their images under $\Phi$) form a spanning set for morphisms in $\Fund^+$. The fact that ladders span is a topological
statement, akin to a Morse decomposition, and the proof in \cite{CKM} runs essentially along those lines.

Before moving on, let us describe some basic ladder manipulations.

\begin{defn} Consider a ladder in $\Lad^n_m$. To a rung between the $i$-th and $(i+1)$-st uprights, associate the letter $s_i$ if the rung is tilted northeast, and $t_i$ if the rung is
tilted northwest. To the ladder itself we associate the corresponding word in the set $\{s_i, t_i\}_{i=1}^{m-1}$, which we call its \emph{tilt-word}. (We ignore all rungs where the
crossbar has label $0$. This parenthetical will not be mentioned again.) \end{defn}

\begin{remark} There is a rung from $(0,a)$ to $(a,0)$, which after applying $\Phi$ is just the identity map of the underlying object $a$. However, with this rigid ladder description,
this rung is tilted northwest! The rung back from $(a,0)$ to $(0,a)$ is tilted northeast. When one fixes the number of uprights in a ladder, one can not ignore these identity rungs, and
must count them in the tilt-word. Similarly, the rung from $(n,a)$ to $(a,n)$ also is sent by $\Phi$ to the identity map of $a$, thanks to tag removal \eqref{tagremoval}. Note that a
rung where the crossbar has label $n$ is necessarily a morphism between the objects $(0,n)$ and $(n,0)$, and is sent by $\Phi$ to the empty diagram, but it still has a well-defined
tilt. \end{remark}

\begin{lemma}(Rung shuffling lemma) Choose a reduced expression for each element of the symmetric group $S_m$. After applying $\Phi$, any ladder in $\Lad^n_m$ can be written as a linear
combination of ladders in $\Lad^n_m$ whose tilt-words have the form $t_{i_1} t_{i_2} \cdots t_{i_d} s_{j_1} s_{j_2} \cdots s_{j_e}$, where $t_{i_1} \cdots t_{i_d}$, respectively
$s_{j_1} \cdots s_{j_e}$, are one of the chosen reduced expressions. \end{lemma}

\begin{proof} The relations \eqref{sln-web-relations} can be used to replace a ladder with a linear combination of ladders. Let us examine the mutations which are performed on the
corresponding tilt-words. To begin with, when $|i-j| \ge 2$, the mutations \[s_i s_j \to s_j s_i, \quad s_i t_j \to t_j s_i, \quad t_i t_j \to t_j t_i \] are realized by the
commutations of distant rungs.
\begin{equation} \label{ladderscommute} \ig{1}{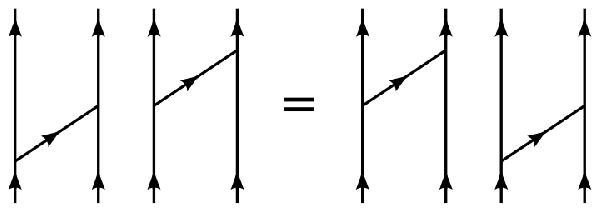} \end{equation}

Associativity \eqref{associativity} replaces $t_{i+1} s_i$ with $s_i t_{i+1}$, and vice versa.
\begin{equation} \label{laddersassoc} \ig{1}{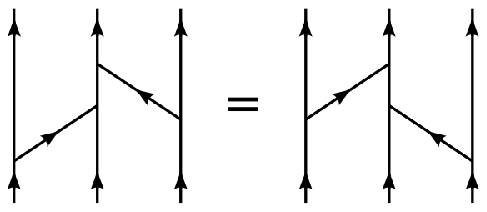} \end{equation}
Similarly, it replaces $t_i s_{i+1}$ with $s_{i+1} t_i$ and vice versa. Note that, when certain inputs and outputs are $0$ or $n$, one uses a tag slide \eqref{tagslide} instead of associativity. (In fact, tag sliding is a special case of associativity, being one of the redundant relations.)

Rung swap \eqref{rungswap} replaces $s_i t_i$ with $t_i s_i$ or a subword thereof (because some of the crossbars in the RHS of \eqref{rungswap} could have label $0$). Special cases
(when various labels are $0$ or $n$) include \eqref{bigonremoval} and \eqref{circleremoval}, which also replace $s_i t_i$ with a subword of $t_i s_i$. These operations can be applied in reverse, sending $t_i s_i$ to a subword of $s_i t_i$.

Thus $s_i$ and $t_k$ commute, modulo shorter words, for all $i$ and $k$. We can assume that all instances of $t$ occur before all instances of $s$, and treat them separately.

The rung squash relation \eqref{rungsquash} replaces $s_i s_i$ with $s_i$ (and similarly with $t_i$). The following lemma proves that $s_i s_{i+1} s_i$ can be replaced with $s_{i+1}
s_i s_{i+1}$ and subwords thereof, and vice versa. Thus the Coxeter relations for $S_m$ hold for the indices $s_i$, modulo shorter words. This is sufficient to imply the proposition.
\end{proof}

\begin{lemma} The following analog of the Coxeter relation $s_i s_{i+1} s_i = s_{i+1} s_i s_{i+1}$ holds in $\sl_n$-webs. We call it the \emph{R3 relation}. In this relation, $m = r+t$ is fixed, and $n$ and $l$ vary in the sum but satisfy $s = l+n$.

\begin{equation} \label{R3} {
\labellist
\tiny\hair 2pt
 \pinlabel {$a$} [ ] at 12 7
 \pinlabel {$b$} [ ] at 36 7
 \pinlabel {$c$} [ ] at 60 7
 \pinlabel {$a$} [ ] at 124 7
 \pinlabel {$b$} [ ] at 148 7
 \pinlabel {$c$} [ ] at 172 7
 \pinlabel {$x$} [ ] at 11 93
 \pinlabel {$y$} [ ] at 36 93
 \pinlabel {$z$} [ ] at 59 93
 \pinlabel {$x$} [ ] at 124 93
 \pinlabel {$y$} [ ] at 148 93
 \pinlabel {$z$} [ ] at 172 93
 \pinlabel {$r$} [ ] at 48 35
 \pinlabel {$s$} [ ] at 22 43
 \pinlabel {$t$} [ ] at 48 75
 \pinlabel {$n$} [ ] at 136 35
 \pinlabel {$m$} [ ] at 160 63
 \pinlabel {$l$} [ ] at 136 75
 \pinlabel {\normalsize $\displaystyle = \sum_{n} \qbinom{m-s}{t-n}$} [ ] at 90 53
\endlabellist
\centering
\ig{1.4}{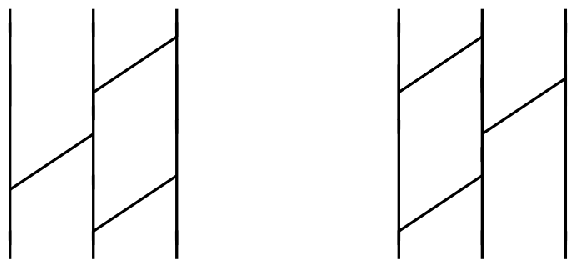}
} \end{equation}
 \end{lemma}

\begin{proof} The manipulations below are, in order: associativity \eqref{associativity} on the upright ending in $z$, followed by a rung swap \eqref{rungswap}, and then associativity on the upright starting with $a$. One sets $n = t-i$. Note that $r = b-d$, so that $b-s+t-d = r+t-s = m-s$.

\begin{equation} {
\labellist
\tiny\hair 2pt
 \pinlabel {$a$} [ ] at 64 110
 \pinlabel {$b$} [ ] at 88 110
 \pinlabel {$c$} [ ] at 112 110
 \pinlabel {$x$} [ ] at 64 197
 \pinlabel {$y$} [ ] at 88 197
 \pinlabel {$z$} [ ] at 112 197
 \pinlabel {$a$} [ ] at 152 110
 \pinlabel {$b$} [ ] at 200 110
 \pinlabel {$c$} [ ] at 225 110
 \pinlabel {$x$} [ ] at 152 197
 \pinlabel {$y$} [ ] at 176 197
 \pinlabel {$z$} [ ] at 225 197
 \pinlabel {$a$} [ ] at 35 6
 \pinlabel {$b$} [ ] at 84 6
 \pinlabel {$c$} [ ] at 109 6
 \pinlabel {$x$} [ ] at 35 92
 \pinlabel {$y$} [ ] at 60 92
 \pinlabel {$z$} [ ] at 109 92
 \pinlabel {$a$} [ ] at 176 6
 \pinlabel {$b$} [ ] at 200 6
 \pinlabel {$c$} [ ] at 224 6
 \pinlabel {$x$} [ ] at 176 92
 \pinlabel {$y$} [ ] at 200 92
 \pinlabel {$z$} [ ] at 224 92
 \pinlabel {$r$} [ ] at 99 136
 \pinlabel {$s$} [ ] at 75 149
 \pinlabel {$t$} [ ] at 100 178
 \pinlabel {$d$} [ ] at 85 137
 \pinlabel {$f$} [ ] at 91 160
 \pinlabel {$r$} [ ] at 205 148
 \pinlabel {$s$} [ ] at 161 133
 \pinlabel {$t$} [ ] at 187 160
 \pinlabel {$d$} [ ] at 186 131
 \pinlabel {$f$} [ ] at 172 146
 \pinlabel {$m$} [ ] at 208 180
 \pinlabel {$m$} [ ] at 90 76
 \pinlabel {$s$} [ ] at 48 21
 \pinlabel {$l$} [ ] at 55 56
 \pinlabel {$k$} [ ] at 89 57
 \pinlabel {$t-i$} [ ] at 73 39
 \pinlabel {$d-i$} [ ] at 72 74
 \pinlabel {$l$} [ ] at 187 75
 \pinlabel {$m$} [ ] at 212 62
 \pinlabel {$n$} [ ] at 188 34
 \pinlabel {$k$} [ ] at 204 44
 \pinlabel {$d-i$} [ ] at 193 62
 \pinlabel {\normalsize $=$} [ ] at 132 156
 \pinlabel {\normalsize $\displaystyle = \sum_i \qbinom{b-s+t-d}{i}$} [ ] at 4 50
 \pinlabel {\normalsize $\displaystyle = \sum_n \qbinom{m-s}{t-n}$} [ ] at 140 53
\endlabellist
\centering
\ig{1.5}{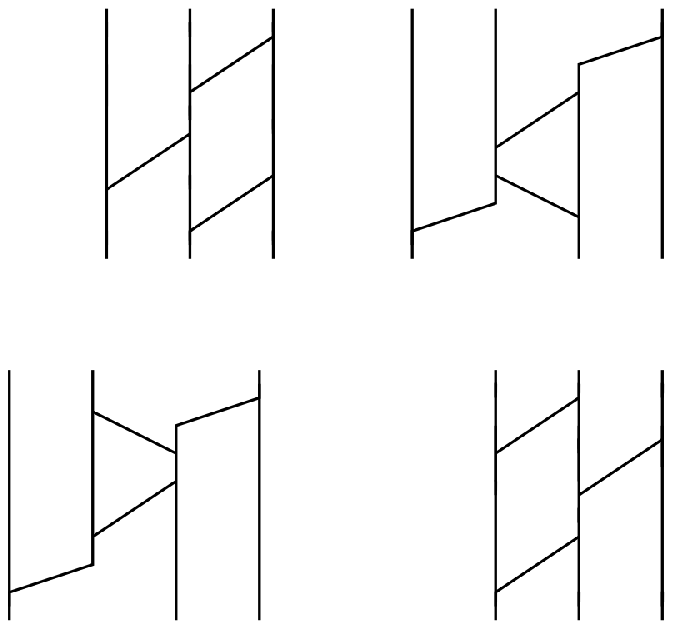}
} \label{R3proof} \end{equation}
\end{proof}

\begin{remark} We make a warning to the reader, recalling Remark \ref{negativechoose}. When $m-s$ is positive, the fact that $m-s - (t-n) = r-l$ implies that we may write
$\qbinom{m-s}{r-l}$ instead of $\qbinom{m-s}{t-n}$, so that the relation \eqref{R3} has extra symmetry. Thus when $m-s$ is positive, one must have $n \le t$ and $l \le r$ in order to
obtain a non-zero coefficient. However, when $m-s$ is negative, either $t-n$ or $r-l$ must be negative, but this is no contradiciton. \end{remark}

We call a rung of the form \eqref{ladderexample} \emph{neutral} if $a=d$ and $b=c$. In other words, it is the unique rung sending inputs $(a,b)$ to outputs $(b,a)$, and is the identity
if $a=b$. In the special case of \eqref{R3} where $a>b>c$ and all three rungs are neutral, it takes the following form.

\begin{equation} \label{R3neutral} {
\labellist
\tiny\hair 2pt
 \pinlabel {$a$} [ ] at 12 7
 \pinlabel {$b$} [ ] at 36 7
 \pinlabel {$c$} [ ] at 60 7
 \pinlabel {$a$} [ ] at 124 7
 \pinlabel {$b$} [ ] at 148 7
 \pinlabel {$c$} [ ] at 172 7
 \pinlabel {$c$} [ ] at 11 93
 \pinlabel {$b$} [ ] at 36 93
 \pinlabel {$a$} [ ] at 59 93
 \pinlabel {$c$} [ ] at 124 93
 \pinlabel {$b$} [ ] at 148 93
 \pinlabel {$a$} [ ] at 172 93
 \pinlabel {$b-c$} [ ] at 48 35
 \pinlabel {$a-c$} [ ] at 22 45
 \pinlabel {$a-b$} [ ] at 48 75
 \pinlabel {$a-b$} [ ] at 136 35
 \pinlabel {$a-c$} [ ] at 160 63
 \pinlabel {$b-c$} [ ] at 136 75
 \pinlabel {\normalsize $\displaystyle = $} [ ] at 95 53
\endlabellist
\centering
\ig{1.4}{R3}
} \end{equation}

\subsection{Ladder decompositions}
\label{subsec-ladderdecomp}

Our goal in this chapter is to describe a cellular basis for $\Fund^+$. This improves drastically on the aforementioned result of \cite{CKM}, by finding a combinatorial set of ladders
which form a basis. The first ingredient in a cellular structure is a poset, which will be the poset of dominant weights.

\begin{defn} Given an object $\un{w}$ of $\Fund^+$, let $\wt(\un{w}) = \sum w_i$, the \emph{weight} of $\un{w}$, denote the dominant $\sl_n$ weight which is the sum of the fundamental
weights in the sequence $\un{w}$. We say that $\un{w} \in P(\l)$ when $\wt(\un{w}) = \l$. Weights are equipped with the usual partial order, where $\l > \mu$ if $\l - \mu$ is the sum of
positive roots.

For an object $\un{w}$ in $\Lad^n$, we define $\wt(\un{w})$ analogously by saying that the indices $0$ and $n$ have weight zero. \end{defn}

We classify rungs as being \emph{inward}, \emph{outward}, or \emph{neutral}, as follows. Consider the diagrams in \eqref{ladderexample}. Clearly, $a + b = c + d$. If the pair $\{a,b\}$
and the pair $\{c,d\}$ are equal, the rung is \emph{neutral}; either it is the identity web and $s=0$, or it swaps the positions of $a$ and $b$. For example, the LHS of
\eqref{neutralundo} is a composition of two neutral ladders.

If the pair $\{a,b\}$ is inside the pair $\{c,d\}$ within the interval $[0,n]$, we call the rung \emph{outward}. Here are some examples of outward rungs and their images after applying
$\Phi$, which include as special cases merging trivalent vertices and caps. Outward rungs can not include splitting trivalent vertices or cups.

\begin{example} Suppose that $a=2$ and $b=4$. There are three outward rungs, which are the first three diagrams below. Some rungs take on special forms when $n$ is small. 
\vskip .1cm
\begin{equation} \label{downwardladderexample} {
\labellist
\small\hair 2pt
 \pinlabel {$2$} [ ] at 1 -5
 \pinlabel {$4$} [ ] at 26 -5
 \pinlabel {$1$} [ ] at 1 45
 \pinlabel {$5$} [ ] at 26 45
 \pinlabel {$2$} [ ] at 54 -5
 \pinlabel {$4$} [ ] at 78 -5
 \pinlabel {$5$} [ ] at 53 45
 \pinlabel {$1$} [ ] at 78 45
 \pinlabel {$2$} [ ] at 101 -5
 \pinlabel {$4$} [ ] at 126 -5
 \pinlabel {$6$} [ ] at 114 45
 \pinlabel {$2$} [ ] at 186 -5
 \pinlabel {$4$} [ ] at 212 -5
 \pinlabel {$1$} [ ] at 198 45
 \pinlabel {$2$} [ ] at 230 -5
 \pinlabel {$4$} [ ] at 253 -5
 \pinlabel {$1$} [ ] at 242 45
 \pinlabel {$2$} [ ] at 294 -5
 \pinlabel {$4$} [ ] at 318 -5
 \pinlabel {$n=5$} [ ] at 221 -15
 \pinlabel {$n=6$} [ ] at 305 -15
\endlabellist
\centering
\ig{1}{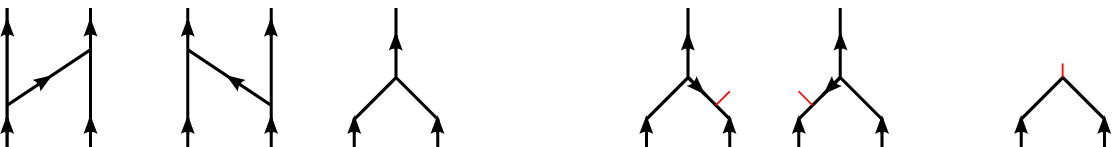}
} \end{equation}
\vskip .5cm
\noindent When $n \le 4$, all the outward rungs are zero. This is to say, they must have edges with labels less than $0$ or greater than $4$; such a thing can not exist in $\Lad^4$, and it is zero as an extended $\sl_4$-web. When $n=5$, the first two outward rungs become trivalent vertices, while the third becomes zero. Note that the two trivalent vertices are equal up to a sign, using \eqref{tagflip} and \eqref{tagslide}. When $n=6$, the third outward rung becomes a cup.  \end{example}

Conversely, if $\{c,d\}$ is inside $\{a,b\}$ within $[0,n]$, then the rung is \emph{inward}. For example, if $a=2$ and $b=4$ then there is one inward rung, for which
$c=d=3$. The duality involution swaps outward and inward rungs.

\begin{remark} \label{0nnotthere} Note that, while $0$ and $n$ are permitted to label any of the edges in a rung, they will never label an edge on the bottom of a (nonzero) outward
rung, or the top of a (nonzero) inward rung. \end{remark}

The reader should confirm that an outward rung will lower the weight (i.e. the target has lower weight than the source), an inward rung will raise it, and a neutral rung will preserve
it. This is to say that the idea of ``moving indices outward/inward'' is the same as the dominance order on $\sl_n$-weights.

\begin{remark} One could also consider outward rungs to be ``downward,'' and inward rungs to be ``upward,'' coinciding with their action on weights. This use of up and down is
consistent with the terminology surrounding light leaves in \cite{EWGr4sb}. However, most directional terminology is overloaded when discussing diagrams, and we hope that inward and
outward will cause the least confusion. The author will happily entertain suggestions for better terminology! \end{remark}

Let us note an easy observation about the dominance order.

\begin{lemma} \label{lem:maximal} Let $\un{x} = (x_1, \ldots, x_k)$ and $\un{y}=(y_1, \ldots, y_k)$ be objects in $\Fund^+$, and suppose that some $x_i$ is strictly larger than every $y_i$ and every other $x_j$. Then $\wt(\un{x}) \ngeq \wt(\un{y})$. \end{lemma}

\begin{proof} If $x_i$ is larger than all other $x_j$, then applying an outward or a neutral rung to the $i$-th upright can only increase (or preserve) the maximal label. Thus no sequence
of neutral or outward rungs can possibly transform $\un{x}$ into $\un{y}$. This is equivalent to the condition that one cannot subtract positive roots from $\wt(\un{x})$ to obtain
$\wt(\un{y})$. \end{proof}

We call a ladder \emph{neutral} if it is composed only of neutral rungs. We call a ladder \emph{outward} if it is composed only of outward and neutral rungs, and \emph{strictly
outward} if it has at least one outward rung. An \emph{inward} ladder is defined similarly. We say that a ladder has an \emph{in-out decomposition}, or we call it an \emph{IO
ladder}, if it can be written as an inward ladder above an outward ladder.  We say it has a \emph{strict in-out decomposition} if either the
outward or the inward ladder is strict. Here is an example of a strict IO ladder.

\[ {
\labellist
\small\hair 2pt
 \pinlabel {$2$} [ ] at 1 -5
 \pinlabel {$4$} [ ] at 26 -5
 \pinlabel {$1$} [ ] at -5 35
 \pinlabel {$5$} [ ] at 32 35
 \pinlabel {$3$} [ ] at 1 71
 \pinlabel {$3$} [ ] at 26 71
\endlabellist
\centering
\ig{1}{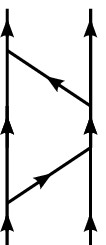}
} \]
\vskip .2cm

\begin{prop} \label{prop:UDdecomp} Morphisms in $\Fund^+$ are spanned by (the images under $\Phi$ of) IO ladders. \end{prop}

We prove this proposition in the next subsection, after some motivational preliminaries.

In contrast to general ladders, IO ladders are very restrictive. If both the source and the sink of an IO ladder are sent by $\Phi$ to the empty sequence, then the ladder itself is sent
to the empty diagram. A given object in $\Fund^+$ has many lifts to objects in $\Lad^n$ by adding instances of $0$ and $n$, but as mentioned in Remark \ref{0nnotthere}, allowing
additional $0$- and $n$-labeled strands in the source and the sink of a diagram will not truly broaden the set of IO ladders.

Intuition from representation theory motivates this proposition. When $\Bbbk = \QM(q)$ so that $\Fund$ represents morphisms in a semisimple category, we know that all morphisms factor
through projection to a common summand of the source and target. As the summands in a tensor product have lower weight than the input sequence, one expects a outward map to project from
the tensor product to the summand, and an inward map to include into the target. Soon, we will describe a set of outward ladders, called light ladders, which enumerate the projections
to summands.

\begin{remark} The reader may be interested in the interaction between the rung shuffling lemma and the labeling of rungs as inward, outward, or neutral. For example, the rung swap
\eqref{rungswap} can be used to replace a ladder with two rungs in the ``wrong order'' OI with a sum of ladders in the right order IO. This essentially guarantees the proposition for
ladders with only two uprights. However, associativity \eqref{associativity} is less pleasant, and will often produce ladders in the ``wrong order.'' \end{remark}

Because of the inherent confusion between ladders in $\Lad^n$ and their images in $\Fund^+$, we must make some technical arguments. Namely, it should be the case that the composition
$\a \circ \b$ in $\Fund^+$, where $\a$ is the image under $\Phi$ of an inward ladder, and $\b$ of an outward ladder, is itself the image of an IO ladder. However, $\a$ and $\b$ may come
from ladders with different numbers of uprights, and the positions of strands labeled $0$ and $n$ may not match well, so that the original ladders need not be composable. This problem
is easily dealt with.

\begin{lemma} \label{lem:composability} Let $\a$ and $\b$ denote two ladders in $\Lad^n$, and suppose that their images under $\Phi$ are composable. In other words, the target of $\b$
and the source of $\a$ are sequences in $\hat{I}_n$ which, after removing each instance of $0$ and $n$, are equal. Then there are ladders $\a'$ and $\b'$ living in $\Lad^n_m$ for some
$m$, for which $\a'$ and $\b'$ are composable, and for which $\a'$ and $\a$ (resp. $\b'$ and $\b$) have the same image under $\Phi$. \end{lemma}

\begin{proof} Let $a$ be the sum of the input labels of $\a$, and $b$ the sum for $\b$. Note that $a$ and $b$ are equal modulo $n$. Choose an integer $c$ which is equal to $a$ and $b$
modulo $n$, and satisfies $c \ge a$ and $c \ge b$. Choose an integer $m$ ``large enough'': it is at least as big as the number of uprights in $\a$ or $\b$, and may need to be even
larger, depending on the choice of $c$. The choice of $c$ and $m$ is ultimately irrelevant. Extend $\a$ and $\b$ to ladders of the same width $m$ by adding new uprights on the right
labeled $0$ and/or $n$, in such a way that the sum of the input labels is now equal to $c$ for both ladders. Continue to call these ladders $\a$ and $\b$. It must now be the case that
the output of $\b$ and the input of $\a$ agree when all copies of $0$ and $n$ are removed, and have the same number of $0$ and the same number of $n$. There is a neutral ladder in
$\Lad^n_m$ which goes from the output of $\b$ to the input of $\a$, and only consists of neutral rungs of the form $(0,x) \leftrightarrow (x,0)$ and $(n,x) \leftrightarrow (x,n)$. That
is, it consists of rungs which are sent to the identity map by $\Phi$. Then, for instance, we may let $\a' = \a$ and $\b'$ be the composition of $\b$ below the neutral ladder.
\end{proof}

Let us note one further implication of IO decompositions.

\begin{defn} For an object $\un{w}$ in $\Fund^+$, let its \emph{width} be the length of $\un{w}$ viewed as a sequence in $I_n$. For an object $\un{w}$ in $\Lad^n_m$, let its
\emph{width} be the width of $\Phi(\un{w})$, which is equal to the length of the sequence obtained by deleting all instances of $0$ or $n$. For a dominant weight $\l$, let its
\emph{width} be the width of any $\un{w} \in P(\l)$. \end{defn}

Width is a coarser indicator than weight. As noted above, any outward rung in $\Lad^n_m$ can only (weakly) decrease the width, as $0$ and $n$ are not inputs in any nonzero outward rung, but they
can be outputs. Thus, in an IO ladder, the width will decrease monotonically until the middle, and then increase again. The existence of
IO decompositions implies that $\sl_n$-webs need not be any wider than necessary. That is, if $m$ is the larger of the widths of two objects $\un{w}$ and $\un{y}$, then every morphism
from $\un{w}$ to $\un{y}$ is in the span of $\Phi(\Lad^n_m)$.

So, IO ladders have an ``hourglass figure,'' visually seen in \eqref{doubleintro}. Unfortunately, the terminology is suboptimal in this regard: outward ladders make things thinner, and
inward ones make them wider. Forgive me!

\subsection{Filtrations and neutral ladders}
\label{subsec-filtrations}

Given a ladder with an in-out decomposition, the object which is the target of the initial outward ladder and the source of the inward ladder is well-defined up to neutral ladders, and
thus it has a well-defined weight. We call this weight the \emph{middle weight} of the IO ladder.

\begin{defn} For a dominant weight $\l$, let $I_{\le \l}$ denote the span, within any morphism space in $\Fund^+$, of those IO ladders with middle weight $\mu$ such that $\mu \le \l$. Let $J_{\le \l}$ denote the span of all morphisms factoring through any object $\un{x}$ with $\wt(\un{x}) \le \l$. Define $I_{< \l}$ and $J_{< \l}$ similarly.
\end{defn}

Clearly $J_{\le \l}$ is a 2-sided (non-monoidal) ideal, almost by construction. If we assume Proposition \ref{prop:UDdecomp}, then one can show that $I_{\le \l} = J_{\le \l}$ and
$I_{\le \l}$ is also a 2-sided ideal.

\begin{claim} \label{cor:filtration} Assume Proposition \ref{prop:UDdecomp}. Then the set $I_{\le \l}$ is a 2-sided (non-monoidal) ideal in $\Fund^+$. \end{claim}

\begin{proof} Let $\un{x} \in P(\mu)$ for $\mu \le \l$. Let $\phi \co \un{w} \to \un{x}$ be an outward ladder, and $\psi \co \un{x} \to \un{y}$ be an inward ladder. If $a \co \un{y} \to
\un{z}$ is any morphism, then $a \circ \psi$ can be rewritten as a linear combination of IO ladders, i.e. $a \psi = \sum \a_i \b_i$ where $\b_i$ is an outward ladder and $\a_i$ is an
inward ladder (we have ignored the coefficients). Thus $a \psi \phi = \sum \a_i (\b_i \phi)$ where $\b_i \phi$ is an outward ladder. Moreover, $\b_i$ is an outward ladder with source of
weight $\mu$, so its target has weight $\nu \le \mu \le \l$. Thus, $a \psi \phi \in I_{\le \l}$. A similar argument shows that $I_{\le \l}$ is closed under right multiplication.
\end{proof}

We have tacitly used Lemma \ref{lem:composability} in this proof, as we assumed that the composition of an inward ladder and an outward ladder is an IO ladder. We will continue to use
Lemma \ref{lem:composability} in this way, without further comment.

\begin{claim} Assume Proposition \ref{prop:UDdecomp}. Then $I_{\le \l} = J_{\le \l}$. \end{claim}

\begin{proof} Let $\phi \co \un{w} \to \un{x}$ and $\psi \co \un{x} \to \un{y}$ be arbitrary morphisms, where $\wt(\un{x}) = \mu \le \l$. We seek to prove that $\psi \phi$ is inside
$I_{\le \l}$. Rewrite $\psi = \sum_i \a_i \b_i$ as a linear combination of IO ladders (we have ignored the coefficients). The middle weight of each IO ladder $\a_i \b_i$ is some weight
$\nu_i \le \mu$. Now, rewrite $\b_i \phi = \sum_j \g_{ij} \d_{ij}$ as a linear combination of IO ladders, which have middle weight $\rho_{ij} \le \nu_i \le \mu$. Then $\psi \phi =
\sum_i \sum_j (\a_i \g_{ij}) \d_{ij}$ is a sum of IO ladders, and lies within $I_{\le \l}$. \end{proof}

We prove these claims before proving Proposition \ref{prop:UDdecomp} to give an idea of the style of argument, and introduce the filtrations in question. When $\l$ is understood, the
ideal $J_{< \l}$ will be thought of as \emph{lower terms}, since its morphisms (are sums of morphisms which) factor through objects of lower weight.

Now we show that neutral ladders are isomorphisms in the associated graded of this filtration. Technically, we have not yet shown that they are nonzero in the associated graded; they are still isomorphisms, but may be the zero isomorphism between zero objects.  That neutral ladders are nonzero will be proven in \S\ref{subsec-lightladdersindep}. The following proof does not assume Proposition \ref{prop:UDdecomp}.

\begin{lemma} (Neutral ladder lemma) \label{lem:neutralladderlemma} For any dominant weight $\l$ and $\un{w}, \un{x} \in P(\l)$, any two neutral ladders $\un{w} \to \un{x}$ are equal modulo $J_{< \l}$. \end{lemma}

\begin{cor} \label{cor:identityspansgr} For any dominant weight $\l$ and $\un{w} \in P(\l)$, $\End(\un{w})$ is spanned by the identity map modulo $J_{< \l}$. (We have not yet proven
that the identity of $\End(\un{w})$ is nonzero modulo $J_{< \l}$.) \end{cor}

In fact, Corollary \ref{cor:identityspansgr} implies the Neutral Ladder Lemma. After all, if $\phi$ and $\psi$ are the two neutral ladders with different source and target,
and we know Corollary \ref{cor:identityspansgr}, then $\phi \equiv \phi \overline{\psi} \psi \equiv \psi$ modulo lower terms, because $\overline{\psi} \psi \equiv \1$ and $\phi
\overline{\psi} \equiv \1$ modulo lower terms.

\begin{proof} The proof is a modification of the rung shuffling lemma. A neutral ladder (on $m$ uprights) is effectively just a permutation of the labels on the uprights, and we can
record it with a sequence of simple reflections in $S_m$, analogous to our previous recording of tilt-words. After all, a neutral rung with inputs $(a,b)$ must either tilt northwest if
$b>a$ or northeast if $a>b$ (or be the identity map and disappear if $a=b$), so keeping track of tilt is not important. Instead, we just record a word which associates $s_i$ to a rung
between the $i$-th and $(i+1)$-st strands. The word uniquely determines the neutral ladder (given the input sequence $\un{w} \in P(\l)$).
	
It is enough to show that one can apply analogs of the Coxeter relations for the symmetric group to neutral ladders, modulo $I_{< \l}$. For example, \eqref{neutralundo} corresponds to
the relation $s_i^2 = 1$ in the symmetric group, because it implies that a ``bigon'' of neutral ladders is equal to the identity modulo lower terms. The distant braid relation $s_i s_j
= s_j s_i$ follows from the monoidal structure.

In \eqref{R3neutral} above we have checked the Reidemeister III relation $s_i s_{i+1} s_i = s_{i+1} s_i s_{i+1}$ for a particular case, where the input strands $(a,b,c)$ satisfied
$a>b>c$. The remaining cases follow from this one by adding bigons. They can also be dealt with on a case by case basis.

For example, let us check Reidemeister III when the input labels are $(a,c,b)$, still with $a<b<c$. We wish to show that the following diagrams are equal modulo lower terms.
\begin{equation} \label{RIII} {
\labellist
\small\hair 2pt
 \pinlabel {$a$} [ ] at 12 6
 \pinlabel {$c$} [ ] at 36 6
 \pinlabel {$b$} [ ] at 60 6
 \pinlabel {$a$} [ ] at 89 6
 \pinlabel {$c$} [ ] at 112 6
 \pinlabel {$b$} [ ] at 136 6
 \pinlabel {$b$} [ ] at 12 92
 \pinlabel {$c$} [ ] at 36 92
 \pinlabel {$a$} [ ] at 60 92
 \pinlabel {$b$} [ ] at 89 92
 \pinlabel {$c$} [ ] at 112 92
 \pinlabel {$a$} [ ] at 136 92
\endlabellist
\centering
\ig{1}{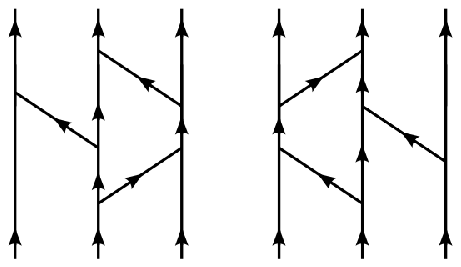}
} \end{equation}
Consider the first diagram. We apply associativity \eqref{associativity} to the bottom, and then a rung swap \eqref{rungswap}.
\begin{equation} {
\labellist
\small\hair 2pt
 \pinlabel {$a$} [ ] at 11 7
 \pinlabel {$c$} [ ] at 36 7
 \pinlabel {$b$} [ ] at 60 7
 \pinlabel {$b$} [ ] at 11 92
 \pinlabel {$c$} [ ] at 36 92
 \pinlabel {$a$} [ ] at 60 92
 \pinlabel {$a$} [ ] at 89 7
 \pinlabel {$c$} [ ] at 111 7
 \pinlabel {$b$} [ ] at 138 7
 \pinlabel {$b$} [ ] at 89 92
 \pinlabel {$c$} [ ] at 111 92
 \pinlabel {$a$} [ ] at 138 92
 \pinlabel {$\displaystyle \sum_t \xi_t$} [ ] at 171 50
 \pinlabel {$a$} [ ] at 191 7
 \pinlabel {$c$} [ ] at 216 7
 \pinlabel {$b$} [ ] at 240 7
 \pinlabel {$t$} [ ] at 228 58
 \pinlabel {$b$} [ ] at 191 92
 \pinlabel {$c$} [ ] at 216 92
 \pinlabel {$a$} [ ] at 240 92
 \pinlabel {$a$} [ ] at 268 7
 \pinlabel {$c$} [ ] at 292 7
 \pinlabel {$b$} [ ] at 316 7
 \pinlabel {$b$} [ ] at 268 92
 \pinlabel {$c$} [ ] at 292 92
 \pinlabel {$a$} [ ] at 316 92
 \pinlabel {$\equiv$} [ ] at 253 52
 \pinlabel {\tiny $c+a-b$} [ ] at 300 50
\endlabellist
\centering
\ig{1.2}{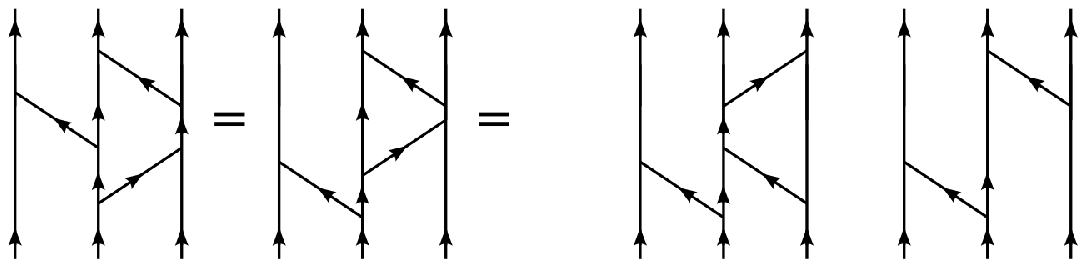}
} \end{equation}

\noindent The coefficients arising from the rung swap (some binomial coefficients) are denoted $\xi_t$, but are ultimately irrelevant. After applying the rung swap, each term with $t > 0$ will factor through lower terms.
One can confirm that the $t=0$ term has coefficient $1$, which gives the rightmost diagram. An analogous computation with the second diagram in \eqref{RIII} gives the same result.
\end{proof}

Now we finally prove that IO ladders span.

\begin{proof}[Proof of Proposition \ref{prop:UDdecomp}] We seek to prove the following auxiliary statement. Let $\un{w} \in P(\l)$ and $\un{x} \in P(\mu)$ for dominant weights $\l,
\mu$. Then outward ladders span $\Hom(\un{w}, \un{x})$ modulo $J_{< \mu}$. We call this statement $S(\l, \mu)$. Note that this statement is quite strong: $S(0,0)$ already implies that
any diagram with empty boundary reduces to the empty diagram, the only outward ladder, as $J_{< 0} = 0$.

Let us rigidify the situation by proving a refinement $S(\l, \mu, m)$: all morphisms in $\Hom(\un{w}, \un{x})$ in the image of $\Lad^n_m$ are within the span of outward ladders from
$\Lad^n_m$ modulo $J_{< \mu}$. The number of uprights $m$ is bounded below, but not above, by the size of $\l$ and $\mu$. Any ladder with $m$ uprights can be viewed as a ladder with
$m+1$ uprights by adding an upright labeled by $0$ or $n$, without changing its image under $\Phi$. Clearly, showing $S(\l, \mu, m)$ for all $m \ge 0$ will prove $S(\l, \mu)$.

We prove $S(\l, \mu, m)$ the result by induction on $\l$ and $\mu$. We will not need $S(0,0, m)$ as a base case; the argument we use will apply equally to $S(0,0, m)$, showing that
outward ladders span morphisms (modulo lower terms, but there are no lower terms). Our argument requires some caution as there are a number of ``invisible ladders'' one may bump into
(see the remark before the rung shuffling lemma). We have already been very careful to ensure that all our previous statements hold even when labels involve $0$ and $n$.

Let us fix two objects of $\Lad^n_m$, $\un{\hat{w}} \in P(\l)$ and $\un{\hat{x}} \in P(\mu)$. Using neutral ladders we may reorder the output sequence $\un{\hat{x}}$. By the neutral
ladder lemma, these reorderings are isomorphisms modulo $J_{< \mu}$; since an outward ladder composed with a neutral ladder is still outward, it is sufficient to prove the result for
any reordering of $\un{\hat{x}}$. So we assume henceforth that $\un{\hat{x}} = (x_1, x_2, \ldots, x_m)$ satisfies $x_1 \le x_2 \le \ldots \le x_m$. Note that all instances of $x_i=0$
happen at the start, and all instances of $x_i = n$ happen at the end.

Consider a ladder $L$ from $\un{\hat{w}}$ to $\un{\hat{x}}$. If the topmost rung is inward, then the ladder must lie in $J_{< \mu}$. Let the bottommost rung be $R$, so that $L = L'
\circ R$. If $R$ is outward, then the source of $L'$ has weight strictly less than $\l$, and induction implies that $L'$ is in the span of outward ladders modulo $J_{< \mu}$; thus so is
$L' \circ R$. Meanwhile, if the bottommost rung $R$ is neutral, then $L$ is in the span of outward ladders if $L'$ is, so we may as well replace $L$ with $L'$. Thus we can assume that
the bottommost rung is inward, while the topmost rung is either outward or neutral.

We may also use the rung shuffling lemma to assume that, within $L$, all northwest tilting rungs appear above all northeast tilting rungs. However, a northwest tilting rung with target
$(x_i, x_{i+1})$ for which $x_i \le x_{i+1}$ is necessarily an inward rung, contradicting our assumption above. Thus, there are no northwest tilting rungs, and every rung tilts
northeast. Of course, if there are no northeast tilting rungs either, then we are done.

We now induct on the number of rungs in our northeast-tilting ladder. The induction step involves an operation to the bottommost rung $R$ of the diagram. This rung has inputs $(a,b)$ and outputs $(c,d)$, tilts northeast, and is inward. This implies that $a > b$. Recall that the $x=0$ term of the RHS of \eqref{neutralundo}
is just the identity map. We use \eqref{neutralundo} to replace the identity of the $(a,b)$ uprights with a sum of more complicated diagrams, to be dealt with in turn.
\begin{equation} \label{reductionstep}
{
\labellist
\tiny\hair 2pt
 \pinlabel {$a$} [ ] at 12 34
 \pinlabel {$b$} [ ] at 36 34
 \pinlabel {$c$} [ ] at 12 104
 \pinlabel {$d$} [ ] at 36 104
 \pinlabel {$a$} [ ] at 76 7
 \pinlabel {$b$} [ ] at 100 7
 \pinlabel {$b$} [ ] at 70 51
 \pinlabel {$a$} [ ] at 96 52
 \pinlabel {$a$} [ ] at 73 82
 \pinlabel {$b$} [ ] at 95 83
 \pinlabel {$c$} [ ] at 76 122
 \pinlabel {$d$} [ ] at 100 124
 \pinlabel {$\displaystyle + \sum_{x > 0} \xi_x$} [ ] at 122 69
 \pinlabel {$a$} [ ] at 148 7
 \pinlabel {$b$} [ ] at 170 7
 \pinlabel {$a+x$} [ ] at 145 50
 \pinlabel {$a-x$} [ ] at 183 51
 \pinlabel {$a$} [ ] at 153 77
 \pinlabel {$b$} [ ] at 178 91
 \pinlabel {$c$} [ ] at 148 125
 \pinlabel {$d$} [ ] at 172 125
 \pinlabel {$\displaystyle \sum_{t \ge 0} \gamma_t$} [ ] at 215 70
 \pinlabel {$a$} [ ] at 235 7
 \pinlabel {$b$} [ ] at 260 6
 \pinlabel {$b$} [ ] at 231 45
 \pinlabel {$a$} [ ] at 256 54
 \pinlabel {$b-t$} [ ] at 232 85
 \pinlabel {$a+t$} [ ] at 269 82
 \pinlabel {$c$} [ ] at 236 124
 \pinlabel {$d$} [ ] at 259 126
 \pinlabel {$\displaystyle + \sum_{x > 0} \xi'_x$} [ ] at 288 67
 \pinlabel {$a$} [ ] at 320 7
 \pinlabel {$b$} [ ] at 343 6
 \pinlabel {$a+x$} [ ] at 317 55
 \pinlabel {$b-x$} [ ] at 355 58
 \pinlabel {$c$} [ ] at 320 124
 \pinlabel {$d$} [ ] at 344 124
\endlabellist
\centering
\ig{1}{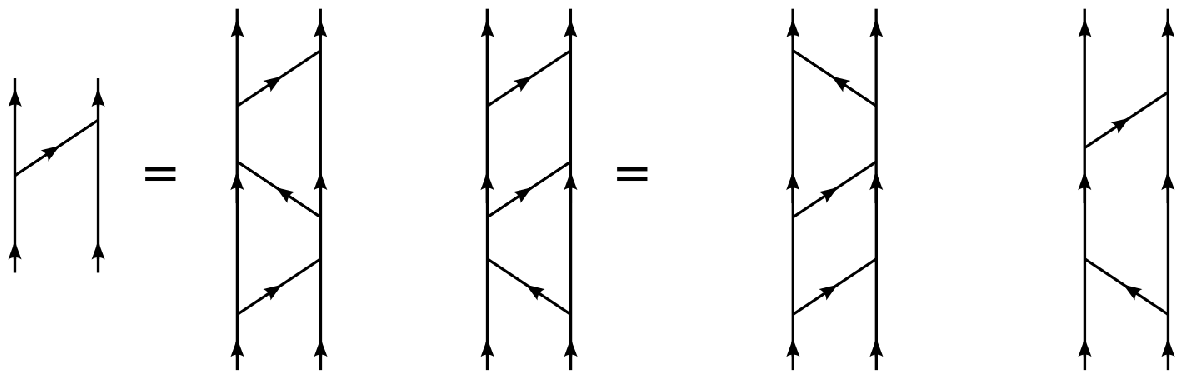}
} \end{equation}
The first equality is \eqref{neutralundo}, while the second equality is \eqref{rungswap} and \eqref{rungsquash} applied to the two terms, respectively. The exact coefficients $\xi_x$, $\xi'_x$ and $\gamma_t$ are irrelevant for this argument.

In each term where $x > 0$ or where $t > 0$, the diagram begins with an outward ladder from $\l$ to a lower weight, and we may use induction as before to assume these terms are in the
span of outward ladders. The remaining term, where $t=0$, is a neutral rung followed by the \emph{seesawed} version of our original rung: the inputs are now $(b,a)$ instead of $(a,b)$, and the tilt is northwest rather than northeast. Call this seesawed rung $R'$.

\begin{equation} \label{seesaw1} \quad {
\labellist
\small\hair 2pt
 \pinlabel {$a$} [ ] at 8 4
 \pinlabel {$b$} [ ] at 31 4
 \pinlabel {$b$} [ ] at 73 4
 \pinlabel {$a$} [ ] at 96 4
 \pinlabel {$c$} [ ] at 8 71
 \pinlabel {$d$} [ ] at 31 71
 \pinlabel {$c$} [ ] at 73 71
 \pinlabel {$d$} [ ] at 96 71
 \pinlabel {$R = $} [ ] at -10 40
 \pinlabel {$ = R'$} [ ] at 110 40
\endlabellist
\centering
\ig{1}{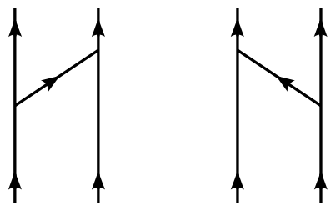}
} \end{equation}
The point is that $R$ and $R'$ are equal, modulo neutral ladders on bottom (required for the sources to match up), and modulo lower terms (diagrams with an outward rung on bottom).

Let us now take $R'$ and commute it past all northeast-tilting rungs, as we would in the rung shuffling lemma. When it commutes past distant rungs as in \eqref{ladderscommute} or past
adjacent rungs using associativity as in \eqref{laddersassoc}, nothing much happens. When it commutes past a rung on the same two uprights using the rung swap \eqref{rungswap}, the
label on both rungs may be lowered, possibly to zero. If some northwest-tilting rung makes its way to the top unscathed, then by the ordering on $\un{\hat{x}}$ it must be inward, and
the result would lie in $J_{< \mu}$. Thus, in every surviving relevant diagram, the northwest-tilting rung had its label lowered to zero from various rung swaps.

However, the remaining diagram has the same shape as $L'$ (the original ladder $L$ without $R$), possibly with some of its rungs removed (also with labels lowered to zero from the rung
swaps). In particular, it too is purely northeast-tilting. Now induction concludes the proof of $S(\l, \mu)$.

Passing from $S(\l, \mu)$ for all $\l, \mu$ to a proof of Proposition \ref{prop:UDdecomp} is straightforward, along the lines of Claim \ref{cor:filtration}. Let us prove inductively
that $J_{\le \l} = I_{\le \l}$. When $\l=0$, $J_{< 0} = 0$, and $S(\mu, 0)$ implies that any morphism from $\un{w} \in P(\mu)$ to the empty sequence is in the span of outward ladders.
Flipping this conclusion upside-down, $S(\nu, 0)$ implies that any morphism from the empty sequence to $\un{x} \in P(\nu)$ is in the span of inward ladders. Thus any morphism $\un{w}
\to \un{x}$ which factors through the empty sequence is in the span of IO ladders. This proves the base case.

The proof of the inductive step is similar. Consider a morphism from $\un{w} \in P(\mu)$ to $\un{x} \in P(\nu)$ which factors through an object of weight $\l$. Rewrite to map $\un{w}
\to \l$ as a sum of outward ladders modulo $J_{< \l}$, and do the same for the map $\l \to \un{x}$ (we abuse notation and write $\l$ for any sequence with weight $\l$). Then the result
is in the span of IO ladders with middle weight $\l$, plus terms in $J_{< \l}$. By induction, the terms in $J_{< \l}$ are themselves in the span of IO ladders with middle weight $< \l$.
Thus the whole morphism is in the span of IO ladders with middle weight less than or equal to $\l$. This concludes the proof. \end{proof}

\subsection{Light ladders}
\label{subsec-lightladders}

In the associated graded of our filtration by dominant weights, one need only consider IO ladders whose middle is a given weight $\l$. Next, we choose a special set of IO ladders with
this middle. In this section, we will not distinguish between a ladder and its image under $\Phi$.

\begin{defn} \label{def:weightsubseq} Let $\un{w} = (\w_{i_1}, \ldots, \w_{i_d})$ be a sequence of fundamental weights. A \emph{weight subsequence} $\eb \subset \un{w}$ is a sequence
$\eb = (e_0, e_1, \ldots, e_d)$ of weights, for which $e_0 = 0$ is the trivial weight, and $e_k - e_{k-1}$ is some weight appearing in the fundamental representation for $\w_{i_k}$ for
each $k \ge 1$. It is called a \emph{dominant weight subsequence} or a \emph{miniscule Littelmann path} if each $e_i$ is a dominant weight. Otherwise, to emphasize that $\eb$ is not
dominant, we may call it a \emph{false weight subsequence}. Regardless, we say that $\eb$ \emph{expresses} the weight $e_d$, and we also write $\un{w}^{\eb} = e_d$. We let $E(\un{w})$
denote the set of dominant weight subsequences, and let $E(\un{w}, \l)$ denote the set of dominant weight subsequences expressing the dominant weight $\l$. Note that, if $\un{w} \in
P(\l)$ then $E(\un{w}, \l)$ consists of a single element, which we call the \emph{full weight subsequence}.\end{defn}

One can think of a weight subsequence $\eb \subset \un{w}$ as being a choice, for each fundamental weight $\w_{i_k}$, of one of its weight spaces (which is $e_k - e_{k-1}$). However, not
all such choices are viable as dominant weight subsequences. Having chosen the first $m$ weights, the next weight can only be chosen such that their sum is still a dominant weight.

Although the definitions above make sense for arbitrary semisimple Lie algebras, the following motivational proposition relies on special features in type $A$.

\begin{prop} \label{prop:motivational} Inside $\Rep_q'$ (i.e. after base change to $\QM(q)$), the tensor product $V_{\un{w}} = \ot V_{\w_{i_k}}$ is isomorphic to $\oplus_{\eb
\subset \un{w}} V_{\un{w}^\eb}$. Thus, for two sequences $\un{w}$ and $\un{y}$, the dimension of $\Hom(V_{\un{w}}, V_{\un{y}})$ is equal to the number of pairs $\eb \subset \un{w}, \fb
\subset \un{y}$ such that $\un{w}^\eb = \un{y}^\fb$. \end{prop}

\begin{proof} Observe that when a weight $\mu$ of a fundamental representation $V_{\w}$ is added to a dominant weight $\l$, the result $\l + \mu$ is either dominant, or lies on a
wall of the shifted dominant chamber (shifted by $-\rho$, the half-sum of positive roots). The principles of plethyism then state that $V_\l \ot V_\w \cong \oplus V_{\l + \mu}$, where
the sum is over weights $\mu$ of $V_{\w}$ with $\l + \mu$ dominant. Iterating this, we get the desired result. \end{proof}

\begin{remark} Outside of type $A$, $\l + \mu$ could lie on the other side of a wall of the shifted dominant chamber, and the principles of plethyism would cause this weight to cancel out
one of the summands coming from a dominant weight $\l + \mu'$. Enumerating the summands of a tensor product of fundamental representations is somewhat more complicated. \end{remark}

We cannot use this proposition directly in $\Fund$ because we work over an arbitrary $\ZM[\d]$-algebra $\Bbbk$, where semisimplicity and standard plethyism may fail. However, it indicates
how big one should expect morphism spaces to be. We will now explicitly construct a basis of $\Hom(\un{w}, \un{y})$ in bijection with such pairs $(\eb, \fb)$ such that $\un{w}^\eb =
\un{y}^\fb$, which we call the \emph{double ladders basis}.

For each dominant weight subsequence $\eb \subset \un{w}$ we will follow an algorithm to construct a map from $\un{w}$ to $\un{x}$ for some $\un{x} \in P(\un{w}^\eb)$, which we call a
\emph{light ladder} associated to $\eb$. This map will be a outward ladder. There is ambiguity in the algorithm, in the choice of $\un{x}$ and in other choices along the way, but we
will just choose arbitrarily one map obtainable by the algorithm. Composing this light ladder with the duality involution applied to another light ladder will produce a double ladder,
which naturally has an in-out decomposition.

Light ladders will be defined inductively, and will be built in \emph{tiers}. Suppose we have already constructed some light ladder $L_d$ from $\un{w}$ to $\un{x}$, associated to $\eb \in E(\un{w})$. We draw $L_d$ as a trapezoid, which is justified by the fact that the width will (weakly) decrease in an outward ladder. The subscript $d$ indicates that $\un{w}$ has width $d$.
\[ {
\labellist
\small\hair 2pt
 \pinlabel {$\un{w}$} [ ] at 95 3
 \pinlabel {$L_d$} [ ] at 45 19
 \pinlabel {$\un{x}$} [ ] at 72 39
\endlabellist
\centering
\ig{1}{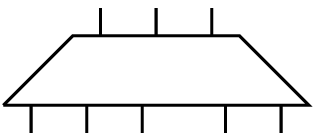}
} \]

Now, consider the sequence obtained by concatenating $\un{w}$ with a new strand labeled $a$, for $1 \le a \le n-1$. There are ${n \choose a}$ weights in the fundamental representation for
$\w_a$, and thus ${n \choose a}$ options for extending $\eb$ to a weight subsequence of $\un{w}a$, although not all of them are viable (that is, not all will produce an element of
$E(\un{w}a)$). In each case we apply some neutral ladder to $\un{x}$ to obtain a desired configuration on the right; we apply a particular outward ladder called an \emph{elementary light
ladder} to this configuration with the new strand $a$; and then we apply some neutral ladder at the end. The result is $L_{d+1}$, a light ladder for $\un{w}a$. The elementary light
ladder, together with the neutral ladders before and after, constitute one \emph{tier} of the multi-tier light ladder.

\begin{example} Consider the example where $n=4$ and $a=2$. In each of the diagrams below, a rectangle represents an arbitrary neutral ladder with the given top and bottom, and
unlabeled strands take on arbitrary values. Each of these diagrams represents a possible choice of the trapezoid $L_{d+1}$, corresponding to the six $\sl_4$-weights in $V_2$. The
elementary light ladder is the part of each diagram ignoring $L_d$ and the neutral ladders, and the last tier is the part ignoring $L_d$.

\begin{equation} \label{LL42} {
\labellist
\tiny\hair 2pt
 \pinlabel {$2$} [ ] at 106 93
 \pinlabel {$2$} [ ] at 231 93
 \pinlabel {$2$} [ ] at 187 142
 \pinlabel {$1$} [ ] at 223 167
 \pinlabel {$3$} [ ] at 187 167
 \pinlabel {$2$} [ ] at 368 18
 \pinlabel {$1$} [ ] at 323 67
 \pinlabel {$3$} [ ] at 324 91
 \pinlabel {$2$} [ ] at 367 146
 \pinlabel {$3$} [ ] at 322 192
 \pinlabel {$1$} [ ] at 322 218
 \pinlabel {$2$} [ ] at 493 93
 \pinlabel {$3$} [ ] at 451 140
 \pinlabel {$1$} [ ] at 434 141
 \pinlabel {$1$} [ ] at 462 169
 \pinlabel {$2$} [ ] at 435 188
 \pinlabel {$2$} [ ] at 623 93
 \pinlabel {$2$} [ ] at 579 147
 \pinlabel {\large $(0,1,0)$} [ ] at 52 67
 \pinlabel {\large $(1,-1,1)$} [ ] at 182 67
 \pinlabel {\large $(1,0,-1)$} [ ] at 375 238
 \pinlabel {\large $(-1,0,1)$} [ ] at 387 43
 \pinlabel {\large $(-1,1,-1)$} [ ] at 455 67
 \pinlabel {\large $(0,-1,0)$} [ ] at 570 67
 \pinlabel {$\un{w}$} [ ] at 0 88
 \pinlabel {$L_d$} [ ] at 51 103
 \pinlabel {$\un{x}$} [ ] at 10 118
\endlabellist
\centering
\ig{.75}{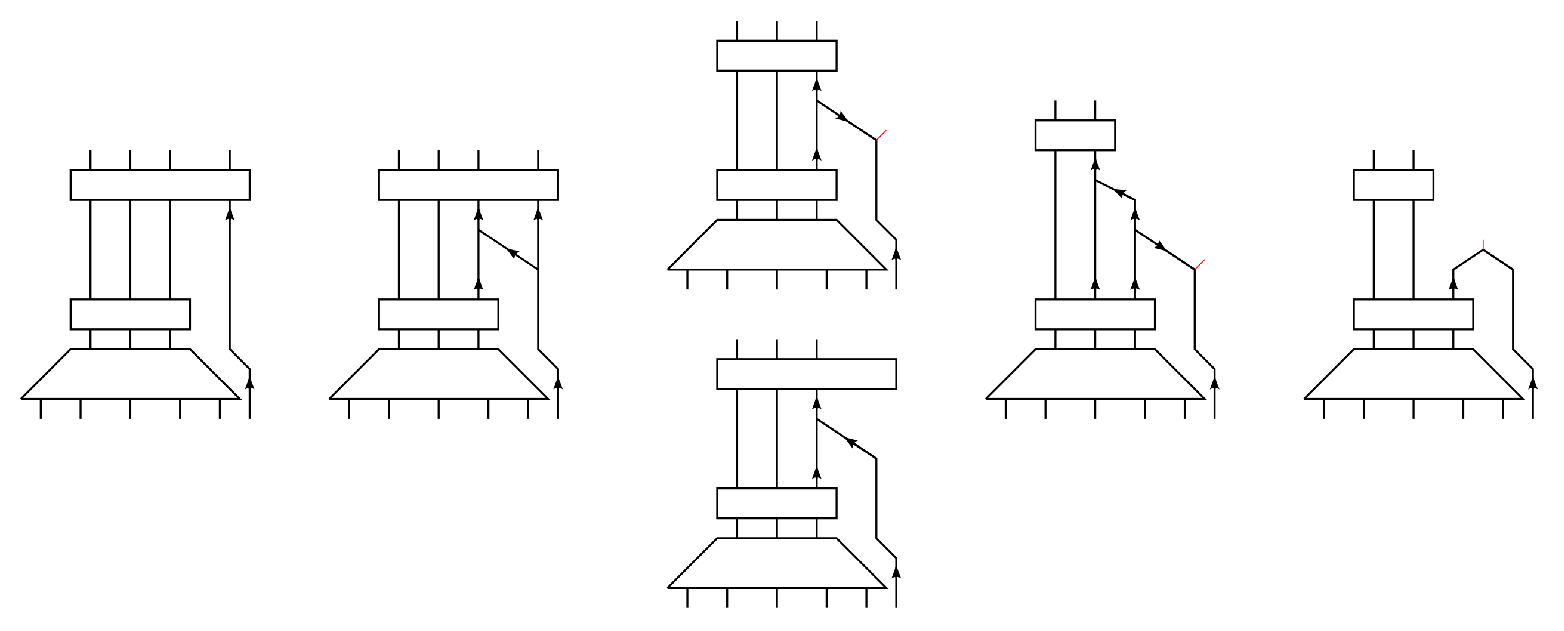}
} \end{equation}

The six elementary light ladders are not all necessarily viable for a given sequence $\un{x}$. For example, if $\un{x}$ does not have any strands labeled $3$, then two of the maps
pictured are impossible; these correspond to false weight sequences, which pass through a non-dominant weight. Also, note that the elementary light ladders are outward ladders, even
when they look like trivalent vertices or caps with tags, c.f. \eqref{downwardladderexample}. \end{example}

It remains to describe the ${n \choose a}$ possibilities for the elementary light ladder in the general case. Except the identity map, corresponding to the highest weight of $V_{\w_a}$,
each elementary light ladder will be strictly outward.

When discussing weights in fundamental representations, it will be much more convenient to view them as $\gl_n$-weights. Let $\Om(a)$ denote the set of \emph{01-sequences}, sequences of
0s and 1s having length $n$, for which there are $a$ total 1s and $n-a$ total 0s. There is a bijection between $\Om(a)$ and the weights appearing in $V_{\w_a}$, and we will always encode
a weight $\mu$ in $V_{\w_a}$ by an element of $\Om(a)$.

Any 01-sequence alternates between \emph{0-strings}, maximal consecutive subsequences of 0s, and \emph{1-strings}. By convention, we number these strings so that $\mu$ always starts with
a 1-string (possibly empty) and ends with a 0-string (possibly empty). Corresponding to $\mu$ are numbers $0 \le y_1 < x_1 < y_2 < x_2 < \ldots < x_k < y_{k+1} \le n$, where $y_1$ is the
last index in the first 1-string (and $y_1=0$ if the first 1-string is empty), $x_1$ is the last index in the first 0-string, and so forth, so that $y_{k+1}$ is the last index in the last
1-string (and $y_{k+1} = n$ if the last 0-string is empty). Note that $k$ is the number of 0-strings, not counting the last.

\begin{example} Suppose that $a=6$ and $n=11$ and $\mu = (01101001110)$. Then $k=3$, and $y_1 = 0$, $x_1 = 1$, $y_2 = 3$, $x_2 = 4$, $y_3 = 5$, $x_3 = 7$, and $y_4 = 10$. \end{example}

\begin{example} Suppose that $a=6$ and $n=11$ and $\mu = (10001100111)$. Then $k=2$, and $y_1 = 1$, $x_1 = 4$, $y_2 = 6$, $x_2 = 8$, and $y_3 = 11$. \end{example}

We also record numbers $y_1 = \a_1 < \a_2 < \ldots < \a_{k+1} = a$ and $0 < \b_1 < \ldots < \b_k < n$, which keep track of a cumulative count of 0s and 1s. Let $\a_i$ be the number of 1s in the first $i$ 1-strings, and $\b_i$ be the number of 0s in the first $i$ 0-strings. It is easy to observe that $\a_i + \b_i = x_i$, and $\a_i + \b_{i-1} = y_i$.

\begin{example} When $\mu = (01101001110)$, then $\a_1 = 0$, $\a_2 = 2$, $\a_3 = 3$, and $\a_4 = 6$; $\b_1 = 1$, $\b_2 = 2$, and $\b_3 = 4$. \end{example}

\begin{example} When $\mu = (10001100111)$, then $\a_1 = 1$, $\a_2 = 3$, and $\a_3 = 6$; $\b_1 = 3$ and $\b_2 = 5$. \end{example}

Then, corresponding to $\mu$, we have the following ladder, which we call the \emph{elementary light ladder} $E_{\mu}$ associated to $\mu$.
\begin{equation} \label{elementarylightladder} {
\labellist
\tiny\hair 2pt
 \pinlabel {$y_1$} [ ] at 12 101
 \pinlabel {$y_2$} [ ] at 36 101
 \pinlabel {$\ldots$} [ ] at 60 101
 \pinlabel {$y_{k+1}$} [ ] at 108 101
 \pinlabel {$y_k$} [ ] at 84 101
 \pinlabel {$x_1$} [ ] at 19 29
 \pinlabel {$x_2$} [ ] at 44 29
 \pinlabel {$\ldots$} [ ] at 69 29
 \pinlabel {$x_k$} [ ] at 92 29
 \pinlabel {$a$} [ ] at 108 8
 \pinlabel {$\b_1$} [ ] at 20 62
 \pinlabel {$\a_2$} [ ] at 32 55
 \pinlabel {$\b_2$} [ ] at 45 62
 \pinlabel {$\a_k$} [ ] at 80 55
 \pinlabel {$\b_k$} [ ] at 95 62
\endlabellist
\centering
\ig{1.2}{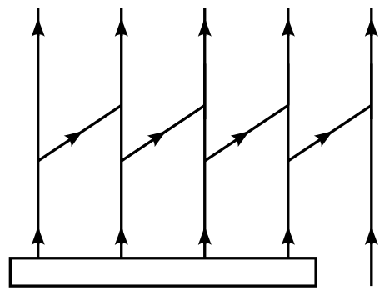}
} \end{equation}
The increasing sequence $\b_i$ labels the crossbars, while the increasing sequence $\a_i$ labels the ``middle segment'' of each upright. The number of rungs is $k$.

One can also think of this construction inductively. Let $\mu^- \in \Om(\a_k)$ denote $\mu$ with the final 1-string removed. Then $E_{\mu^-}$ is $E_{\mu}$ with the rightmost rung removed. Beginning with the zero weight or with a fundamental weight, one can add 1-strings one at a time to obtain $\mu$, and can produce rungs one at a time to get $E_{\mu}$.

\begin{example} When $\mu = (01101001110)$, $\mu^- = (01101000000)$. They have the same values of $x_i$ and $\b_i$ for $i < 3$, and the same values of $y_i$ and $\a_i$ for $i \le 3$.
	
The sequence of weights $\mu = (01101001110) \to (01101000000) \to (01100000000) \to (00000000000)$ produces the elementary light ladder pictured below.
	
\begin{equation} {
\labellist
\small\hair 2pt
	\pinlabel {$6$} [ ] at 84 8
	\pinlabel {$7$} [ ] at 54 28
	\pinlabel {$4$} [ ] at 30 29
	\pinlabel {$1$} [ ] at 5 29
	\pinlabel {$10$} [ ] at 84 100
	\pinlabel {$5$} [ ] at 59 100
	\pinlabel {$3$} [ ] at 35 100
	\pinlabel {$3$} [ ] at 65 52
	\pinlabel {$2$} [ ] at 40 68
\endlabellist
\centering
\ig{1}{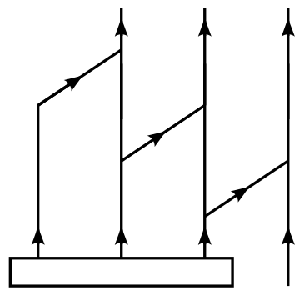}
} \end{equation} \end{example}

Note the following special cases. If $y_{k+1} = n$, so that the final 0-string is empty, then the rightmost rung is secretly a trivalent vertex with tag, thanks to \eqref{remove0n}. If $y_1 = 0$ so that the first 1-string is empty, then the leftmost rung is secretly a trivalent vertex. The most important special case is when $k=0$, or equivalently, $\mu = \w_a = (11111100000)$ is a fundamental weight. In this case, there are no rungs in $E_{\mu}$, and it is just the identity map of a strand labeled $a$.

\begin{example} Let $n=7$ and $\mu = (1101011)$. There are three 1-strings, but only two rungs in the elementary light ladder, because the ``third rung" is the identity map of the upright labeled $2$.

\begin{equation} {
\labellist
\small\hair 2pt
 \pinlabel {$5$} [ ] at 60 6
 \pinlabel {$5$} [ ] at 30 28
 \pinlabel {$3$} [ ] at 6 28
 \pinlabel {$4$} [ ] at 36 100
 \pinlabel {$2$} [ ] at 12 100
 \pinlabel {$3$} [ ] at 40 52
\endlabellist
\centering
\ig{1}{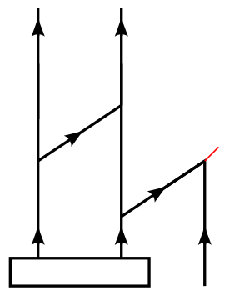}
} \end{equation} \end{example}

The reader should confirm that the diagrams in \eqref{LL42} are all elementary light ladders. The most complicated example is the weight $\mu = (0101)$, denoted there by its $\sl_n$-weight $(-1, 1, -1)$, which is the only example having $k > 1$.

\begin{claim} A ladder of the form
\begin{equation} \label{loweringform} {
\labellist
\small\hair 2pt
 \pinlabel {$y_1$} [ ] at 12 101
 \pinlabel {$y_2$} [ ] at 36 101
 \pinlabel {$\ldots$} [ ] at 60 101
 \pinlabel {$y_{k+1}$} [ ] at 108 101
 \pinlabel {$y_k$} [ ] at 84 101
 \pinlabel {$x_1$} [ ] at 19 29
 \pinlabel {$x_2$} [ ] at 44 29
 \pinlabel {$\ldots$} [ ] at 69 29
 \pinlabel {$x_k$} [ ] at 92 29
 \pinlabel {$a$} [ ] at 108 8
\endlabellist
\centering
\ig{1}{loweringmapclaim}
} \end{equation}
is an elementary light ladder if and only if $0 \le y_1 < x_1 < y_2 < x_2 < \ldots < x_k < y_{k+1} \le n$. The corresponding weight $\mu$ in $V_{\w_a}$ is uniquely determined from $x_i$ and $y_i$. \end{claim}

\begin{proof} The bijection between $\Om(a)$ and sequences $y_1 < x_1 < \ldots < x_k < y_{k+1}$ is obvious. \end{proof}

It is important to note that any northeast-tilted outward rung is an elementary light ladder for which $k=1$. This will be the base case of an inductive argument in a later section.

This concludes the algorithmic construction of a ladder attached to a dominant weight subsequence. To reiterate, at each step the elementary light ladder is deterministic, given by the
choice of weight, but the neutral ladders applied before and after are arbitrary, as is the target of the map.

\begin{defn} Any ladder constructible via the above algorithm, with elementary light ladders determined by a dominant weight subsequence $\eb \in E(\un{w}, \l)$, will be called a
\emph{light ladder} for $\eb$. \end{defn}

In order to obtain a basis, we should make some arbitrary choices to obtain one light ladder for each $\eb$.

\begin{defn} \label{defn:doubleladders} Henceforth, we fix a sequence of fundamental weights $\un{x}_\l$ for each dominant weight $\l$. For each $\eb \in E(\un{w}, \l)$, we choose one light ladder giving a morphism
from $\un{w}$ to $\un{x}_\l$. We call this map $LL_{\eb}$, \emph{the (chosen)} light ladder for $\eb$. We use the convention that, if $\un{w} = \un{x}_\l$ and $\eb \in E(\un{w},\l)$ is the full subsequence, then $LL_{\eb}$ is the identity map. Given $\eb \subset \un{w}$ and $\fb \subset \un{y}$ which express a common weight $\l$, we define the \emph{double (light) ladder} $\LL^\l_{\eb, \fb}$ to be the composition
\[ {
\labellist
\small\hair 2pt
 \pinlabel {$LL_{\eb}$} [ ] at 44 17
 \pinlabel {$\overline{LL}_{\fb}$} [ ] at 44 46
 \pinlabel {$\un{w}$} [ ] at 91 0
 \pinlabel {$\un{x}_\l$} [ ] at 75 32
 \pinlabel {$\un{y}$} [ ] at 91 67
\endlabellist
\centering
\ig{1}{doubleladder}
} \]
where $\overline{LL}_{\fb}$ denotes the map $LL_{\fb}$ upside-down with reverse orientation, i.e. the duality involution applied to $LL_{\fb}$.

We will write $M(\l, \un{y})$ for the same set of dominant weight subsequences as $E(\un{y}, \l)$, but thought of as indexing upside-down light ladders, rather than light ladders
themselves. Thus, the set of double ladders in $\Hom(\un{w}, \un{y})$ is indexed by $M(\l, \un{y}) \times E(\un{w}, \l)$, and this notation mimics the usual notation for the composition
of morphism spaces. \end{defn}

Note that light ladders are themselves double ladders, where the upside-down light ladder happens to be the identity map. Conversely, if a double ladder is outward, then it is just a
light ladder.

We claim that double ladders form a basis for morphism spaces in $\Fund^+$. This will be proven over the next few sections.

\begin{remark} \label{lightisgeneral} The algorithm for constructing light ladders is extremely analogous to that for constructing light leaves in the context of diagrammatic Soergel
bimodules \cite{EWGr4sb}. Neutral ladders are to $\sl_n$-webs what rex moves are to Soergel diagrams. Our proof that double ladders form a basis will be similar to the proof for double
leaves in \cite{EWGr4sb}. That is, we will show linearly independence after evaluating under some functor, and we will show spanning by a terribly convoluted diagrammatic argument.

Both $\sl_n$-webs and diagrammatic Soergel bimodules are \emph{strictly object adapted cellular categories} or SOACCs, as defined in \cite{ELauda}. However, they are both also monoidal
categories. The interplay between cellular structures and monoidal structures has never been fully explored. We would like to suggest that algorithms like the light ladder algorithm are
actually quite general for monoidal SOACCs. \end{remark}

\subsection{Double ladders span: part I}
\label{subsec-lightladdersspan1}

We wish to show that double ladders span. We will prove this theorem in tandem with several auxiliary results.

\begin{thm}[(Strong) Double Ladder Theorem]\label{thm:doubleladdersspan} All morphism spaces in $\Fund^+$ are spanned by double light ladders. \end{thm}

\begin{lemma}[(Strong) Light Ladder Lemma] \label{lem:lightladdersspan} Let $\un{w}$ be an object in $\Fund^+$, and fix $\l$. Then the (chosen) light ladders $LL_{\eb}$ for $\eb \in
E(\un{w}, \l)$ span $\Hom(\un{w}, \un{x}_\l)$ modulo $I_{< \l}$. In addition, for any other $\un{x} \in P(\l)$, the space $\Hom(\un{w}, \un{x})$ modulo $I_{< \l}$ is spanned by $N \circ
LL_{\eb}$, where $N$ is any fixed neutral ladder from $\un{x}_\l$ to $\un{x}$. \end{lemma}

There is also a weaker version of the above results, allowing for more flexibility in the choice of neutral ladders.

\begin{lemma}[(Weak) Light Ladder Lemma] \label{lem:lightladdersweak} Let $\un{w}$ be an object in $\Fund^+$, and fix $\un{y} \in P(\l)$. Then $\Hom(\un{w}, \un{y})$ is spanned by light
ladders (not the chosen ones, but anything constructible by the algorithm) for $\eb \in E(\un{w}, \l)$, modulo $I_{< \l}$. \end{lemma}

The statement of the Weak Double Ladder Theorem is analogous. The proof is long and convoluted, but can be divided
into two major parts. The first part will be proving the following proposition, which is the real crux of the proof of the Double Ladder Theorem. It relies in a very concrete way on the exact definitions of the elementary light ladders.

\begin{prop} \label{thecrux} Fix sequences $\un{w}$ and $\un{z}$ with $\un{z} \in P(\l)$. Then, modulo $I_{< \l}$, any morphism $L \co \un{w} \to \un{z}$ is in the span of
morphisms of the form $N \circ E \circ X$
\begin{equation} \label{eq:thecrux} {
\labellist
\small\hair 2pt
 \pinlabel {$X$} [ ] at 44 20
 \pinlabel {$E$} [ ] at 62 47
 \pinlabel {$N$} [ ] at 63 74
\endlabellist
\centering
\ig{1}{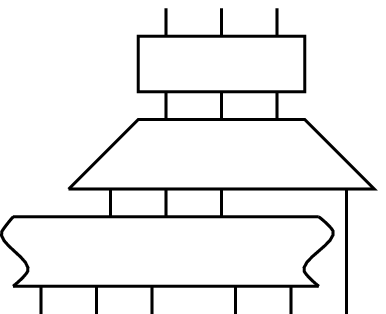}
} \end{equation}
where $N$ is a neutral ladder, $E$ is an elementary light ladder, and $X$ is an arbitrary morphism with one fewer input strand. \end{prop}

The remainder of the proof is a convoluted induction, which bootstraps Proposition \ref{thecrux} into a proof of many other results. However, it does not rely in any way on the
definition of the elementary light ladders, relying on the algorithmic ``tiered'' construction of light ladders without caring how the tiers are actually defined. It also
relies on the filtration by ideals $I_{< \l}$, and on the existence of IO decompositions. The same proof would apply to any class of morphisms which satisfied Proposition \ref{thecrux}.
In essence, this proof should be broadly applicable to monoidal categories with an (object-adapted) cellular filtration, to construct a cellular basis adapted to the monoidal structure.
This paragraph expands on Remark \ref{lightisgeneral}.

This section constitutes the proof of this crucial proposition. The next section will contain the general bootstrapping arguments. For the rest of this section, fix sequences $\un{w} \in P(\mu)$ and $\un{z} \in P(\l)$. We will be considering a morphism $L \co \un{w} \to \un{z}$, modulo $I_{< \l}$.

\begin{lemma} \label{lem:cruxlemma} Modulo $I_{< \l}$, $L$ is a linear combination of ladders of the form
\begin{equation} \label{cruxlemma} {
\labellist
\small\hair 2pt
 \pinlabel {$X$} [ ] at 48 19
 \pinlabel {$N$} [ ] at 60 112
\endlabellist
\centering
\ig{1}{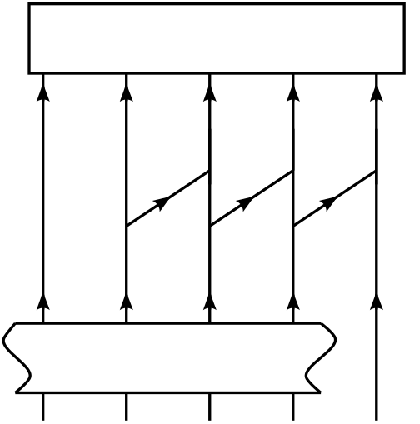}
}, \end{equation}
where $X$ is an arbitrary ladder, and $N$ is a neutral ladder.  \end{lemma}

\begin{proof} Because ladders span, we may assume that $L$ is already expressed as a ladder with $m$ rungs. By the Neutral Ladder Lemma, we may assume without loss of generality that
the diagram ends with a neutral ladder from $\un{y}$ to $\un{z}$, where $\un{y}$ is in order, with $y_1 \le y_2 \le \cdots$. Now we apply the rung shuffling lemma to the map from
$\un{w}$ to $\un{y}$, asserting that all northwest-tilting rungs occur on top. Such rungs must be inward, so if there is any northwest-tilting rung, the diagram lives in $I_{< \l}$.
Hence we can assume that there are only northeast-tilting rungs. Moreover, we can choose any reduced expression for each element of the symmetric group $S_m$, which will determine the
position of the northeast-tilting rungs. We choose our reduced expression so that it begins with an expression for an element of $S_{m-1}$, and concludes with a minimal coset
representative of $S_{m-1}$. The corresponding ladder has exactly the form of \eqref{cruxlemma} (without the neutral ladder on top), where the visible rungs correspond to the minimal
coset representative. \end{proof}

When we discuss a ladder of the form \eqref{cruxlemma}, we let $E$ denote the part of the ladder which is neither $N$ nor $X$. The target of $E$ is $\un{y} \in P(\l)$, and the source is $\un{x}a$. Thus $E$ is as follows.
\begin{equation} \label{lfredux} {
	\labellist
	\small\hair 2pt
	 \pinlabel {$y_1$} [ ] at 12 101
	 \pinlabel {$y_2$} [ ] at 36 101
	 \pinlabel {$\ldots$} [ ] at 60 101
	 \pinlabel {$y_{k+1}$} [ ] at 108 101
	 \pinlabel {$y_k$} [ ] at 84 101
	 \pinlabel {$x_1$} [ ] at 19 29
	 \pinlabel {$x_2$} [ ] at 44 29
	 \pinlabel {$\ldots$} [ ] at 69 29
	 \pinlabel {$x_k$} [ ] at 92 29
	 \pinlabel {$a$} [ ] at 108 8
	\endlabellist
	\centering
	\ig{1}{loweringmapclaim}
} \end{equation}
though possibly with an identity morphism to the left. We can also assume, as in the proof above, that $\un{y}$ is in order.

\begin{lemma} \label{lem:cruxlemma2} Modulo $I_{< \l}$, $L$ is a linear combination of ladders of the form \eqref{cruxlemma}, with the additional restriction that all the rungs in $E$ are outward or neutral. \end{lemma}

\begin{proof} We claim that if $E$ has an inward rung, then the diagram lies within $I_{< \l}$. Suppose that the rightmost rung of $E$ is inward, and the rest are outward or neutral. All other cases reduce to this one, by replacing $E$ with a subdiagram of $E$. We continue to assume that $\un{y}$ is ordered.

If all but the last rung is outward or neutral then $x_i \le y_{i+1}$ for $1 \le i < k$, and if the last rung is inward, then $x_k > y_{k+1}$ and $x_k > a$. But since $\un{y}$ is
ordered, we see that $x_k$ is larger than any other input $x_i$ or $a$, and is larger than any output $y_i$. Hence, by Lemma \ref{lem:maximal}, we see that the source of $E$ is not a
larger weight than $\l$, so that any IO ladder must factor through a weight strictly less than $\l$. \end{proof}

Now to prove the proposition, we need only prove the following lemma.

\begin{lemma} \label{rewritelemma} Any ladder $E$ as in \eqref{lfredux}, built of outward or neutral rungs, can be rewritten as a linear combination of elementary light ladders, up to
some operations: removing a neutral rung from the top (it becomes part of $N$), ignoring diagrams with an inward rung on top (they are in $I_{< \l}$), removing a neutral or an outward
rung from the bottom on the first $m-1$ uprights (it becomes part of $X$). (We will not actually need to add any inward rungs to $X$.) \end{lemma}

\begin{proof} We argue by induction on the number of rungs of $E$, or equivalently, the number of uprights involved. The identity map is an elementary light ladder, handing the
case of no rungs. Every northeast-tilted outward rung is an elementary light ladder, which handles the case of one outward rung. If the topmost rung of \eqref{lfredux} is
neutral (including the case of one rung which is neutral) then it can be absorbed into $N$, leaving behind a ladder of the form \eqref{lfredux} with fewer rungs, to which we can apply
induction. Thus we can assume that there is more than one rung, and the topmost rung is outward.

The fact that each of the rungs is outward or neutral gives some inequalities on the labels $x_i$ and $y_i$. For example, we know that $x_i \le y_{i+1}$, with equality if and only if the
$i$-th rung is neutral. We know that $y_1 < x_1 < y_2$ because the topmost rung is outward. There are several inequalities involving $\a_i$ as well. It is straightforward to deduce that
$y_1$ is strictly smaller than all the other inputs and outputs, and all the $\a_i$ as well.

The ladder is an elementary light ladder if $y_1 < x_1 < y_2 < x_2 < \ldots < x_k < y_{k+1}$. Problems can occur if either $y_i \ge x_i$, or if $x_i = y_{i+1}$ and a rung is neutral. We
progress by finding the first inequality which is not as desired, and applying a seesawing argument similar to that used in the proof of Proposition \ref{prop:UDdecomp}.

Rewriting \eqref{reductionstep} upside-down, we see that
\begin{equation} \label{seesaw} {
\labellist
\tiny\hair 2pt
 \pinlabel {$c$} [ ] at 11 31
 \pinlabel {$d$} [ ] at 36 31
 \pinlabel {$a$} [ ] at 11 110
 \pinlabel {$b$} [ ] at 36 110
 \pinlabel {\small $\displaystyle \sum_{t \ge 0}$} [ ] at 77 70
 \pinlabel {$c$} [ ] at 103 8
 \pinlabel {$d$} [ ] at 128 8
 \pinlabel {$b+t$} [ ] at 113 53
 \pinlabel {$a-t$} [ ] at 138 53
 \pinlabel {$a$} [ ] at 103 124
 \pinlabel {$b$} [ ] at 128 124
 \pinlabel {\small $\displaystyle \sum_{x > 0}$} [ ] at 174 67
 \pinlabel {$c$} [ ] at 203 8
 \pinlabel {$d$} [ ] at 229 8
 \pinlabel {$\equiv$} [ ] at 251 66
 \pinlabel {$a-x$} [ ] at 214 53
 \pinlabel {$b+x$} [ ] at 239 53
 \pinlabel {$a$} [ ] at 203 124
 \pinlabel {$b$} [ ] at 228 124
 \pinlabel {$c$} [ ] at 272 8
 \pinlabel {$d$} [ ] at 295 8
 \pinlabel {$b$} [ ] at 277 53
 \pinlabel {$a$} [ ] at 301 53
 \pinlabel {$a$} [ ] at 272 124
 \pinlabel {$b$} [ ] at 295 124
\endlabellist
\centering
\ig{1}{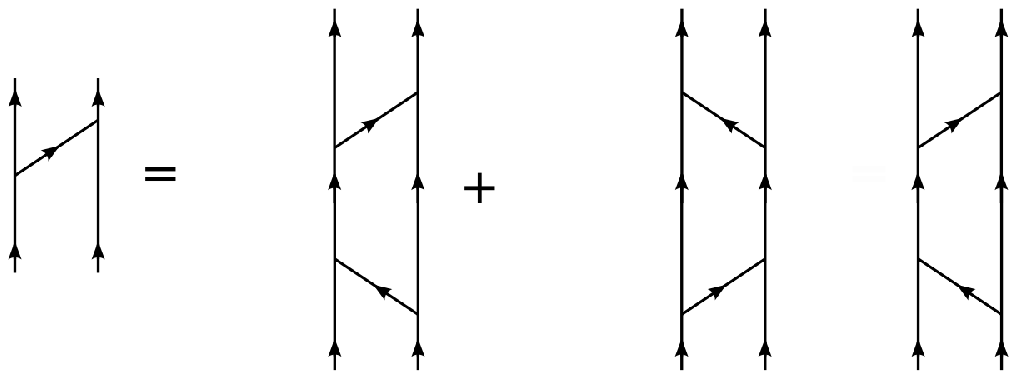}
}. \end{equation}
That is, we assume the rung on the LHS is outward, and $a < b$. We apply \eqref{reductionstep} upside-down to get the middle term, with some coefficients which we ignore. If we work modulo lower terms, i.e. we ignore ladders with an inward rung on top, then we obtain the RHS, which has only the $t=0$ term.

Suppose that $y_2 \ge x_2$ in $E$. Using \eqref{seesaw} on the first rung, and ignoring the terms which vanish modulo $I_{< \l}$, we obtain the LHS below.
\begin{equation} \label{toprungseesaw} {
\labellist
\tiny\hair 2pt
 \pinlabel {$x_1$} [ ] at 17 41
 \pinlabel {$x_2$} [ ] at 41 41
 \pinlabel {$y_2$} [ ] at 6 86
 \pinlabel {$y_1$} [ ] at 30 86
 \pinlabel {$y_1$} [ ] at 12 125
 \pinlabel {$y_2$} [ ] at 35 125
 \pinlabel {$y_3$} [ ] at 59 109
 \pinlabel {$x_1$} [ ] at 158 37
 \pinlabel {$x_2$} [ ] at 184 37
 \pinlabel {$y_2$} [ ] at 145 76
 \pinlabel {$y_1$} [ ] at 170 76
 \pinlabel {$y_1$} [ ] at 151 125
 \pinlabel {$y_2$} [ ] at 176 125
 \pinlabel {$y_3$} [ ] at 199 109
\endlabellist
\centering
\ig{1}{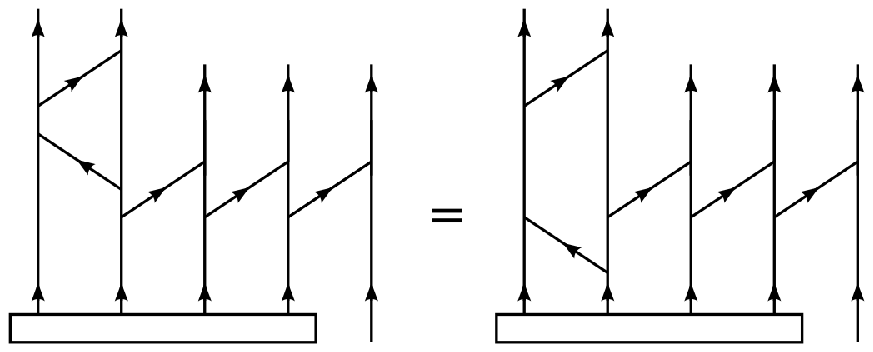}
} \end{equation}
We apply associativity \eqref{associativity} to obtain the RHS. In the RHS, the topmost rung is neutral (absorbing into $N$), and the bottommost rung is outward or neutral (absorbing into $X$) because $y_2 > x_1$ and $y_2 \ge x_2$. Also, because $y_1$ is minimal amongst all vertical segment labels, it is easy to confirm that the remaining rungs are all outward or neutral.  The remaining rungs are again of the form \eqref{lfredux} but with fewer uprights, and induction handles the rest.

Thus we can assume that $y_2 < x_2$. Suppose then that $x_2 = y_3$ and the second rung is neutral. Let us apply \eqref{neutralundo} to the identity of the input strands $(x_1, x_2)$ in \eqref{lfredux}. 
\begin{equation} \label{mixingitup}
{
\labellist
\tiny\hair 2pt
 \pinlabel {$x_1$} [ ] at 6 35
 \pinlabel {$x_2$} [ ] at 44 35
 \pinlabel {$x_2$} [ ] at 6 63
 \pinlabel {$x_1$} [ ] at 44 63
 \pinlabel {$x_1$} [ ] at 6 80
 \pinlabel {$x_2$} [ ] at 44 85
 \pinlabel {$y_1$} [ ] at 11 125
 \pinlabel {$y_2$} [ ] at 36 125
 \pinlabel {$y_3$} [ ] at 59 125
 
 \pinlabel {$x_1$} [ ] at 150 35
 \pinlabel {$x_2$} [ ] at 187 35
 \pinlabel {$d$} [ ] at 150 63
 \pinlabel {$e$} [ ] at 187 63
 \pinlabel {$x_1$} [ ] at 150 84
 \pinlabel {$x_2$} [ ] at 187 84
 \pinlabel {$y_1$} [ ] at 156 125
 \pinlabel {$y_2$} [ ] at 179 125
 \pinlabel {$y_3$} [ ] at 204 125
  \pinlabel {\small $\displaystyle \sum_{r>0}$} [ ] at 136 62
  \pinlabel {$\a_3$} [ ] at 54 78
  \pinlabel {$\a_3$} [ ] at 200 78
\endlabellist
\centering
\ig{1}{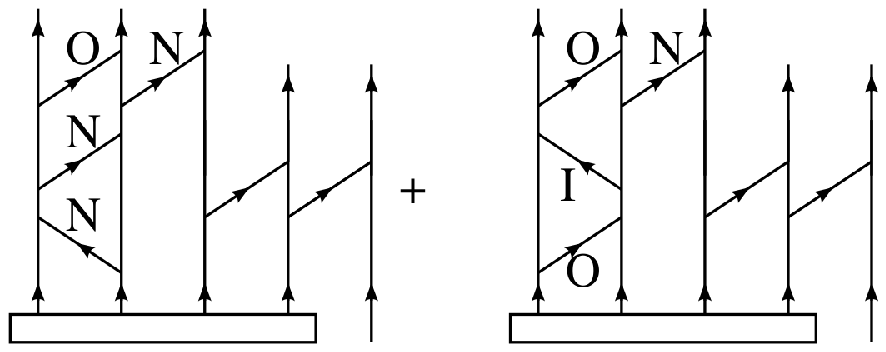}
} \end{equation}
In this diagram, $d = x_1 - r$ and $e = x_2 + r$. We have labeled the rungs I, N, O for inward, neutral, outward, and we have ignored coefficients.
In each diagram, the bottommost rung is neutral or outward, and will absorb into $X$, so we ignore it. We claim that each of the diagrams with $r>0$ is actually inside $I_{< \l}$. After all, $e$ is larger than $y_3=x_2$, which is larger than $y_1$ and $y_2$ and $d$ and $\a_3$; no sequence of outward and neutral ladders could transform $(d,e,\a_3)$ into $(y_1, y_2, y_3)$ by Lemma \ref{lem:maximal}. Thus we need only consider the first term in \eqref{mixingitup}.

Now we can apply ``rung Reidemeister III'' to the top three rungs. We rewrite \eqref{R3} in the version we will use.
\begin{equation} \label{R3special} {
\labellist
\tiny\hair 2pt
 \pinlabel {$a$} [ ] at 7 7
 \pinlabel {$b$} [ ] at 31 7
 \pinlabel {$c$} [ ] at 55 7
 \pinlabel {$b$} [ ] at 13 49
 \pinlabel {$c$} [ ] at 28 63
 \pinlabel {$p$} [ ] at 7 92
 \pinlabel {$q$} [ ] at 32 93
 \pinlabel {$a$} [ ] at 55 93
 \pinlabel {$a$} [ ] at 84 7
 \pinlabel {$b$} [ ] at 107 7
 \pinlabel {$c$} [ ] at 133 7
 \pinlabel {$p$} [ ] at 103 32
 \pinlabel {$q$} [ ] at 127 55
 \pinlabel {$p$} [ ] at 83 93
 \pinlabel {$q$} [ ] at 108 92
 \pinlabel {$a$} [ ] at 132 92
 \pinlabel {$a$} [ ] at 36 44
 \pinlabel {$a$} [ ] at 112 57
 \pinlabel {\small $\sum$} [ ] at 163 51
 \pinlabel {$a$} [ ] at 176 7
 \pinlabel {$b$} [ ] at 199 7
 \pinlabel {$c$} [ ] at 222 7
 \pinlabel {$p$} [ ] at 175 93
 \pinlabel {$q$} [ ] at 200 93
 \pinlabel {$a$} [ ] at 223 93
\endlabellist
\centering
\ig{1}{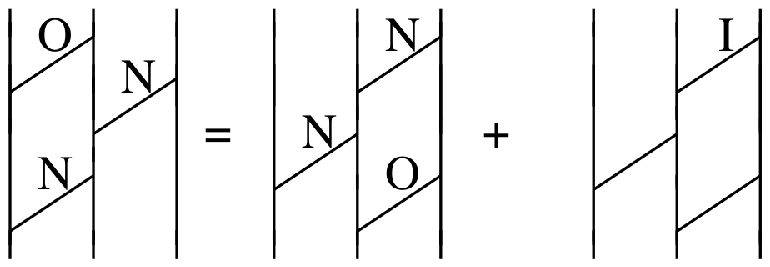}
} \end{equation}
(The first term has coefficient $1$, which will be relevant in Proposition \ref{gammaEasyprop}; we ignore the remaining coefficients.) So, applying \eqref{R3special} to
the first term in \eqref{mixingitup}, and ignoring the neutral rungs on top and bottom, and ignoring terms with an inward rung on top, we end up with another ladder of the form \eqref{lfredux}, with fewer uprights. Induction handles
this case.

Thus we can assume that $x_2 < y_3$. Suppose next that $y_3 \ge x_3$. Now we apply \eqref{seesaw} to the second rung, ignoring coefficients as usual.
\begin{equation} \label{secondrungseesaw} {
\labellist
\tiny\hair 2pt
 \pinlabel {$x_1$} [ ] at 17 27
 \pinlabel {$x_2$} [ ] at 41 27
 \pinlabel {$x_3$} [ ] at 64 27
 \pinlabel {$y_1$} [ ] at 11 112
 \pinlabel {$y_2$} [ ] at 36 112
 \pinlabel {$y_3$} [ ] at 60 112
 \pinlabel {$x_1$} [ ] at 158 27
 \pinlabel {$x_2$} [ ] at 183 27
 \pinlabel {$x_3$} [ ] at 206 27
 \pinlabel {$y_1$} [ ] at 151 112
 \pinlabel {$y_2$} [ ] at 177 112
 \pinlabel {$y_3$} [ ] at 200 112
 \pinlabel {$x_1$} [ ] at 298 27
 \pinlabel {$x_2$} [ ] at 322 27
 \pinlabel {$x_3$} [ ] at 346 27
 \pinlabel {$y_1$} [ ] at 292 112
 \pinlabel {$y_2$} [ ] at 316 112
 \pinlabel {$y_3$} [ ] at 340 112
 \pinlabel {$y_3$} [ ] at 31 74
 \pinlabel {$\a_2$} [ ] at 55 74
 \pinlabel {$c$} [ ] at 172 74
 \pinlabel {$d$} [ ] at 196 74
 \pinlabel {$e$} [ ] at 311 68
 \pinlabel {$f$} [ ] at 345 66
 \pinlabel {\small $\sum$} [ ] at 135 58
 \pinlabel {\small $\sum$} [ ] at 275 58
\endlabellist
\centering
\ig{1}{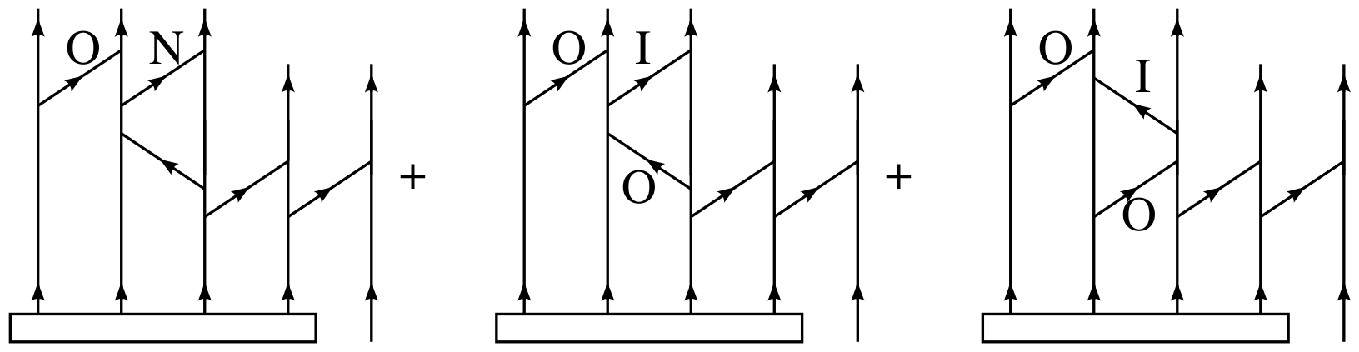}
} \end{equation}
Once again, any term with an inward rung is actually in $I_{< \l}$. This can be seen by observing that both $c$ and $f$ are larger than $y_1$, $y_2$, $y_3$, $x_1$, and their counterpart $d$ or $e$, and using Lemma \ref{lem:maximal}. The remaining term has the form
\begin{equation} \label{weirdo} {
\labellist
\tiny\hair 2pt
 \pinlabel {$x_1$} [ ] at 19 25
 \pinlabel {$x_2$} [ ] at 42 25
 \pinlabel {$x_3$} [ ] at 66 25
 \pinlabel {$y_1$} [ ] at 10 114
 \pinlabel {$y_2$} [ ] at 35 114
 \pinlabel {$y_3$} [ ] at 60 114
 \pinlabel {$y_3$} [ ] at 30 73
 \pinlabel {$x_1$} [ ] at 156 25
 \pinlabel {$x_2$} [ ] at 179 25
 \pinlabel {$x_3$} [ ] at 202 25
 \pinlabel {$y_1$} [ ] at 146 114
 \pinlabel {$y_2$} [ ] at 173 114
 \pinlabel {$y_3$} [ ] at 195 114
 \pinlabel {$y_3$} [ ] at 167 59
\endlabellist
\centering
\ig{1}{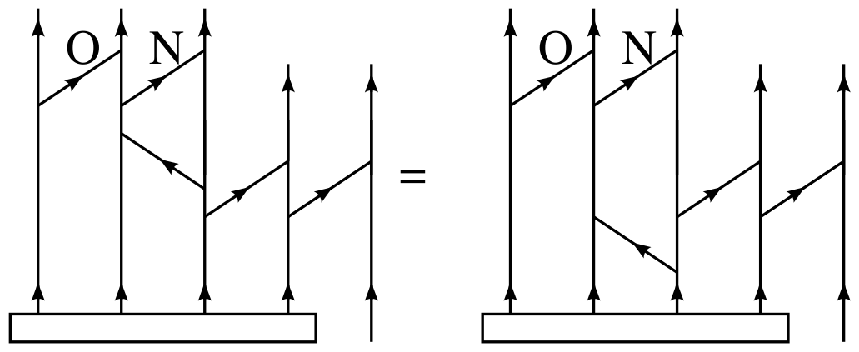}
} \end{equation}
Once again, the bottommost rung in the new ladder is outward or neutral, because we assumed $y_3 \ge x_3$. The remaining rungs are outward or neutral, but the second rung is neutral, so this ladder can be dealt with as in the $x_2 = y_3$ case above.
 
Similarly, if $x_3 = y_4$, we apply \eqref{neutralundo} to the identity of $(x_2, x_3)$. We obtain
\begin{equation} \label{mixingitup2} {
\labellist
\tiny\hair 2pt
 \pinlabel {\small $\sum$} [ ] at 138 65
 \pinlabel {$d$} [ ] at 209 62
 \pinlabel {$\a_4$} [ ] at 223 76
 \pinlabel {$c$} [ ] at 174 62
\endlabellist
\centering
\ig{1}{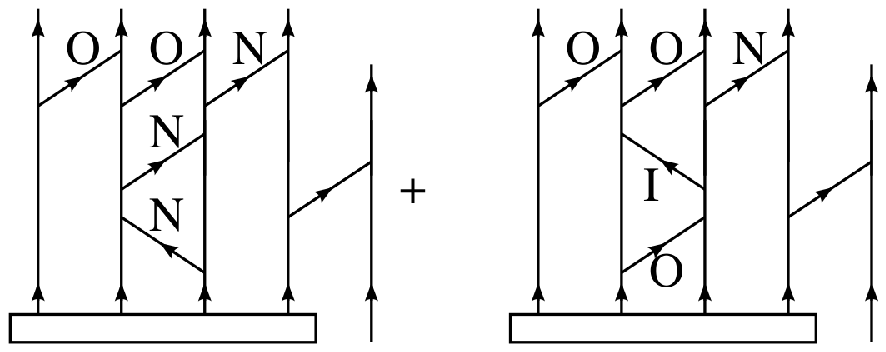}
} \end{equation}
Once again, any diagram with an inward rung is in $I_{< \l}$ using Lemma \ref{lem:maximal}, because $d$ is larger than $y_i$ for $i \le 4$, and larger than $x_1$, $c$, or $\a_4$. Meanwhile, the remaining term can be dealt with using \eqref{R3special}, reducing it to a previous case.

Similar arguments handle every possible problem inequality, either $y_i \ge x_i$ or $x_i = y_{i+1}$. If there are no problem inequalities, then the original ladder was already an
elementary light ladder, as desired. This concludes the proof.

Note that there is no inequality relating $y_{k+1}$ and $a$. In particular, the simplification algorithm proposed by this proof will never introduce an outward or neutral rung on the
bottom involving the last upright. \end{proof}

\subsection{Double ladders span: part II}
\label{subsec-lightladdersspan2}

Now we continue our proof of the Light Ladder Lemma, using Proposition \ref{thecrux} as a black box. The proofs of the strong versions of the Light Ladder Lemma and Double Ladder Theorem are very similar to the proofs of the weak versions, but slightly more technical in a distracting way. We concern ourselves with the weak versions first.

Fix a dominant weight $\l$. For an object $\un{w} \in \Fund^+$, let $A(\un{w},\l)$ be the statement that all ladders from $\un{w}$ to a sequence $\un{z} \in P(\l)$ are in the (weak)
span of light ladders modulo $I_{< \l}$. Then $A(\un{w}, \l)$ is just a restatement of the weak Light Ladder Lemma for that particular $\un{w}$ and $\l$. Note that terms in $I_{< \l}$
can be expressed as IO ladders where the inward part is strict. Let $B(\un{w}, \l)$ be the same statement, with the additional requirement that the lower terms in $I_{< \l}$ can be
expressed as IO ladders where the outward part is actually a light ladder.

\begin{lemma}$A(\un{w}, \mu)$ for all $\mu \le \l$ implies $B(\un{w}, \mu)$ for all $\mu \le \l$. \end{lemma}

\begin{proof} For the base case, $A(\un{w}, 0)$ and $B(\un{w}, 0)$ are equivalent, since there are no lower terms. So, let us assume $B(\un{w}, \mu)$ for $\mu < \l$ and $A(\un{w}, \l)$, and prove $B(\un{w}, \l)$.
	
Given an arbitrary morphism $\un{w} \to \un{z}$, $A(\un{w}, \l)$ lets us write it as a linear combination of light ladders and lower terms, which are IO ladders $\a \circ \b$ with a
strict inward part $\a$. Any lower term $\a \circ \b$ factors through some weight $\mu < \l$. By $B(\un{w}, \mu)$, we can rewrite $\b$ as a linear combination of light ladders $\g$,
plus lower terms $\d \circ \e$ where $\e$ is a light ladder. Thus both $\a \circ \g$ and $\a \circ \d \circ \e$ are IO ladders where the outward part is a light ladder. Hence $B(\un{w},
\l)$ is proven. \end{proof}
 
Let $A(\nu, \l)$ be the statement that $A(\un{w}, \l)$ holds for all $\un{w} \in P(\nu)$. Let $A(m, \l)$ be the statement that $A(\un{w}, \l)$ holds for all $\un{w}$ with width $\le m$.
Note that many of these statements are trivial because no outward ladders are possible, such as when if $\nu \nleq \l$, or when the width of $\l$ is larger than $m$. Let $A(m)$ denote
$A(m, \l)$ for all $\l$; this is really only a statement about $\l$ with width $\le m$. Clearly $A(m)$ implies the corresponding statement $B(m)$. Note that $A(0)$ and $A(1)$ are trivial, and $A(2)$ follows directly from Proposition \ref{thecrux}.

\begin{lemma} \label{mtomp1} $A(m)$ implies $A(m+1)$. \end{lemma}

\begin{proof} Let $\un{w}$ have weight $\mu$ and width $m$, and consider the sequence $\un{w}a$ which adjoins a new strand labeled $a$. We wish to deduce $A(\un{w}a, \l)$ for any given
$\l$. We inductively assume $A(m)$ and $B(m)$, and $A(\un{w}' a, \l)$ for any $\un{w}' \in P(\rho)$, $\rho < \mu$. Consider any ladder from $\un{w}a$ to $\un{z}$ where $\un{z} \in
P(\l)$. By Proposition \ref{thecrux}, we can assume that our ladder has the form \eqref{eq:thecrux}. Using $B(m)$, we can assume that $X$ is either a light ladder, or has a strict IO
decomposition with its outward part being a light ladder.

\begin{equation} \label{inductiveoptions} {
\labellist
\small\hair 2pt
 \pinlabel {$\mu + \w_a$} [ ] at -7 30
 \pinlabel {$LL$} [ ] at 42 48
 \pinlabel {$\nu + \w_a$} [ ] at 3 62
 \pinlabel {$E$} [ ] at 58 75
 \pinlabel {$\l$} [ ] at 29 91
 \pinlabel {$N$} [ ] at 58 103
 \pinlabel {or} [ ] at 124 71
 \pinlabel {$\mu + \w_a$} [ ] at 125 4
 \pinlabel {$LL$} [ ] at 179 24
 \pinlabel {$\t + \w_a$} [ ] at 152 41
 \pinlabel {$I$} [ ] at 182 48
 \pinlabel {$\nu + \w_a$} [ ] at 140 61
 \pinlabel {$E$} [ ] at 195 75
 \pinlabel {$\l$} [ ] at 162 90
 \pinlabel {$N$} [ ] at 195 102
\endlabellist
\centering
\ig{1}{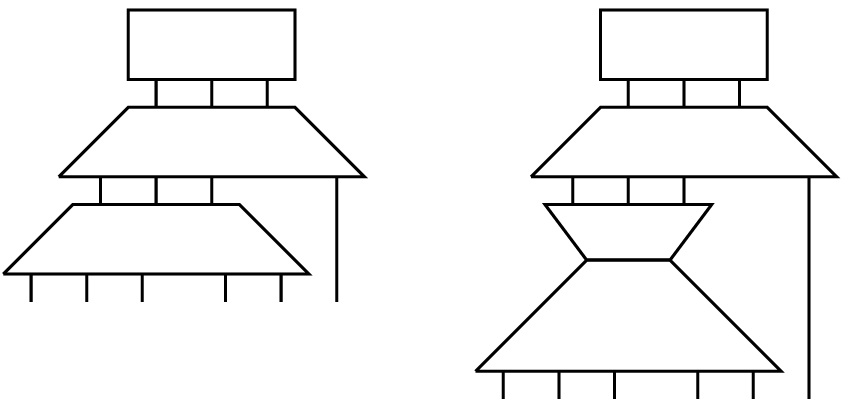}
} \end{equation}

By the algorithmic construction, the first diagram above is a light ladder. Meanwhile, for the second diagram above, $\t < \mu$ so we may inductively simplify the morphism from $\t + \w_a$ to $\l$, obtaining the following schema.
\begin{equation} \label{inductiveproblem} {
\labellist
\small\hair 2pt
 \pinlabel {$\mu + \w_a$} [ ] at -7 3
 \pinlabel {$LL$} [ ] at 42 19
 \pinlabel {$\t + \w_a$} [ ] at 6 33
 \pinlabel {$LL$} [ ] at 59 46
 \pinlabel {$\l$} [ ] at 28 62
\endlabellist
\centering
\ig{1}{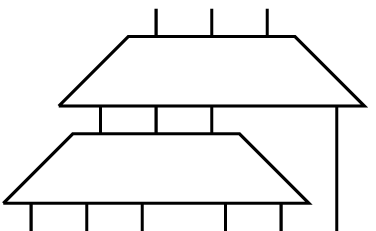}
} \end{equation}
The problem is that the upper light ladder may have multiple tiers, unlike the first diagram of \eqref{inductiveoptions} where it was a single tier. Thus we can not yet conclude that such a diagram is in the span of light ladders. However, we may resolve it with the same style of induction.

Let $C(\mu, \t, a, \l)$ denote the statement that any diagram of the form \eqref{inductiveproblem} is in the span of light ladders, modulo $I_{< \l}$. In such a diagram, the upper light
ladder has multiple tiers; let $\s + \w_a$ be the weight at the source of the last tier $T$. Then again our diagram is exactly of the form \eqref{eq:thecrux}, where $E$ and $N$ form the
last tier $T$, and $X$ is the rest (a composition of a light ladder with most of another light ladder). So we can use $B(m)$ again to rewrite it in terms of the diagrams in
\eqref{inductiveoptions}, but replacing $\nu$ with $\s$, and $\t$ with some weight $\rho < \s$.

Thus our problem alternates between needing to reduce the second diagram of \eqref{inductiveoptions}, and needing to reduce \eqref{inductiveproblem}. However, each time we reach
\eqref{inductiveproblem}, the weight $\t$ is lowered. Hence, $C(\mu, \t, a, \l)$ reduced to $C(\mu, \rho, a, \l)$ for some $\rho < \t$. Of course, if $\rho = \l$ or $\rho < \l$ then the
statement is trivial. Induction yields $C(\mu, \t, a, \l)$ for all $\t$, which implies $A(\un{w} a, \l)$ as above. \end{proof}

Thus, we have proven $A(\un{w}, \l)$ for all $\un{w}$ and $\l$, which is the Weak Light Ladder Lemma.

\begin{lemma} \label{weaktostrong} The Weak Light Ladder Lemma implies the Weak Double Ladder Theorem. \end{lemma} 
	
\begin{proof} Let $D(\un{w}, \un{x}, \l)$ be the statement that double ladders of the form $\overline{L}_{\fb} \circ L_\eb$ for $\eb \in E(\un{w}, \l)$ and $\fb \in E(\un{x}, \l)$ span
the part of $\Hom(\un{w}, \un{x})$ which is in the ideal $I_{\le \l}$, modulo the ideal $I_{< \l}$. If we can prove $D(\un{w}, \un{x}, \l)$ for all $\l$ less than the (minimum of the)
weights of $\un{w}$ and $\un{x}$, then one knows that the entire Hom space $\Hom(\un{w}, \un{x})$ is spanned by double ladders, by an obvious induction over the filtration. Proving $D$
for all cases then gives the Weak Double Ladder Lemma.

A morphism in $\Hom(\un{w}, \un{x})$ in the ideal $I_{\le \l}$ is by definition a morphism $\un{w} \to \un{y}$ composed with a morphism $\un{y} \to \un{x}$, where $\un{y} \in P(\mu)$
for $\mu \le \l$, and linear combinations of these. Modulo $I_{< \l}$, we may assume that $\un{y} \in P(\l)$. Applying $A(\un{w}, \l)$, we can assume that the map $\un{w} \to \un{y}$ is
a light ladder. Applying $A(\un{x}, \l)$, we may assume that the map $\un{y} \to \un{x}$ is an upside-down light ladder. Thus we have proven $D(\un{w}, \un{x}, \l)$. \end{proof}

Now, let us deduce the Strong Light Ladder Lemma and the Strong Double Ladder Theorem from their weak counterparts. The issue is that we need to rigidify our choice of neutral ladders.
Thankfully, the have the Neutral Ladder Lemma \ref{lem:neutralladderlemma}, which says that any two choices of neutral ladders differs by lower terms.

\begin{proof}[Proof of the Strong Light Ladder Lemma \ref{lem:lightladdersspan}] Let $A'(m)$ denote the statement that, for any $\un{w}$ of width $\le m$ and any $\un{z} \in P(\l)$, the
morphism space $\Hom(\un{w}, \un{z})$ is spanned by (chosen) light ladders $LL_{\eb} \co \un{w} \to \un{x}_\l$ for $\eb \in E(\un{w}, \l)$, composed with a fixed neutral ladder $N \co
\un{x}_\l \to \un{z}$. Once again, $A'(0)$ and $A'(1)$ are both trivial.

We use $A'(m)$ to prove $A'(m+1)$ in exactly the same way as Lemma \ref{mtomp1}. When reducing an arbitrary map $\un{w} \to \un{z}$, we can use Proposition \ref{thecrux}, and then
simplify the diagram $X$ using $A'(m)$, to obtain a linear combination of diagrams as in \eqref{inductiveproblem}. This time, we can assume that the $LL$ maps which appear are the
chosen ones, thanks to $A'(m)$ and $B'(m)$.

In the proof of Lemma \ref{mtomp1}, the first diagram in \eqref{inductiveproblem} was already a light ladder. Now it may differ from our chosen light ladder in the choice of concluding
neutral ladder. The difference between these two neutral ladders factors through lower terms, which shows that this first diagram is not a problem.

(Technically, the first diagram in \eqref{inductiveproblem} may also differ in a neutral ladder before the elementary light ladder $E$. But because $A'(m)$ allows for any fixed neutral
ladder to be chosen in our spanning set, we can choose the neutral ladder as desired.)

To deal with the second diagram in \eqref{inductiveproblem}, we reduce to $C'(\mu, \t, a, \l)$, exactly as in the proof of Lemma \ref{mtomp1}. The rest of the argument works verbatim,
replacing arbitrary light ladders with chosen light ladders. \end{proof}

\begin{proof}[Proof of the Strong Double Ladder Theorem] Let $D'(\un{w}, \un{x}, \l)$ be the statement that double ladders $\LL^\l_{\eb, \fb}$ (see Definition \ref{defn:doubleladders})
span $\Hom(\un{w}, \un{x})$ inside $I_{\le \l}$, modulo $I_{< \l}$. Once again, if we can prove $D'(\un{w}, \un{x}, \l)$ for all inputs, one obtains the Strong Double Ladder Theorem by
an obvious induction over the filtration.

A morphism in $\Hom(\un{w}, \un{x})$ in the ideal $I_{\le \l}$ is a composition of a map $\un{w} \to \un{y}$ with a map $\un{y} \to \un{x}$, for $\un{y} \in P(\mu)$, $\mu \le \l$, and
linear combinations of these. Again, modulo $I_{< \l}$ we can assume that $\un{y} \in P(\l)$, and using $A'(\un{w}, \l)$ and $A'(\un{x}, \l)$ we can assume that our morphism is the
composition \[\overline{LL}_{\fb} \circ \overline{N} \circ N \circ LL_{\eb}.\] Here, $N$ is a fixed neutral ladder from $\un{x}_\l$ to $\un{y}$. By the Neutral Ladder Lemma,
$\overline{N} \circ N$ can be replaced with the identity of $\un{x}_\l$ modulo $I_{< \l}$, which proves $D'(\un{w}, \un{x}, \l)$. \end{proof}

This concludes the proof that double ladders span. In the next section we prove that double ladders are linearly independent, and are thus a basis for morphisms in $\Fund^+$.

We encourage the reader to look at \cite[\S 2]{ELauda} for the definition of a strictly object-adapted cellular category (SOACC), and to verify that $\Fund^+$ satisfies the conditions of
\cite[Lemma 2.8]{ELauda}. Therefore, $\Fund^+$ is an SOACC, for which double ladders form a cellular basis. Since SOACCs are closed under categorical equivalence, $\Fund$ is also an SOACC.

\subsection{Light ladders are linearly independent}
\label{subsec-lightladdersindep}

Suppose that $\Fund$ is defined with base ring $\Bbbk = \ZM[\d]$. Let us recall the evaluation functor $\Gamma_n$ defined in \cite[\S 3]{CKM}. It is a functor from $\Fund \ot_{\ZM[\d]}
\ZM[q,q^{-1}]$ to representations of Lusztig's $\ZM[q,q^{-1}]$ integral form of the quantum group $U_q(\sl_n)$. In \cite{CKM} the functor is only defined over the base ring $\CM(q)$. It
is clear, however, that their definition makes sense over $\ZM[q,q^{-1}]$; what is no longer clear is whether $\G_n$ is fully faithful.

Let $V_{\w_a}$ denote the free $\ZM[q,q^{-1}]$-module with basis $x_S$ given by subsets $S \subset \{1, 2, \ldots, n\}$ of size $a$. The integral form of $U_q(\sl_n)$ acts on
$V_{\w_a}$ by the usual formulas (not relevant for us, see \cite[\S 3]{CKM}).

Given two disjoint subsets $S,T \subset \{1, \ldots, n\}$, let $\ell(S,T)$ denote the number of pairs $i < j$ with $i \in S$ and $j \in T$. We now define some $\ZM[q,q^{-1}]$-linear maps
between tensor products of various $V_{\w_a}$. The multiplication map $M \co V_{\w_a} \ot V_{\w_b} \to V_{\w_{a+b}}$ is defined by
\begin{equation} \label{mergedef} M(x_S \ot x_T) = \begin{cases} (-q)^{\ell(S,T)} x_{S \cup T} & \textrm{if } S \cap T = \emptyset \\ 0 & \textrm{else.} \end{cases} \end{equation}
The comultiplication map $M' \co V_{\w_{a+b}} \to V_{\w_a} \ot V_{\w_b}$ is given by
\begin{equation} \label{splitdef} M'(x_S) = (-1)^{ab} \sum_{T \subset S} (-q)^{-\ell(S \setminus T,T)} x_T \ot x_{S \setminus T}. \end{equation}
Finally, the duality map $D \co V_{\w_a} \to V_{\w_{n-a}}^*$ is given by
\begin{equation} \label{dualdef} D(x_S)(x_T) = \begin{cases} (-q)^{\ell(S,T)} & \textrm{if } S \cap T = \emptyset \\ 0 & \textrm{else.} \end{cases} \end{equation}
These maps are the images under $\Gamma_n$ of the merging trivalent vertex, the splitting trivalent vertex, and the tag respectively. It was proven in \cite[\S 3]{CKM} that these maps are 
$U_q(\sl_n)$-intertwiners, and that they satisfy the $\sl_n$-web relations.

The signs and powers of $q$ in the formulas above will be entirely irrelevant for our results below. Suffice it to say that these coefficients are invertible.

\begin{remark} \label{howtocompute} To compute how a web in $\Fund^+$ acts, assign to the $i$-th input (resp. output) strand a subset $S_i \subset \{1, \ldots, n\}$ (resp. $T_i$), with
size equal to the label on that strand. The coefficient of $x_{T_1} \ot x_{T_2} \ot \cdots$ in the web applied to $x_{S_1} \ot x_{S_2} \ot \cdots$ is a sum of invertible scalars. Each
term comes from a labelling of all the strands in the web by subsets of $\{1, \ldots, n\}$ of the appropriate size, such that for each trivalent vertex, the bigger set is the disjoint
union of the two smaller sets. We call this a \emph{consistent labeling} by subsets. In the examples we compute below, there will either be one or zero consistent labelings, so we need
not worry about any cancellation which may occur when adding invertible scalars. \end{remark}

To a subset $S \subset \{1, \ldots, n\}$ of size $a$, we associate the obvious 01-sequence which records the entries of $S$, and we think of this sequence as a weight $\mu$ in
$V_{\w_a}$. Thus we also write $x_\mu$ for the basis element $x_S$. The tensor product $V_{\un{w}}$ has a basis given by $x_{\un{w}, \eb} = x_{S_1} \ot x_{S_2} \ot \cdots \ot x_{S_d}$
for various weight subsequences $\eb \subset \un{w}$ (not just dominant ones), in the obvious way. We will be primarily interested in the \emph{dominant subbasis} given by $x_{\un{w},
\eb}$ for $\eb \in E(\un{w})$ a dominant weight subsequence. Clearly, the linear independence of various maps between tensor products can be checked by their linear independence on
their respective dominant subbases.

Let $\top$ denote the highest weight in $V_{\w_a}$, the 01-sequence $(11111100000)$, so that $x_{\top}$ is a basis element. We also let $\top$ denote the full weight subsequence of
$\un{w}$, so that $x_{\un{w}, \top}$ denotes the tensor product $x_{\top} \ot x_{\top} \ot \cdots \ot x_{\top}$. Corresponding to the weight $\top$ is the subset $[1,a] = \{1, \ldots, a\} \subset [1,n]$, where we use interval notation.

\begin{lemma} \label{lem:whatNdoes} A neutral ladder from $\un{w}$ to $\un{y}$ will send $x_{\un{w},\top}$ to $x_{\un{y},\top}$ with invertible coefficient. Moreover, $x_{\un{w},\top}$ is the only basis element which, after applying a neutral ladder, has a nonzero coefficient for $x_{\un{y}, \top}$. (Note that ``basis element'' refers to any $x_{\un{w}, \eb}$, not just for $\eb$ dominant.) \end{lemma}

\begin{proof} It is enough to check this for a neutral rung. Label each strand in a neutral rung: $S_1$ and $S_2$ for the inputs, $T_1$ and $T_2$ for the outputs, and $C$ for the
crossbar. We assume that $T_1 = [1,a]$ and $T_2 = [1,b]$ are the highest weights, so that $(T_1, T_2) = (\top,\top)$. Without loss of generality, assume that $a < b$. Then $C \subset
T_2$ has size $b-a$, and is disjoint from $T_1$, so there is only one possibility, $C = [a+1,b]$. Then $S_1 = C \coprod T_1 = T_2$ and $S_2 = T_2 \setminus C = T_1$. Thus if $(T_1, T_2)
= (\top,\top)$ then in a consistent labeling by subsets, we must also have $(S_1, S_2) = (\top,\top)$. This is an ``if and only if,'' by the same argument run in reverse. This implies
the desired result, by Remark \ref{howtocompute}. \end{proof}

\begin{lemma} \label{lem:whatLdoes} Let $L_\mu$ be an elementary light ladder from $\un{w} \ot a$ to $\un{y}$, associated to a weight $\mu$ in the fundamental representation $V_{\w_a}$. Then one has
	\begin{equation} \label{whatLdoes} L_\mu(x_{\un{w}, \top} \ot x_\nu) = \begin{cases}
		0 & \textrm{if } \nu > \mu, \\
		\xi x_{\un{y}, \top} & \textrm{if } \nu = \mu,\\
		\textrm{lower terms} & \textrm{if } \nu \ngeq \mu. \end{cases} \end{equation}
Here, $\xi$ is some invertible coefficient. Lower terms refers to elements of $V_{\un{y}}$ with zero coefficient for $x_{\un{y}, \top}$. \end{lemma}

\begin{proof} First consider the case of a single rung, a diagram as in \eqref{loweringform} except with $k=1$. Following Remark \ref{howtocompute}, let us attempt to label the strands
of our diagram with subsets of $[1,n]$. Let $A, X_1, Y_1, Y_2$ denote the sets on the various inputs and outputs, and let $C$ denote the set on the crossbar.
\begin{equation} {
\labellist
\small\hair 2pt
 \pinlabel {$X_1$} [ ] at 8 29
 \pinlabel {$A$} [ ] at 40 7
 \pinlabel {$Y_1$} [ ] at 16 99
 \pinlabel {$Y_2$} [ ] at 40 99
 \pinlabel {$C$} [ ] at 28 52
\endlabellist
\centering
\ig{1}{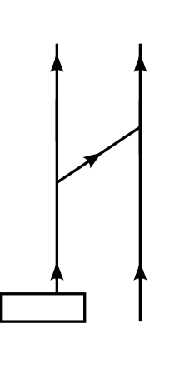}
} \end{equation}
Let $M = [1,y_1] \cup [x_1+1,y_2]$. This is the set of size $a$ which corresponds to the weight $\mu$ in $V_{\w_a}$ which gives rise to this elementary light ladder. For a set $A$ to be less than $M$ in the dominance order is equivalent to the statement that $A \cap [1,x_1]$ has size $\le y_1$. (Note: everything makes sense when $y_1 = 0$, interpreting the interval $[1,y_1]$ as the empty set. It also makes sense when $y_2 = n$.)
	
We need to check three things. \begin{itemize}
	\item If $X_1 = [1,x_1]$ and $A \cap [1,x_1]$ has size $>y_1$, then there are no consistent labelings for any $Y_1, Y_2$.
	\item If $X_1 = [1,x_1]$ and $A = M$, then there is one consistent labeling, and it satisfies $Y_1 = [1,y_1]$ and $Y_2 = [1,y_2]$.
	\item If $X_1 = [1,x_1]$ and $Y_1 = [1,y_1]$ and $Y_2 = [1,y_2]$ in a consistent labeling, then $A=M$.
\end{itemize}
Together, the second and third point say that when $X_1 = [1,x_1]$, then $(Y_1, Y_2) = (\top,\top)$ if and only if $A = M$, in which case there is a unique consistent labeling (i.e. $C$ is determined).

Observe that $Y_1$ contains $A \cap X_1$, since it contains every element of $X_1$ not in $C$, and $C$ and $X_1$ are disjoint. Thus in a consistent labeling, $y_1$ is greater than or
equal to the size of $A \cap X_1$. This confirms the first point. Moreover, if $A \cap X_1$ has size exactly $y_1$, then $Y_1 = A \cap X_1$. From this, it is easy to deduce the second
point. The third point is obvious, since $C$ is forced to be $X_1 \setminus Y_1$ and $A$ is forced to be $Y_2 \setminus C$.

Let us generalize slightly. Fix some particular subset $Z$ inside $[1,x_1]$ of size $y_1$, as an alternative to $[1,y_1]$ for a standard label of $Y_1$. We claim that $(Y_1, Y_2) =
(Z,\top)$ if and only if $A = Z \cup [x_1 + 1, y_2]$. The argument is identical.

Now consider the general case, a diagram as in \eqref{loweringform} with $y_1 < x_1 < \ldots < x_k < y_{k+1}$. Consider a consistent labeling with inputs $X_i$ and $A$, and outputs
$Y_i$. Let $Z_i$ denote the set which labels the interior strand on the $i$-th upright, for $2 \le i \le k$, having size $\a_i$. We set $Z_1 = Y_1$ and $Z_{k+1} = A$. We are interested
in evaluating this elementary light ladder on $x_{\un{w}, \top} \ot x_\nu$, which means we can assume that $X_i = [1,x_i]$ for all $1 \le i \le k$.

By the case of one rung, applied to the leftmost rung, we see that the diagram is inconsistent unless $Z_2 \cap [1,x_1]$ has size $\le y_1$; and that $Y_1 = [1,y_1]$ and $Y_2 = [1,y_2]$
if and only if $Z_2 = [1, y_1] \cup [x_1+1, y_2]$.

Now let $Z \subset [1,y_2]$ be arbitrary of size $\a_2$, and consider the second rung. We see that the labeling is inconsistent unless $Z_3 \cap [1,x_2]$ has size $\le \a_2$, and that
$Z_2 = Z$ and $Y_3 = [1,y_3]$ if and only if $Z_3 = Z \cup [x_2+1,y_3]$.  This argument can be repeated for every rung.

Let $M = [1, y_1] \cup [x_1+1,y_2] \cup [x_2+1,y_3] \cup \cdots \cup [x_k+1,y_{k+1}]$. This is the subset of $[1,n]$ which corresponds to the weight $\mu$ of $V_{\w_a}$ giving rise to
this elementary light ladder. If $A = M$ then $Z_k = M \setminus [x_k+1,y_{k+1}]$, and $Y_{k+1} = [1,y_{k+1}]$ by the analysis of the last rung. By the analysis of the $i$-th rung, each
$Z_i$ is obtained from $Z_{i+1}$ by removing the last interval, and each $Y_{i+1} = [1,y_{i+1}]$. Conversely, if $Y_i = [1,y_i]$ for all $i$, then $Z_2 = [1,y_1] \cup [x_1+1,y_2]$, and
$Z_3 = Z_2 \cup [x_2+1, y_3]$, and so on until $A = M$. Thus $A = M$ if and only if $Y_i = [1,y_i]$ for all $i$, and moreover the labeling is unique.

Meanwhile, considering the last rung, the diagram is inconsistent unless $A \cap [1, x_k]$ has size $\le \a_k$. Moreover, $A \cap [1,x_k] \subset Z_k$, so that $A \cap [1,x_{k-1}] \subset Z_k \cap [1,x_{k-1}]$, which must have size $\le \a_{k-1}$ in a consistent diagram. Repeating this we see that $A \cap [1,x_i]$ has size $\le \a_i$ in a consistent diagram, which is precisely the condition that $A$ is less than or equal to $M$ in the dominance order. \end{proof}

\begin{lemma} \label{lem:whattierdoes} Now let $L_{\mu}$ be the elementary light ladder attached to $\mu$, precomposed with a neutral ladder on $\un{w}$, and postcomposed with a neutral ladder as well. In other words, $L_{\mu}$ is one tier in the algorithmic construction of a light ladder. Then the formula \eqref{whatLdoes} still holds. \end{lemma}
	
\begin{proof} This follows immediately from Lemma \ref{lem:whatLdoes} and Lemma \ref{lem:whatNdoes}. \end{proof}

\begin{defn} Let us equip the set $E(\un{w})$ with its \emph{path dominance order}. We say that $\eb = (e_0, e_1, \ldots, e_d) \ge \fb = (f_0, f_1, \ldots, f_d)$ if $e_i \ge f_i$ for all $i$. Recall from Definition \ref{def:weightsubseq} that $e_i$ records a list of dominant weights, not a list of miniscule weights. Clearly, $\top$ is the maximal element in $E(\un{w})$. \end{defn}

\begin{prop} \label{prop:whatLLdoes} The light ladder $LL_{\un{w}, \eb}$ (with target $\un{y}$) acts on $x_{\un{w}, \fb}$ as follows.
	\begin{equation} \label{whatLLdoes} LL_{\un{w}, \eb}(x_{\un{w}, \fb}) = \begin{cases}
		0 & \textrm{if } \fb > \eb, \\
		\xi x_{\un{y}, \top} & \textrm{if } \fb = \eb,\\
		\textrm{lower terms} & \textrm{if } \fb \ngeq \eb. \end{cases} \end{equation} \end{prop}

\begin{proof} This follows immediately from an iterated application of Lemma \ref{lem:whatLdoes} and Lemma \ref{lem:whatNdoes}, together with the inductive algorithm for constructing $LL_{\un{w}, \eb}$. \end{proof}

Now let us provide the analogous statements for upside-down light ladders. The proofs use the same exact analysis of consistent labelings of an elementary light ladder by subsets, only
interpreted for the upside-down map. We leave this interpretation as an exercise to the reader.

\begin{lemma} \label{lem:whatLbardoes} Let $L_\mu$ be a elementary light ladder in $\Fund^+$ from $\un{w} \ot a$ to $\un{y}$, associated to a weight $\mu$ in the fundamental representation $V_{\w_a}$. Then one has \begin{equation}
\label{whatLbardoes} \overline{L}_\mu(x_{\un{y}, \top}) = \xi x_{\un{w}, \top} \ot x_\mu + \textrm{ lower terms}, \end{equation} where $\xi$ is invertible. Now, lower terms refer to linear combinations of tensors in $V_{\un{w}} \ot V_{\w_a}$ which do not have $x_{\un{w}, \top}$ as their $\un{w}$-component. Moreover, for an arbitrary weight subsequence $\eb \subset \un{y}$, $\overline{L}_\mu(x_{\un{y}, \eb})$ will only have a nonzero coefficient for $x_{\un{w}, \top} \ot x_\nu$ if $\nu \le \mu$, with equality implying that $\eb = \top$. (These elements with $\nu < \mu$ can also be considered as lower terms, if $\mu$ is understood.) \end{lemma}

\begin{prop} \label{prop:whatLLbardoes} The dual light ladder $\overline{LL}_{\un{w}, \eb}$ (with target $\un{y}$) will send $x_{\un{y}, \top}$ to $\xi x_{\un{w}, \eb}$ plus lower
terms. Lower terms are sent to lower terms. \end{prop}

Putting together these two propositions, we have our linear independence theorem.

\begin{thm} \label{linindepthm} Inside $\Hom(\un{w}, \un{z})$, the double ladders $\LL_{\eb, \fb}$ for $\eb \in E(\un{w}, \l)$ and $\fb \in M(\l, \un{z})$ for various $\l$ are all
linearly independent. \end{thm}

\begin{proof} Applying the propositions above, we see that $\LL_{\eb, \fb}$ sends $x_{\un{w}, \eb}$ to $x_{\un{z}, \fb}$ with invertible coefficient plus lower terms, acts as zero on
$x_{\un{w}, \gb}$ for $\gb > \eb$, and otherwise has image consisting of lower terms. Here, lower terms are spanned by $x_{\un{z}, \hb}$ for $\hb < \fb$, together with parts of the
basis which are not in our subbasis. Thus, by a straightforward upper-triangularity argument, these maps must be linearly independent. \end{proof}

We already know that double ladders span $\Fund^+$, so that this theorem implies that $\Gamma_n$ is faithful (on $\Fund^+$, which then by isomorphism implies it is faithful on
$\Fund$). Another way of stating this equivalence is that all $U_q(\sl_n)$ intertwiners are determined by their action on the dominant subbasis.

In fact, our theorem implies that $\Gamma_n$ is faithful after specialization from $\ZM[q,q^{-1}]$ to any field $\Bbbk$, because this does not change the invertibility or the upper
triangularity argument! Let us use this to sketch a proof that $\Fund_{\Bbbk}$ is actually an integrable form for quantum group representations, a fact which has not yet appeared in the
$\sl_n$-web literature. The proof used here is an adaptation of a proof of Webster for the $\sl_2$ case, found in the appendix to \cite{ELib}.

\begin{thm} \label{EquivThm} Let $\Bbbk$ be any $\ZM[q,q^{-1}]$-algebra. Let $U_q(\sl_n) \ot \Bbbk$ denote the $\Bbbk$-integral form of quantum $\sl_n$, and let $\CC$ denote the full
subcategory of its module category whose objects are tensor products of $V_{\w_a} \ot \Bbbk$. Then $\Gamma_n \co \Fund_{\Bbbk} \to \CC$ is an equivalence of categories. \end{thm}

\begin{proof} The proof uses some knowledge of the representation theory of $U_q(\sl_n)$, which we will simply quote. These facts are quite standard, though difficult to find in the
literature; the paper \cite{AST} certainly helps. The missing proofs are straightforward, analogous to the detailed proofs in the appendix to \cite{ELib}.

Let the \emph{Weyl module} $W(\l)$ be the submodule of $V_{\un{w}} \ot \Bbbk$ for $\un{w} \in P(\l)$, generated by the highest weight vector $x_{\top}$. It is free over $\Bbbk$ and
finite rank. There is a duality functor on modules over $U_q(\sl_n) \ot \Bbbk$ which are free of finite rank over $\Bbbk$, sending $W$ to $\Hom_{\Bbbk}(W, \Bbbk)$, and letting
$U_q(\sl_n)$ act via the antipode. Thus the \emph{dual Weyl module} $W(\l)^*$ is defined.

One can show that $\Hom_{U_q}(W(\l), W^*(\mu))$ is free of rank 1 over $\Bbbk$ when $\l = \mu$, and is zero otherwise. One can also show that $\Ext^1(W(\l), W^*(\mu)) = 0$ for all $\l,
\mu$. Both of these can be shown using the universal properties of Weyl and dual Weyl modules.

A module is \emph{tilting} if it has a filtration whose subquotients are Weyl modules, and another filtration whose subquotients are dual Weyl modules. One can show that tensoring with
a fundamental representation will send a tilting module to another tilting module. Thus $V_{\un{w}} \ot \Bbbk$ is tilting. One can also compute the multiplicities of Weyl and dual Weyl
modules in this filtration, which agrees with the plethyism rules for tensoring representations in the semisimple version of $\sl_n$-representations.

By the Hom and Ext vanishing results above, one can compute the size of a morphism space between two tilting modules by pairing the multiplicities in a Weyl filtration with the
multiplicities in a dual Weyl filtration. The purpose of this entire discussion is to deduce that $\Hom(V_{\un{w}}, V_{\un{y}})$ is a free $\Bbbk$-module whose rank is equal to the number of pairs $\eb \subset \un{w}$ and $\fb \subset \un{y}$ such that $\un{w}^\eb = \un{y}^\fb$.

By the computations of this chapter, we know that the functor $\Gamma_n$ is faithful after specialization from $\Bbbk$ to any field. We also know that it is full after specialization,
because both sides are vector spaces of the same dimension. It is a standard application of Nakayama's lemma in algebraic geometry, that a map of finitely-generated $\Bbbk$-modules which
becomes a surjective after every specialization to a field was already surjective. Thus $\Gamma_n$ was full (and faithful) before specialization. In particular, $\Gamma_n$ is an
equivalence. \end{proof}

\section{Clasps}
\label{sec-clasps}

Suppose that $\un{w}, \un{x} \in P(\l)$. There is a unique double ladder $\LL^\l_{\eb, \fb}$ from $\un{w}$ to $\un{x}$ in cell $\l$, for which both $\eb$ and $\fb$ are the full weight
subsequence. This map $\LL^\l_{\eb, \fb}$ is a neutral ladder. We call it the \emph{neutral double ladder}. For an arbitrary morphism $f \in \Hom(\un{w}, \un{x})$ expressed in the
double ladders basis, the coefficient of the neutral double ladder is called the \emph{neutral coefficient} of $f$.

Recall the Neutral Ladder Lemma, which stated that any two neutral ladders with the same source and target are equal modulo lower terms. For this reason, the neutral coefficient of $f$
does not depend on the choices of neutral ladders made in the construction of the double ladders basis, so it is actually intrinsic to the morphism $f$. Moreover, composing $f$ with a neutral ladder will not change the neutral coefficient.

\subsection{Definition and basic properties}
\label{subsec-definition}

As discussed in the introduction, a clasp is a morphism which represents projection from a sequence $\un{w}$ of fundamental weights (i.e. an object in $\Fund^+$) to the irreducible
summand corresponding to the total weight of the sequence. It is orthogonal to the other summands, which all have lower weights. This motivates the formal definition.

\begin{defn} Let $\l$ be a dominant weight, and let $\un{w}, \un{x} \in P(\l)$. Then a morphism $\phi \co \un{w} \to \un{x}$ is a \emph{clasp} if it is killed by postcomposition with any
strictly outward ladder, and if the neutral coefficient is $1$. We call it a \emph{$\l$-clasp} when $\l$ is significant but the specific choice of $\un{w}, \un{x} \in P(\l)$ is not.
\end{defn}

\begin{prop} Let $\un{w}, \un{x} \in P(\l)$. A clasp from $\un{w}$ to $\un{x}$, if it exists, is unique. It is also the unique morphism with neutral coefficient $1$ that is killed by
precomposition with any strictly inward ladder. The postcomposition of a clasp $\un{w} \to \un{x}$ with a neutral ladder $\un{x} \to \un{y}$ is a clasp $\un{w} \to \un{y}$, and similarly
for precomposition. Thus if any $\l$-clasp exists, then all $\l$-clasps exist. The composition of two clasps is a clasp. Clasps are sent to clasps by the duality involution. \end{prop}

\begin{proof} Let $T = T_{\un{w}, \un{x}}$ be the subspace of $\Hom(\un{w}, \un{x})$ which is killed by placing a strictly outward ladder on top. Let $B = B_{\un{w}, \un{x}}$ be the
subspace of $\Hom(\un{w}, \un{x})$ killed by placing a strictly inward ladder on bottom. The existence of a clasp will imply the existence of an element of $T$ with invertible neutral
coefficient. Our first task is to demonstrate that $T_{\un{y}, \un{z}}$ and $B_{\un{y}, \un{z}}$ also contain an element with invertible neutral coefficient, for all $\un{y}, \un{z} \in P(\l)$.

It is clear that precomposing an element $f \in T_{\un{w}, \un{x}}$ with any morphism $\un{y} \to \un{w}$ will produce an element (possibly zero) of $T_{\un{y}, \un{x}}$. Precomposing
with a neutral ladder will not change the neutral coefficient, yielding a nonzero element of $T_{\un{y}, \un{x}}$. Meanwhile, postcomposing with a neutral ladder $N \co \un{x} \to \un{z}$
will produce an element of $T_{\un{w}, \un{z}}$, also with the same neutral coefficient. After all, if $D$ is a strictly outward ladder beginning at $\un{z}$, then $D \circ (N \circ f) =
(D \circ N) \circ f = 0$ because $(D \circ N)$ is a strictly outward ladder. Thus if $T_{\un{w}, \un{x}}$ has an element with an invertible neutral coefficient, so does every $T_{\un{y},
\un{z}}$, obtained by applying neutral ladders to either side.

The duality involution will interchange $T_{\un{w}, \un{x}}$ and $B_{\un{x}, \un{w}}$, without changing the neutral coefficient. Thus the arguments above apply equally to $B$.

In fact, for a fixed element $f \in T_{\un{w}, \un{x}}$, postcomposing with arbitrary neutral ladders can yield only one morphism in any space $T_{\un{w}, \un{z}}$. This is because any
two neutral ladders agree modulo lower terms, and these lower terms are orthogonal to $T$. A similar statement holds when precomposing an element of $B$ with a neutral ladder.

Now, assuming that a clasp exists, we will show that $T$ is one-dimensional, that any element is determined by its neutral coefficient, and that $T = B$. The trick is to multiply $f \in T$ with $g \in B$, and show that the result is equal to both $f$ and $g$ (after translating by neutral ladders).

Let $f \in T_{\un{w}, \un{x}}$, with $f = c N + L$, where $N$ is the neutral double ladder, $L$ is the lower terms, and $c$ is the neutral coefficient. Let $g \in B_{\un{x}, \un{y}}$, with $g = d M + K$ where $M$ is neutral and $K$ is lower terms. Then
\[g \circ f = d (M \circ f) \in T_{\un{w}, \un{y}},\]
because $K \circ f = 0$. However, we also have
\[g \circ f = c (g \circ N) \in B_{\un{w}, \un{y}}, \]
since $g \circ L = 0$.  In particular, if the neutral coefficient $d$ is invertible, then $f$ can be recovered from $g$ by composing further with some neutral ladder $P \co \un{y} \to \un{x}$, since $P \circ M \circ f = f$. Thus 
\[f = d^{-1} P \circ g \circ f = cd^{-1} P \circ g \circ N.\]
In particular, $f$ is obtained from $g$ by applying neutral ladders, and is thus also in $B$. Since $f$ was arbitrary, this implies that $T \subset B$. A similar argument shows that $B \subset T$, and thus $B = T$.

Moreover, if $f' = c N + L'$ is another element of $T_{\un{w}, \un{x}}$ with the same neutral coefficient, then we also have $f' = c d^{-1} O \circ g \circ N$, from which we deduce that
$f = f'$. \end{proof}

When a clasp exists in $\Hom(\un{w}, \un{x})$, we shall draw it as a rounded rectangle, labeled with the element $\l$, having inputs $\un{w}$ and outputs $\un{x}$. We encode the basic
properties of clasps in the following relations, which also hold after applying the duality involution (i.e. flipping upside-down).

\begin{equation} \label{claspout} {
\labellist
\small\hair 2pt
 \pinlabel {$\l$} [ ] at 41 13
 \pinlabel {$= \quad 0$} [ ] at 100 39
\endlabellist
\centering
\ig{1}{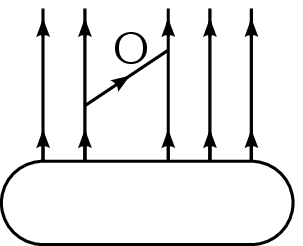}
} \end{equation}
Here $O$ is any outward rung, and the tilt is irrelevant.

\begin{equation} \label{claspN} {
\labellist
\small\hair 2pt
 \pinlabel {$\l$} [ ] at 38 13
 \pinlabel {$\l$} [ ] at 166 27
 \pinlabel {$\un{w}$} [ ] at 80 29
 \pinlabel {$\un{x}$} [ ] at 79 67
 \pinlabel {$\un{x}$} [ ] at 206 47
\endlabellist
\centering
\ig{1}{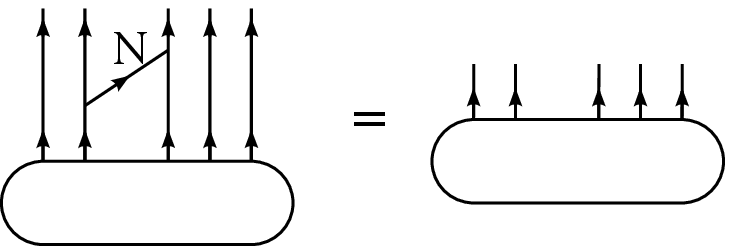}
} \end{equation}
Here $N$ is any neutral rung, and the tilt is irrelevant.

Because of \eqref{claspout} and \eqref{claspN}, we can freely seesaw any rung which meets a clasp, as in \eqref{seesaw1}.

\begin{equation} \label{claspclasp} {
\labellist
\small\hair 2pt
 \pinlabel {$\l$} [ ] at 41 29
 \pinlabel {$\nu$} [ ] at 41 69
 \pinlabel {$\l$} [ ] at 179 44
 \pinlabel {$\un{w}$} [ ] at 80 8
 \pinlabel {$\un{x}$} [ ] at 80 48
 \pinlabel {$\un{y}$} [ ] at 80 92
 \pinlabel {$\un{w}$} [ ] at 219 21
 \pinlabel {$\un{y}$} [ ] at 219 69 
\endlabellist
\centering
\ig{1}{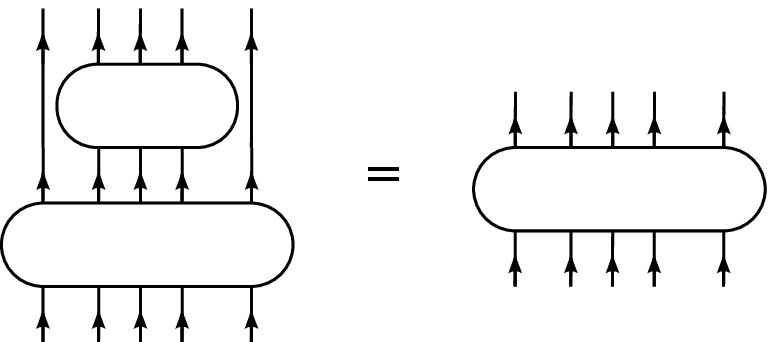}
} \end{equation}
Here, $\nu$ comes from some consecutive subsequence of $\un{x}$ and of $\un{y}$. The case $\nu = \l$ is allowed.

The $\l$-clasp satisfies the three properties (compatibility, orthogonality, unitality) mentioned in the introduction. One major implication of these properties is that the $\l$-clasp picks out a unique new object in the Karoubi envelope of $\Fund$.

\begin{defn} Let $\{X_i\}$ be a family of objects in an additive category $\CC$, and let $\phi_{i,j} \co X_i \to X_j$ be a family of morphisms. Then $\phi = (\phi_{i,j})$ is a
\emph{compatible system of projections} if $\phi_{j,k} \circ \phi_{i,j} = \phi_{i,k}$ for all $i,j,k$. \end{defn}

In the Karoubi envelope of an additive category $\CC$, there is an object for each pair $(X,e)$ of an object $X \in \CC$ and an idempotent $e \in \End(X)$. Given a compatible system, the
idempotents $\phi_{i,i}$ each contribute objects $(X_i,\phi_{i,i})$. Moreover, the morphisms $\phi_{i,j}$ descend to isomorphisms $(X_i, \phi_{i,i}) \to (X_j, \phi_{j,j})$. In fact, a
compatible system of projections is precisely the data required to pin down a common summand of a family of objects, and it gives rise to a single new object in the Karoubi envelope which
is canonically embedded in each object of the original family.

\begin{defn} When the $\l$-clasp exists, we let $V_\l$ or simply $\l$ denote the corresponding object in the Karoubi envelope of $\Fund$. So, for instance, $\l \un{w}$ is the tensor
product of $\l$ with $\un{w}$. \end{defn}

Recall that, for two objects $(X,e)$ and $(Y,f)$ in a Karoubi envelope, one has \[ \Hom((X,e), (Y,f)) = f \Hom(X,Y) e. \]

\begin{prop} \label{prop:blah} Suppose that the $\l$-clasp exists. Then for any sequences $\un{w}$, $\un{x}$ in $\Fund^+$, the morphism space $\Hom(\l \un{w}, \un{x})$ in the Karoubi
envelope can be described as follows. The space $\Hom(\un{y} \un{w}, \un{x})$, for some $\un{y} \in P(\l)$, has a basis given by double ladders $\LL^\mu_{\eb \fb, \gb}$, where $\eb
\subset \un{y}$, $\fb \subset \un{w}$, the concatenation $\eb \fb$ lives in $E(\un{y} \un{w},\mu)$, and $\gb \in E(\un{x}, \mu)$. Consider the subspace spanned by double ladders where
$\eb$ is the full weight subsequence of $\un{y}$. This subspace projects isomorphically to $\Hom(\l \un{w}, \un{x})$, under precomposition with the clasp on $\un{y}$ (tensored with the
identity of $\un{w}$). \end{prop}

\begin{proof} We know that if $\eb$ is not full, then the light ladder $LL_{\un{y},\eb}$ is orthogonal to the $\l$-clasp on $\un{y}$. By the algorithmic construction of light ladders, then $\LL^\mu_{\eb \fb, \gb}$ is also orthogonal to the $\l$-clasp (tensored with the identity of $\un{w}$).

That the double ladders for which $\eb$ is full remain linearly independent after composition with the clasp can be seen after applying the functor $\Gamma_n$. Because the clasp on $\un{y}$ has neutral coefficient $1$, it sends the highest weight vector $x_{\top}$ of $V_{\un{y}}$ to itself (modulo lower terms). Thus, precomposition with the clasp will not affect the upper-triangularity argument used in Theorem \ref{linindepthm}, for weight subsequences $\eb \fb$ where $\eb$ is full. \end{proof}

\begin{cor} \label{cor:orthogonal} Suppose that the $\l$-clasp and the $\nu$-clasp both exist. Then $\Hom(\l, \nu)$ is spanned by the identity map if $\l = \nu$, and is zero otherwise.
\end{cor}
	
\begin{proof} Using the previous proposition, morphisms between them come from double ladders $\LL_{\eb, \fb}$ for which both $\eb$ and $\fb$ are full. This is not possible unless
$\l = \nu$, in which case the neutral ladder $\LL_{\eb, \fb}$ induces the identity map. \end{proof}

\subsection{Intersection forms and triple clasp formulas}
\label{subsec-tripleclasp}

Let $\Om(a)$ denote the set of weights in $V_{\w_a}$, and $\Om^-(a)$ denote $\Om(a) \setminus \{\w_a\}$.

Suppose that the $\l$-clasp exists. From our knowledge of plethyism for $\sl_n$-representations in the semisimple case, we expect $V_\l \ot V_{\w_a}$ to decompose into irreducible
representations of the form $V_{\l + \mu}$ for various $\mu \in \Om(a)$, each appearing with multiplicity one when $\l + \mu$ is dominant, and zero when $\l + \mu$ is not dominant.
Note that plethyism can be slightly more complicated outside of type $A$, as discussed in Proposition \ref{prop:motivational}.

To measure what happens in various specializations of the integral form, we should compute the so-called local intersection forms. Because each summand appears with multiplicity one,
these local intersection forms are $1 \times 1$ matrices, or just numbers.

\begin{defn} Suppose that $\l$ and $\l + \mu$ are dominant, where $\mu \in \Om(a)$. Let $L_{\l, \mu}$ denote any morphism obtained by taking
the elementary light ladder corresponding to $\mu$, and pre- or postcomposing with neutral ladders, to obtain a map $\un{x} a \to \un{y}$ for some $\un{x} \in P(\l)$ and $\un{y} \in
P(\l + \mu)$. When the $\l$-clasp exists, we let $E_{\l, \mu}$ denote the induced map in $\Hom(\l a, \un{y})$, which is independent of the choice of neutral ladders modulo $I_{<
\l+\mu}$. When the $(\l+\mu)$-clasp also exists, we let $E_{\l, \mu}$ denote the induced map in $\Hom(\l a, \l+\mu)$, which is independent of the choice of neutral ladders. \end{defn}

Proposition \ref{prop:blah} above implies that $E_{\l, \mu}$ descends to a basis of the rank-$1$ $\Bbbk$-module $\Hom(\l a, \l + \mu)$. By construction, the map $E_{\l, \mu}$ is the
following composition: include from $\l$ to a sequence $\un{x} \in P(\l)$ that ends with the inputs $x_1 < x_2 < \ldots < x_k$ to the elementary light ladder $E_\mu$, apply the light
ladder $E_\mu$ tensored with the identity on any extraneous strands, and then project back to $\l+\mu$ (if we are working modulo $I_{< \l+\mu}$, the last step is unnecessary). This can be
seen in the picture \eqref{intersectionform} below. The ladder $E_{\mu}$ was independent of $\l$, but $\l$ is required to be big enough so that some $\un{x}$ ending in the desired
sequence exists; this is equivalent to the fact that $\l+\mu$ is dominant.

\begin{defn} Suppose that the $\l$-clasp and the $(\l + \mu)$-clasp exist, for some $\mu \in \Om(a)$. Let $\k_{\l, \mu}$ be the \emph{local intersection form} of $\l a$ at $\l + \mu$,
obtained as follows. Composing $\overline{E_{\l, \mu}}$ with $E_{\l, \mu}$, one obtains an endomorphism of $\l + \mu$, and this endomorphism space is spanned by the identity map. Let
$\k_{\l, \mu}$ be the coefficient of the identity in this composition. \end{defn}

Pictorially, we have the following description of the local intersection form.

\begin{equation} \label{intersectionform} {
\labellist
\tiny\hair 2pt
 \pinlabel {\small $\l + \mu$} [ ] at 71 22
 \pinlabel {\small $\l + \mu$} [ ] at 71 212
 \pinlabel {\small $\l$} [ ] at 57 118
 \pinlabel {\small $\l + \mu$} [ ] at 310 118
 \pinlabel {\small $\k_{\l,\mu}$} [ ] at 208 118
 \pinlabel {$y_1$} [ ] at 44 42
 \pinlabel {$y_2$} [ ] at 67 42
 \pinlabel {$\dots$} [ ] at 93 42
 \pinlabel {$y_k$} [ ] at 116 42
 \pinlabel {$y_{k+1}$} [ ] at 143 42
 \pinlabel {$x_1$} [ ] at 44 96
 \pinlabel {$x_2$} [ ] at 67 96
 \pinlabel {$\dots$} [ ] at 93 96
 \pinlabel {$x_k$} [ ] at 116 96
 \pinlabel {$a$} [ ] at 140 118
 \pinlabel {$x_1$} [ ] at 44 138
 \pinlabel {$x_2$} [ ] at 67 138
 \pinlabel {$\dots$} [ ] at 93 138
 \pinlabel {$x_k$} [ ] at 116 138
 \pinlabel {$y_1$} [ ] at 44 194
 \pinlabel {$y_2$} [ ] at 67 194
 \pinlabel {$\dots$} [ ] at 93 194
 \pinlabel {$y_k$} [ ] at 116 194
 \pinlabel {$y_{k+1}$} [ ] at 143 194
 \pinlabel {\small $\overline{E}_{\mu}$} [ ] at 150 69
 \pinlabel {\small $E_{\mu}$} [ ] at 150 166
\endlabellist
\centering
\ig{1}{intersectionform}
} \end{equation}

\begin{remark} Note that this intersection form is defined even when the $\l + \mu$ clasp is not defined. Instead of pre- and postcomposing with the $(\l+\mu)$-clasp, one considers morphisms modulo $I_{< \l + \mu}$. The identity map will span the appropriate Hom space. Regardless, we will use the picture above as a convenient mnemonic. \end{remark}

The intersection form measures how close the map $E_{\l, \mu}$ (resp. $\overline{E}_{\l, \mu}$) is to being a projection map (resp. inclusion map). If $\k_{\l, \mu}$ is invertible, then
$E_{\l, \mu}$ and a rescaling of $\overline{E}_{\l, \mu}$ compose to the identity of $(\l+\mu)$, and thus pick out $(\l + \mu)$ as a summand inside $\l a$. If the composition is zero, as
it may be in some non-semisimple specialization of $\ZM[q,q^{-1}]$, then $(\l + \mu)$ is not a summand of $\l a$.

Only one of these local intersection forms is immediately easy to compute: when $\mu = \w_a$ is the highest weight. In this case, $E_{\mu}$ is just the identity map. Thus $\k_{\l, \w_a} =
1$.

\begin{prop} (The triple clasp expansion) Suppose that the $\l$-clasp exists, that the $(\l + \mu)$-clasp exists for each $\mu \in \Om^-(a)$, and that all local intersection forms $\k_{\l, \mu}$ are invertible. Then the $(\l + \w_a)$-clasp exists, and it obeys the following formula. The sum is over $\mu \in \Om^-(a)$.

\begin{equation} \label{eq:tripleclasp} {
\labellist
\small\hair 2pt
 \pinlabel {$\l+\w_a$} [ ] at 66 118
 \pinlabel {$\l$} [ ] at 260 118
 \pinlabel {$\l$} [ ] at 466 25
 \pinlabel {$\l+\mu$} [ ] at 482 118
 \pinlabel {$\l$} [ ] at 466 210
 \pinlabel {$\displaystyle - \sum_{\mu} \frac{1}{\k_{\l, \mu}}$} [ ] at 370 118
 \pinlabel {\tiny $a$} [ ] at 138 101
 \pinlabel {\tiny $a$} [ ] at 138 138
 \pinlabel {\tiny $a$} [ ] at 337 113
 \pinlabel {\tiny $a$} [ ] at 547 9
 \pinlabel {\tiny $a$} [ ] at 547 214
 \pinlabel {$E_{\mu}$} [ ] at 560 75
 \pinlabel {$\overline{E}_{\mu}$} [ ] at 560 167
\endlabellist
\centering
\ig{.73}{tripleclaspformula}
} \end{equation} \end{prop}

\begin{proof} The neutral coefficient of the RHS is clearly $1$, as only the first term contributes. Now we need only show that the RHS is orthogonal to $I_{< \l + \w_a}$. The RHS
manifestly factors through $\l a$ on bottom, so by Proposition \ref{prop:blah}, it remains to show that the RHS is orthogonal to light ladders $L_{\l, \nu}$ (where $\nu$ is not the
heighest weight in $V_{\w_a}$). In fact, it is enough to show that the RHS is orthogonal to the maps $E_{\l, \nu} \co \l a \to \l + \nu$ for each $\nu$ (the difference being that $E_{\l,
\nu}$ is $L_{\l, \nu}$ postcomposed with the $(\l+\nu)$-clasp, but this is equal to $L_{\l, \nu}$ modulo lower terms, which in turn factor through other $L_{\l, \nu'}$, and so forth).

If $\mu \ne \nu$ then $\Hom(\l + \mu, \l + \nu) = 0$ by Corollary \ref{cor:orthogonal}. Thus when composing the RHS of \eqref{eq:tripleclasp} with $E_{\l, \nu}$ we can ignore all the
terms on the RHS attached to $\mu \ne \nu$. Thus \[E_{\l, \nu} \circ (\textrm{RHS}) = E_{\l, \nu} - \frac{1}{\k_{\l, \nu}} E_{\l, \nu} \overline{E}_{\l, \nu} E_{\l, \nu} = E_{\l, \nu} -
\frac{\k_{\l, \nu}}{\k_{\l, \nu}} E_{\l, \nu} = 0,\] where we have used \eqref{intersectionform} to replace $E_{\l, \nu} \overline{E}_{\l, \nu}$ with a factor of $\k_{\l, \nu}$.
\end{proof}

Moving the sum over $\mu \in \Om^-(a)$ from the RHS of \eqref{eq:tripleclasp} to the left, one obtains an expression for the identity map of $\l a$ as a sum of orthogonal idempotents
factoring through $\l + \mu$. This formula has the advantage that it is a uniform sum over all $\mu \in \Om(a)$.

\begin{equation} \label{tripleclaspAlt} {
\labellist
\small\hair 2pt
 \pinlabel {$\l$} [ ] at 59 119
 \pinlabel {$\l$} [ ] at 310 24
 \pinlabel {$\l+\mu$} [ ] at 323 119
 \pinlabel {$\l$} [ ] at 310 212
 \pinlabel {$\displaystyle \sum_{\mu \in \Om(a)} \frac{1}{\k_{\l, \mu}}$} [ ] at 201 119
 \pinlabel {$E_{\mu}$} [ ] at 402 74
 \pinlabel {$\overline{E}_{\mu}$} [ ] at 402 167
\tiny\hair 2pt
 \pinlabel {$a$} [ ] at 139 110
 \pinlabel {$a$} [ ] at 390 23
 \pinlabel {$a$} [ ] at 390 212
\endlabellist
\centering
\ig{.73}{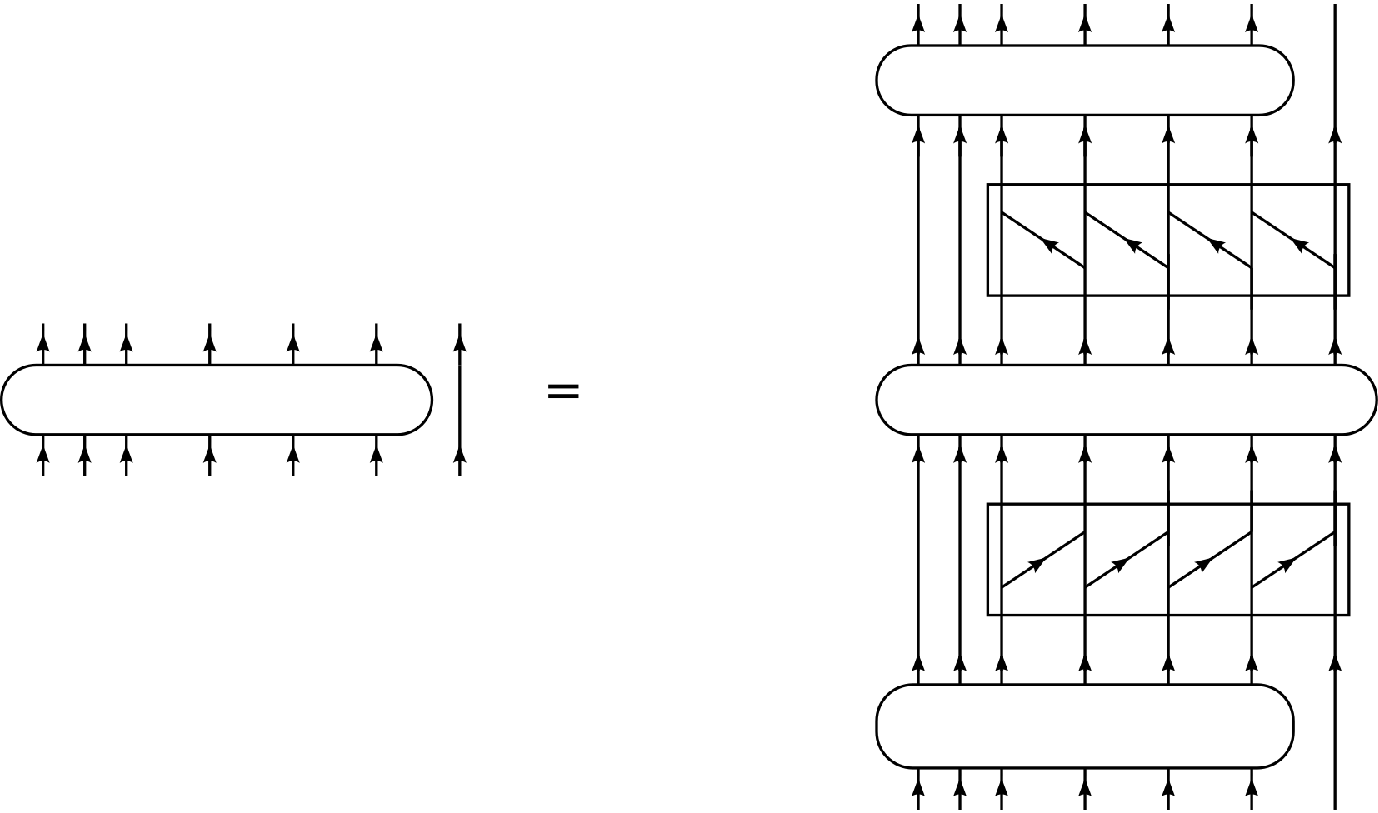}
} \end{equation}

\begin{remark} \label{notsametopbottom} Note that \eqref{eq:tripleclasp} only describes the $(\l + \w_a)$-clasp as a morphism from $\l a$ to $\l a$, or living inside the morphism space
$\Hom(\un{w} a, \un{x} a)$ for any $\un{w}, \un{x} \in P(\l)$. If $\un{y}, \un{z} \in P(\l + \w_a)$ do not end in the same element (e.g. $a$), we can not apply \eqref{eq:tripleclasp} directly to obtain the $(\l + \w_a)$-clasp from $\un{y}$ to $\un{z}$. First, we must apply neutral ladders to $\un{y}$ and $\un{z}$ to obtain sequences which do end in the same element, and then we apply \eqref{eq:tripleclasp}. This works by \eqref{claspN}.  \end{remark}

\subsection{A conjectural formula for the local intersection form}
\label{subsec-conjecture}

The triple clasp expansion is essentially tautological, given the form of the double ladders basis. The interesting part is the computation of the local intersection forms
$\k_{\l, \mu}$. We have computed these for $n \le 4$, and believe that the general answer obeys a particularly nice formula.

Let us write $\l = (\l_1, \l_2, \ldots, \l_{n-1})$ for $\l = \l_1 \w_1 + \l_2 \w_2 + \ldots + \l_{n-1} \w_{n-1}$. So when $n=4$ we may write $\l = (b,c,d)$ and we may label clasps with
this triple of numbers. We continue to describe $\mu \in V_{\w_a}$ with its $\gl_n$-weight, a sequence of zeroes and ones of length $n$. We now state the results for $n \le 4$, which all
involve ratios of quantum numbers.

\begin{align}
& \label{n=2LIF} (n=2, a=1): \quad 
\k_{(b), (10)} = 1, \quad 
\k_{(b), (01)} = \frac{[b+1]}{[b]}.
\\
&\nonumber \\
&\label{n=3LIF} (n=3, a=1): \quad 
\k_{(b,c), (100)} = 1, \quad 
\k_{(b,c), (010)} = \frac{[b+1]}{[b]}, \quad 
\k_{(b,c), (001)} = \frac{[c+1]}{[c]}\frac{[b+c+2]}{[b+c+1]}.
\\
&\nonumber (n=3, a=2): \quad  
\k_{(b,c), (110)} = 1, \quad  
\k_{(b,c), (101)} = \frac{[c+1]}{[c]}, \quad  
\k_{(b,c), (011)} = \frac{[b+1]}{[b]}\frac{[b+c+2]}{[b+c+1]}.
\\
&\nonumber \\
&\label{n=4LIF} (n=4, a=1): \quad  
\k_{(b,c,d), (1000)} = 1, \quad  
\k_{(b,c,d), (0100)} = \frac{[b+1]}{[b]}, \quad  
\k_{(b,c,d), (0010)} = \frac{[c+1]}{[c]}\frac{[b+c+2]}{[b+c+1]}, \quad \\
&\nonumber 
\quad \k_{(b,c,d), (0001)} = \frac{[d+1]}{[d]} \frac{[c+d+2]}{[c+d+1]} \frac{[b+c+d+3]}{[b+c+d+2]}.
\\
&\nonumber (n=4, a=2): \quad  
\k_{(b,c,d), (1100)} = 1, \quad  
\k_{(b,c,d), (1010)} = \frac{[c+1]}{[c]}, \quad  
\k_{(b,c,d), (1001)} = \frac{[d+1]}{[d]}\frac{[c+d+2]}{[c+d+1]}, \quad \\
&\nonumber 
\quad \k_{(b,c,d), (0110)} = \frac{[b+1]}{[b]}\frac{[c+b+2]}{[c+b+1]}, \quad 
\k_{(b,c,d), (0101)} = \frac{[b+1]}{[b]}\frac{[d+1]}{[d]}\frac{[b+c+d+3]}{[b+c+d+2]}, \\
&\nonumber \quad \quad \k_{(b,c,d), (0011)} = \frac{[c+1]}{[c]} \frac{[c+d+2]}{[c+d+1]}\frac{[c+b+2]}{[c+b+1]}\frac{[b+c+d+3]}{[b+c+d+2]}. \\
&\nonumber (n=4, a=3): \quad  
\k_{(b,c,d), (1110)} = 1, \quad  
\k_{(b,c,d), (1101)} = \frac{[d+1]}{[d]}, \quad 
\k_{(b,c,d), (1011)} = \frac{[c+1]}{[c]}\frac{[c+d+2]}{[c+d+1]}, \quad \\
&\nonumber 
\quad \k_{(b,c,d), (0111)} = \frac{[b+1]}{[b]}\frac{[b+c+2]}{[b+c+1]}\frac{[b+c+d+3]}{[b+c+d+2]}.
 \end{align}
 
It should be noted that the denominator has one factor which may vanish (namely $[b]$, $[c]$, or $[d]$), and other factors which can not vanish (like $[b+c+1]$) for a dominant weight
$\l$. Moreover, the denominator vanishes precisely when $\l + \mu$ is not a dominant weight, so that $\k_{\l, \mu}$ is not (and need not be) defined in those cases.

As is evident, each local intersection form is a product of factors of the form $\frac{[n]}{[n-1]}$. Moreover, the numerator $[n]$ of each factor satisfies $n = \langle \l + \rho, \a
\rangle$ for some positive root $\a$, where $\rho = \sum \w_i$ is the half-sum of the positive roots. Let us recall some notation from the introduction.

\begin{defn} Fix $n \ge 2$. Let $\rho$ denote the half-sum of the positive roots, which is also the sum of all fundamental weights. For any dominant weight $\l$ and
positive root $\a$ let $A(\l, \a) = \langle \l + \rho, \a \rangle$. \end{defn}

\begin{defn} For any weight $\mu$ in a fundamental representation $V_{\w_a}$, let $w_{\mu}$ be the unique element of $S_n$ which
sends $\w_a$ to $\mu$, and has minimal length (i.e. is a minimal coset representative for $S_a \times S_{n-a}$, the stabilizer of $\w_a$). Let $\Phi(\mu)$ denote the subset of positive
roots $\a$ for which $w_{\mu}^{-1}(\a)$ is a negative root. \end{defn}

\begin{claim} \label{minusoneclaim} The roots $\a \in \Phi(\mu)$ are in bijection with pairs, inside $\mu$, of a zero occuring before a one. This bijection sends the pair to the root
corresponding to the transposition in $S_n$ which swaps the zero and the one. In particular, for any dominant weight $\l$, $A(\l+\mu, \a) = A(\l,\a)-1$ if and only if $\a \in \Phi(\mu)$.
\end{claim}

\begin{claim} \label{zeroclaim} For a dominant weight $\l$, $A(\l,\a)-1 = 0$ for some $\a \in \Phi(\mu)$ if and only if $\l + \mu$ is not dominant. \end{claim}

\begin{proof} These are straightforward exercises. \end{proof}
	
We now state our main conjecture, which holds for $n \le 4$ by the above formulas.

\begin{conj} \label{MainConj} Let $\l$ be dominant, and $\mu \in \Om(a)$. We conjecture that, whenever $\l + \mu$ is dominant, one has
\begin{equation} \label{mainconj} \k_{\l, \mu} = \prod_{\a \in \Phi(\mu)} \frac{[A(\l, \a)]}{[A(\l, \a)-1]}, \end{equation}
or equivalently, by Claim \ref{minusoneclaim}, one has
\begin{equation} \label{mainconj2} \k_{\l, \mu} = \prod_{\a \in \Phi(\mu)} \frac{[A(\l, \a)]}{[A(\l+\mu, \a)]}. \end{equation}
These products are well-defined by Claim \ref{zeroclaim}. \end{conj}

In the next section, we discuss a strategy to prove this conjecture. In the meantime, let us discuss the relationship between the conjecture and the Weyl dimension formula.

The quantum Weyl dimension formula states that the (graded) dimension of the irreducible representation $V_\l$ (in the semisimple case) is expressible as a product of ratios of quantum numbers. Let $\Phi^+$ denote the set of all positive roots.
\begin{equation} \label{weyldim} \dim_q V_\l = \prod_{\a \in \Phi^+} \frac{[A(\l,\a)]}{[A(0,\a)]}. \end{equation}
Thus the formula \eqref{mainconj2} for $\pm \k_{\l, \mu}$ is similar to the formula for $\frac{\dim_q V_\l}{\dim_q V_{\l + \mu}}$, except that it only involves positive roots in $\Phi(\mu)$, which is a proper subset of $\Phi^+$ so long as $n > 2$. The missing roots are significant (they don't just cancel out in the ratio) already in $\sl_3$.

The (graded) dimension of $V_\l$, also thought of as the ``graded trace'' of the identity map, is seen pictorially by closing up the identity map of $V_{\l}$ around a circle. (These next two diagrams use downward-oriented strands.)
\begin{equation} \label{claspcircle} {
\labellist
\small\hair 2pt
 \pinlabel {$\l$} [ ] at 59 61
 \pinlabel {\large $\dim_q V_\l$} [ ] at -25 66
\endlabellist
\centering
\ig{1}{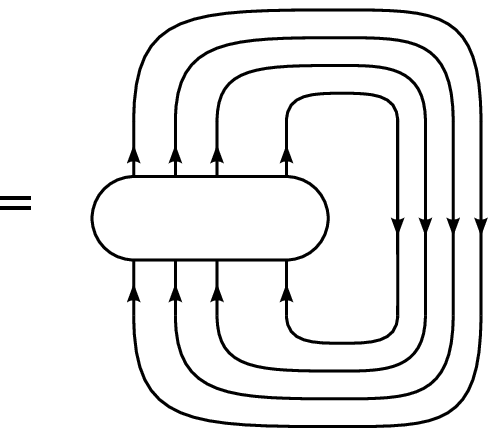}
} \end{equation}
To relate the Weyl dimension formula and the intersection form formula let us close up \eqref{intersectionform} on a circle.
\begin{equation} \quad \label{intersectionformcircle} {
\labellist
\small\hair 2pt
 \pinlabel {$\k_{\l, \mu} \dim_q V_{\l+\mu}$} [ ] at 26 154
 \pinlabel {$\l+\mu$} [ ] at 162 58
 \pinlabel {$\l+\mu$} [ ] at 162 247
 \pinlabel {$\l+\mu$} [ ] at 412 156
 \pinlabel {$\l+\mu$} [ ] at 606 156
 \pinlabel {$\l$} [ ] at 148 153
 \pinlabel {$\l$} [ ] at 424 250
 \pinlabel {$\l$} [ ] at 424 63
 \pinlabel {$\l$} [ ] at 616 63
 \pinlabel {$\l$} [ ] at 616 250
 \tiny
 \pinlabel {$a$} [ ] at 202 153
 \pinlabel {$a$} [ ] at 349 147
 \pinlabel {$a$} [ ] at 560 213
 \pinlabel {$\overline{E}_{\l,\mu}$} [ ] at 213 103
 \pinlabel {$E_{\l, \mu}$} [ ] at 213 202
\endlabellist
\centering
\ig{.6}{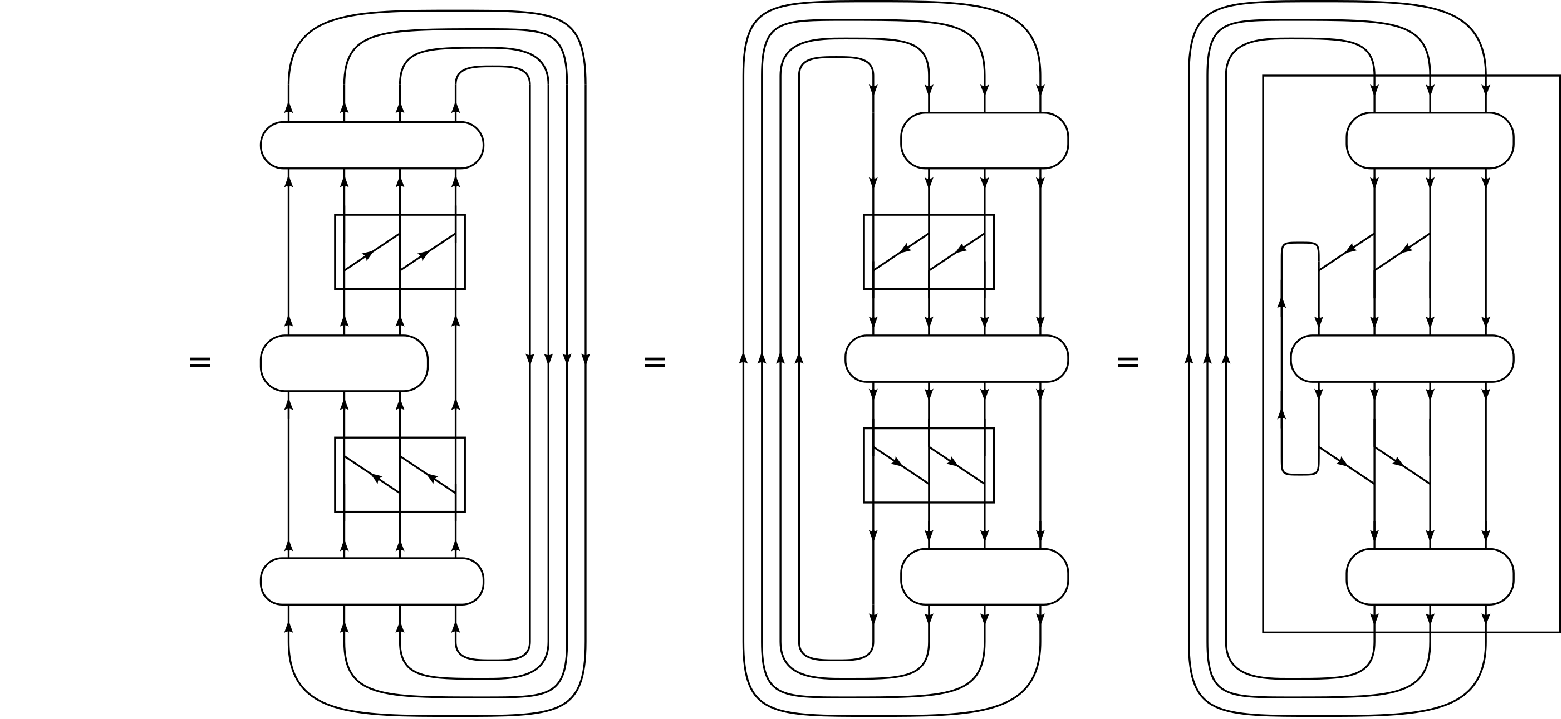}
} \end{equation}
We begin with $\k_{\l, \mu} \dim_q V_{\l+\mu}$, which is the closure of the RHS of \eqref{intersectionform}. We then double the $\l$-clasp (it is an idempotent) and wrap the morphism around the circle, merging the two $(\l+\mu)$-clasps. Finally, we observe that what lies within the big box in the last diagram is an endomorphism of $\l$. Suppose that this box is equal to $\t_{\l, \mu}$ times the identity of $\l$. Then the result is just $\t_{\l, \mu} \dim_q V_{\l}$. We have shown that
\begin{equation} \label{kvst} \frac{\dim_q V_{\l}}{\dim_q V_{\l + \mu}} = \frac{\k_{\l, \mu}}{\t_{\l, \mu}}. \end{equation}
Thus, if Conjecture \ref{MainConj} is correct, then $\t_{\l, \mu}$ has a formula involving a product over positive roots not in $\Phi(\mu)$.

In the next section we outline an inductive method to compute $\k_{\l, \mu}$. To compute $\t_{\l, \mu}$ seems to require a very similar method; it does not appear to be any easier to
compute $\t$ than $\k$. We have explained this connection to quantum dimensions partially in order to discourage the reader from envisaging an easy way out of the computation using
circular closures. However, a representation-theoretic explanation for the formula in Conjecture \ref{MainConj} is sorely lacking, and it is possible that some trace argument could go a
long way towards clarifying the mystery.

One extremely deep implication of the conjecture is that local intersection forms are all positive! This has been discussed in the introduction.

\subsection{Deriving a recursive formula: part I}
\label{subsec-recursive1}

Suppose we have already computed $\k_{\s, \nu}$ for $\s < \l$, and we wish to compute $\k_{\l, \mu}$ by examining the LHS of \eqref{intersectionform}. Using \eqref{eq:tripleclasp} we
can replace the central clasp in \eqref{intersectionform} with a number of more complicated diagrams. The term which appears with coefficient $1$ is
\begin{equation} \label{intersectionterm1} {
\labellist
\tiny\hair 2pt
 \pinlabel {\small $\l + \mu$} [ ] at 71 22
 \pinlabel {\small $\l + \mu$} [ ] at 71 210
 \pinlabel {\small $\l - \w_{x_k}$} [ ] at 50 116
 \pinlabel {$y_1$} [ ] at 44 42
 \pinlabel {$y_2$} [ ] at 67 42
 \pinlabel {$\dots$} [ ] at 93 42
 \pinlabel {$y_k$} [ ] at 116 42
 \pinlabel {$y_{k+1}$} [ ] at 144 42
 \pinlabel {$x_1$} [ ] at 44 96
 \pinlabel {$x_2$} [ ] at 67 96
 \pinlabel {$\dots$} [ ] at 93 96
 \pinlabel {$x_k$} [ ] at 116 96
 \pinlabel {$a$} [ ] at 140 118
 \pinlabel {$x_1$} [ ] at 44 138
 \pinlabel {$x_2$} [ ] at 67 138
 \pinlabel {$\dots$} [ ] at 93 138
 \pinlabel {$x_k$} [ ] at 116 138
 \pinlabel {$y_1$} [ ] at 44 194
 \pinlabel {$y_2$} [ ] at 67 194
 \pinlabel {$\dots$} [ ] at 93 194
 \pinlabel {$y_k$} [ ] at 116 194
 \pinlabel {$y_{k+1}$} [ ] at 144 194
 \pinlabel {\small $\overline{E}_{\mu}$} [ ] at 150 69
 \pinlabel {\small $E_{\mu}$} [ ] at 150 166
\endlabellist
\centering
\ig{.9}{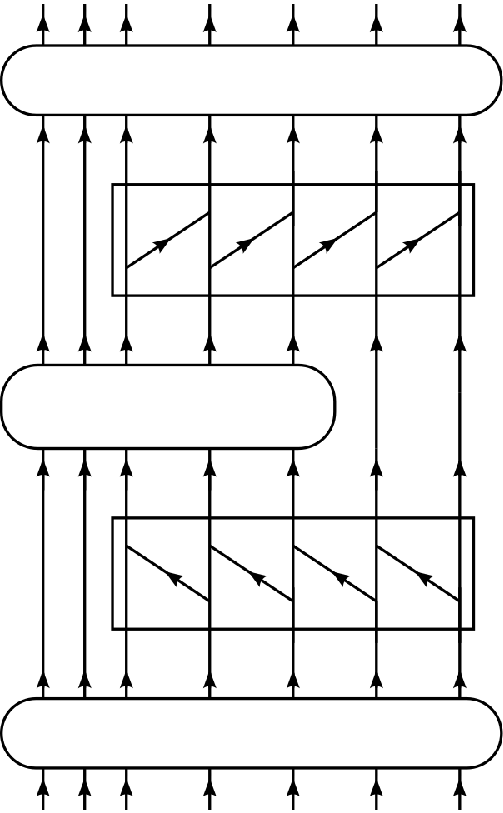}
} \end{equation}
Let us apply a rung swap \eqref{rungswap} to the rightmost rungs. Recall that $\b_i$ denotes the label on the $i$-th crossbar in $E_{\mu}$, so that $\b_k$ is the rightmost crossbar, and that $\a_i$ denotes the label on the middle segment of the $i$-th upright in $E_{\mu}$, sandwiched between $x_i$ and $y_i$. When we apply the rung swap, some number $t$ will be subtracted from the labels on both crossbars, so that $\b_k - t$ remains. The remaining rungs can be slid to the top or bottom of the diagram by associativity. Since $y_k < y_{k+1}$, if $\b_k-t \ne 0$, the rung which reaches the clasp will be orthogonal to it by \eqref{claspout}. Only the $t = \b_k$ term survives. Hence the diagram in \eqref{intersectionterm1} is equal to
\begin{equation} \label{intersectionterm2} {
\labellist
\small\hair 2pt
 \pinlabel {\LARGE $\qbinom{y_{k+1} - \a_k}{\b_k}$} [ ] at 37 117
 \pinlabel {$\l+\mu$} [ ] at 149 22
 \pinlabel {$\l - \w_{x_k}$} [ ] at 120 116
 \pinlabel {$\l+\mu$} [ ] at 148 210
 \pinlabel {\tiny $a$} [ ] at 217 116
 \pinlabel {\tiny $\a_k$} [ ] at 192 116
\endlabellist
\centering
\ig{.8}{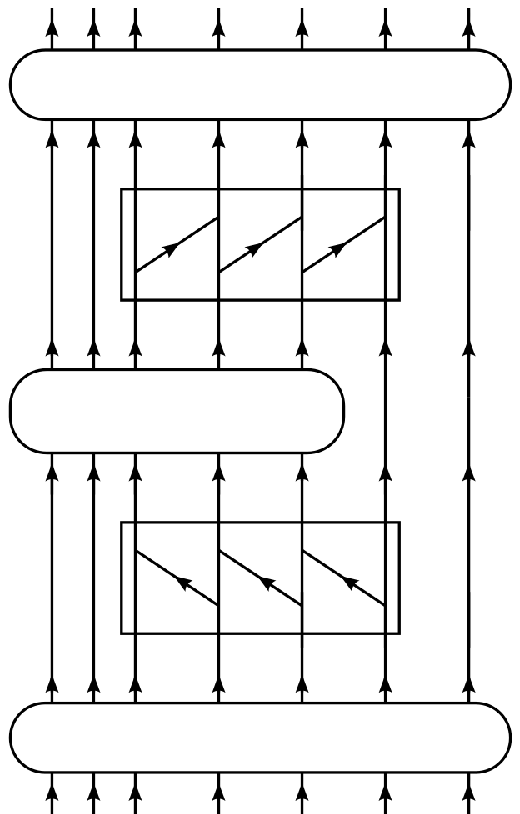}
} \end{equation}
What remains inside the boxes are still elementary light ladders: they are $E_{\mu^-}$ and $\overline{E}_{\mu^-}$, where $\mu^-$ is $\mu$ with the last 1-string removed, a weight inside $V_{\w_{\a_k}}$. Hence, \eqref{intersectionterm2} contains a smaller intersection form, and is equal to \[ \qbinom{y_{k+1}-\a_k}{\b_k} \k_{\l - \w_{x_k}, \mu^-}.\]
Note also that when $\mu$ has a single 1-string, $\a_k = 0$, and $\mu^-$ is the zero weight. While $E_{\mu^-}$ is not technically defined for $\mu^- = 0$, it has the obvious interpretation as an identity diagram. If we extend our notation to let $\k_{\l, 0} = 1$, then the formula above makes sense in all cases.

Meanwhile, when \eqref{eq:tripleclasp} is applied to the middle of \eqref{intersectionform}, a number of other terms appear, parametrized by $\nu \in \Om^-(x_k)$ for which $\l - \w_{x_k} + \nu$ is dominant. We temporarily use $\l^-$ as shorthand for $\l - \w_{x_k}$, so that these diagrams appear with coefficients $- \frac{1}{\k_{\l^-, \nu}}$.
\begin{equation} \label{intersectionterm3} {
\labellist
\small\hair 2pt
 \pinlabel {$\l+\mu$} [ ] at 72 21
 \pinlabel {$\l^-$} [ ] at 47 113
 \pinlabel {$\l^- + \nu$} [ ] at 60 208
 \pinlabel {$\l^-$} [ ] at 47 304
 \pinlabel {$\l+\mu$} [ ] at 73 395
 \pinlabel {$\overline{E}_{\mu}$} [ ] at 150 69
 \pinlabel {$E_{\nu}$} [ ] at 150 163
 \pinlabel {$\overline{E}_{\nu}$} [ ] at 150 259
 \pinlabel {$E_{\mu}$} [ ] at 150 351
 \pinlabel {\tiny $a$} [ ] at 141 302
 \pinlabel {\tiny $x_k$} [ ] at 115 302
 \pinlabel {\tiny $x_k$} [ ] at 115 115
 \pinlabel {\tiny $a$} [ ] at 141 115
\endlabellist
\centering
\ig{.75}{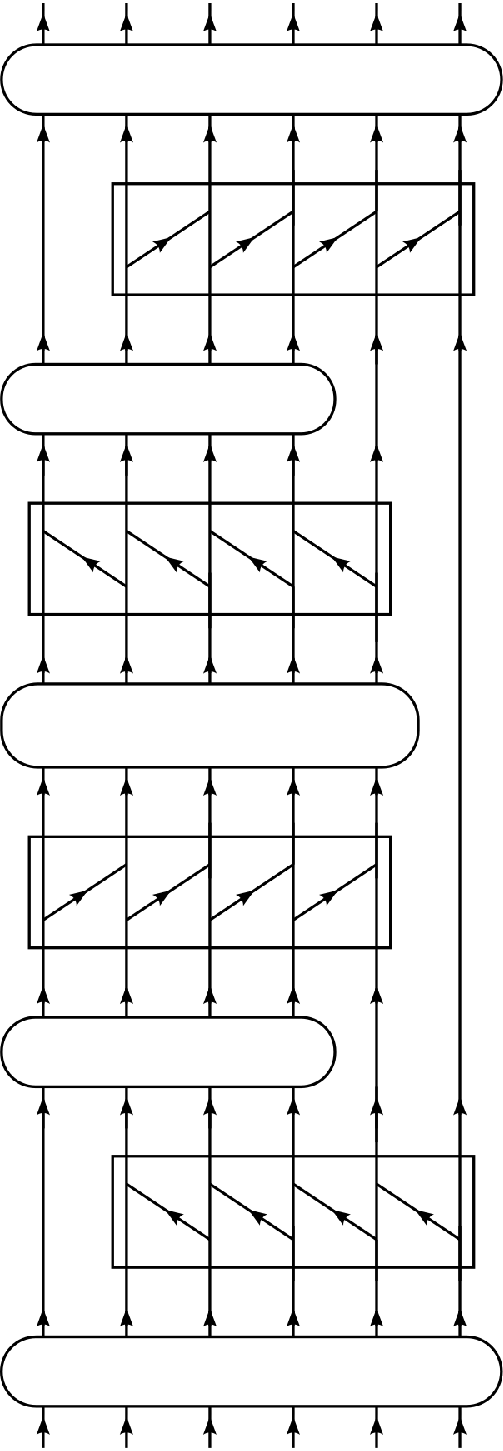}
} \end{equation}
The top half of \eqref{intersectionterm3} is a morphism in $\Hom((\l^- + \nu)a, \l+\mu)$, which is either one or zero dimensional. If it is one dimensional, it is spanned by some $E_{\l^- + \nu, \s}$ for some $\s \in V_{\w_a}$, where necessarily $\l^- + \nu + \s = \l + \mu$ so that $\s = \mu + \w_{x_k} - \nu$. We define the coefficient $\g_{\l, \mu, \nu}$ as below.
\begin{equation} \label{gammaterm} {
\labellist
\tiny\hair 2pt
 \pinlabel {\small $\l - \w_{x_k} + \nu$} [ ] at 59 12
 \pinlabel {\small $\l - \w_{x_k}$} [ ] at 46 107
 \pinlabel {\small $\l+\mu$} [ ] at 70 200
 \pinlabel {\small $\overline{E}_{\nu}$} [ ] at 153 59
 \pinlabel {\small $E_{\mu}$} [ ] at 152 154
 \pinlabel {\small $\g_{\l, \mu, \nu}$} [ ] at 222 101
 \pinlabel {\small $\l - \w_{x_k} + \nu$} [ ] at 306 55
 \pinlabel {\small $\l+\mu$} [ ] at 321 148
 \pinlabel {\small $E_{\s}$} [ ] at 405 101
 \pinlabel {$y_{k+1}$} [ ] at 144 183
 \pinlabel {$y_k$} [ ] at 115 183
 \pinlabel {$\dots$} [ ] at 93 183
 \pinlabel {$y_2$} [ ] at 67 183
 \pinlabel {$y_1$} [ ] at 43 183
 \pinlabel {$x_k$} [ ] at 115 124
 \pinlabel {$\dots$} [ ] at 93 124
 \pinlabel {$x_2$} [ ] at 67 124
 \pinlabel {$x_1$} [ ] at 43 124
 \pinlabel {$a$} [ ] at 138 107
 \pinlabel {$t_{l+1}$} [ ] at 120 32
 \pinlabel {$t_l$} [ ] at 91 32
 \pinlabel {$\dots$} [ ] at 70 32
 \pinlabel {$t_2$} [ ] at 44 32
 \pinlabel {$t_1$} [ ] at 18 32
 \pinlabel {$s_l$} [ ] at 93 92
 \pinlabel {$\dots$} [ ] at 70 92
 \pinlabel {$s_2$} [ ] at 44 92
 \pinlabel {$s_1$} [ ] at 18 92
 \pinlabel {$q_{m+1}$} [ ] at 394 132
 \pinlabel {$q_m$} [ ] at 364 132
 \pinlabel {$\dots$} [ ] at 341 132
 \pinlabel {$q_2$} [ ] at 316 132
 \pinlabel {$q_1$} [ ] at 291 132
 \pinlabel {$p_m$} [ ] at 366 72
 \pinlabel {$\dots$} [ ] at 341 72
 \pinlabel {$p_2$} [ ] at 316 72
 \pinlabel {$p_1$} [ ] at 291 72
 \pinlabel {$a$} [ ] at 388 57
\endlabellist
\centering
\ig{.9}{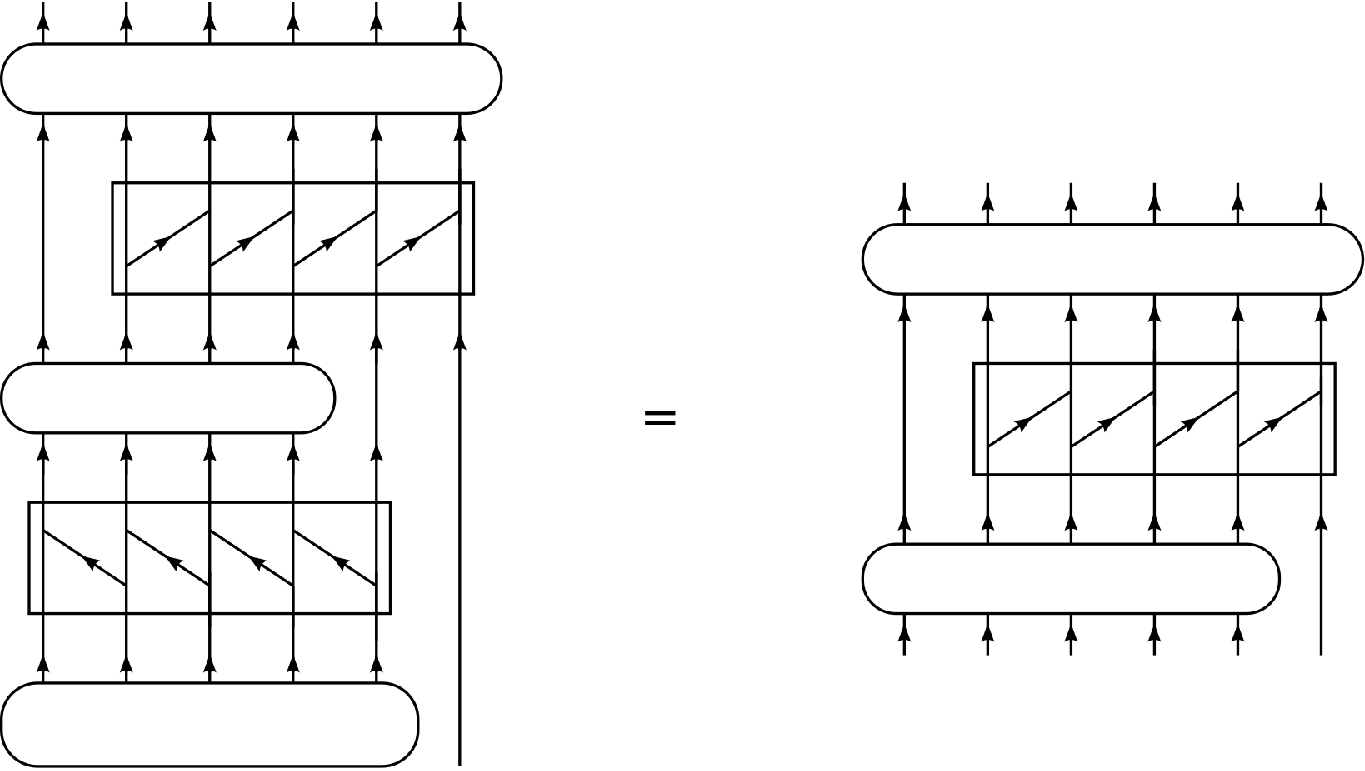}
} \end{equation}

For a given $\mu$ and $\nu$, it is possible that $\s = \mu + \w_{x_k} - \nu$ is not actually a weight in $V_{\w_a}$. In this case, the LHS of \eqref{gammaterm} is just zero, living in a zero-dimensional Hom space. In this case, we say that $\g_{\l, \mu, \nu}=0$, and that is vanishes for trivial reasons.

Resolving the top and bottom half of \eqref{intersectionterm3} using \eqref{gammaterm}, we see that \eqref{intersectionterm3} is equal to \[ \k_{\l^- + \nu, \s} \g_{\l, \mu, \nu}^2. \]
Adding the terms together, we get our first recursive formula.
\begin{equation} \label{recursive1} \k_{\l, \mu} = \qbinom{y_{k+1} - \a_k}{\b_k} \k_{\l^-, \mu^-} - \sum_{\nu \in \Om^-(x_k)} \frac{\k_{\l^- + \nu, \s}}{\k_{\l^-, \nu}} \g_{\l, \mu, \nu}^2. \end{equation}
The sum is over all $\nu \in \Om^-(x_k)$ for which $\l^- + \nu$ is dominant. It remains to produce a recursive formula for $\g_{\l, \mu, \nu}$. Unfortunately, we have not been able to produce a formula for $\g_{\l, \mu, \nu}$ in all cases. Let us explore what can be shown about $\g$, to illustrate the general principles involved in its computation, and some of the subtleties.

There are some hidden factors and assumptions in the picture \eqref{gammaterm}. We have tried to draw our pictures to be fairly general, and handle an arbitrary weight
$\l$, but as a result we have made some sacrifices. Let us emphasize in words what is hard to emphasize in pictures.

The number of rungs involved in $E_{\mu}$, $E_{\nu}$, and $E_{\s}$ can be very different (unlike the picture above, where they all have four rungs). There are many
additional strands involved in these clasps (like the extra strand to the left of $E_{\mu}$), because the strands entering a $\l$-clasp must be a sequence in $P(\l)$. We will tend
not to draw these extra strands when they play no explicit role. However, they do play a significant ``secret" role here. The labels $\{x_i\}$ and $\{s_j\}$ need not match up in any
proscribed fashion (neither will the labels $\{y_i\}$ and $\{q_j\}$, or $\{t_i\}$ and $\{p_j\}$). For the inputs and outputs of the $\l^-$-clasp to match, there must be extra strands to
the left of $E_{\mu}$ labeled by the elements of $\{s_j\}$ which are not elements of $\{x_i\}$, and there must be extra strands to the left of $E_{\nu}$ labeled with the
elements of $\{x_i\}$ which are not elements of $\{s_j\}$. This is the absolute minimum requirement for such a drawing to make sense; in the absence of these extra strands, it would be
impossible to draw this picture, and correspondingly either $\l + \mu$ or $\l^- + \nu$ would not be dominant. There may be even more extra strands for larger $\l$, and these may start to
play a role once the $\l^-$-clasp is resolved.

Let us discuss in more detail the restrictions placed on $\nu$ coming from the fact that $\s = \mu + \w_{x_k} - \nu$ must be an element of $\Om(a)$, in order for $\g$ not to be zero for
trivial reasons. Suppose $\mu = (\dots {\color{red} 001100}{\color{blue}11100})$, so that $y_{k+1}$ is the location of the last 1, $x_k$ is the location of the last 0 before the last 1,
and so forth. The indices $i$ with $i \le x_k$ have been colored red, and the rest blue. We know that $\mu + \w_{x_k} = (\dots {\color{red} 112211}{\color{blue}11100})$. In particular, it
must be the case that $\nu$ and $\s$ must have a red 1 wherever $\mu$ has a red 1, and $\nu$ and $\s$ must have a blue 0 wherever $\mu$ has a blue 0. Thus $t_{l+1}$, the location of the
last 1 in $\nu$, must satisfy $t_{l+1} \le y_{k+1}$. Moreover, because $\nu \in \Om(x_k)$, it must be the case that $t_{l+1} \ge x_k$, with equality if and only if $\nu = \w_{x_k}$. Said
another way, if $\nu \in \Om^-(x_k)$ then $\nu$ has a blue 1 overlapping with some blue 1 in $\mu$, and $\s$ has a blue zero at that location.

\begin{remark} These considerations came from the requirement that $\s$ be a 01-sequence, without any instance of 2 or -1. The reader may worry that we are using $\gl_n$-weights instead
of $\sl_n$-weights; theoretically, the sequence $(122112)$ is a $\gl_n$-weight which agrees as an $\sl_n$-weight with $(011001) \in \Om(3)$. However, because $\s = \mu + \w_{x_k} - \nu$,
the sum of all the indices is still $a$, so there can be no such confusion. \end{remark}

Clearly $\g_{\l, \mu, \nu} = 1$ when $\nu = \w_{x_k}$, so that $\s = \mu$. We continue with another easy case.

\begin{prop} \label{gammaEasyprop} Suppose that $\mu$ is a weight with at most one 0-string before a 1-string. In the notation of \eqref{loweringform}, we have $k=1$. Then $\g_{\l, \mu, \nu} = 1$ when it is not zero for trivial reasons. \end{prop}

\begin{proof} What makes this case easy is that the $\l^-$-clasp pulls directly into the $(\l+\mu)$-clasp, using \eqref{claspclasp}. Applying associativity \eqref{associativity} once, we get the first equality below.
\begin{equation} \label{gammaEasy1} {
\labellist
\small\hair 2pt
 \pinlabel {$\l^- + \nu$} [ ] at 58 12
 \pinlabel {$\l^-$} [ ] at 47 106
 \pinlabel {$\l+\mu$} [ ] at 71 199
 \pinlabel {$\l^- + \nu$} [ ] at 239 53
 \pinlabel {$\l+\mu$} [ ] at 252 163
 \pinlabel {$\l^- + \nu$} [ ] at 420 53
 \pinlabel {$\l+\mu$} [ ] at 432 163
\tiny\hair 2pt
 \pinlabel {$y_1$} [ ] at 116 183
 \pinlabel {$y_2$} [ ] at 139 183
 \pinlabel {$x_1$} [ ] at 116 123
 \pinlabel {$a$} [ ] at 138 122
 \pinlabel {$s_1$} [ ] at 19 91
 \pinlabel {$s_2$} [ ] at 43 91
 \pinlabel {$\dots$} [ ] at 70 91
 \pinlabel {$s_l$} [ ] at 93 91
 \pinlabel {$t_1$} [ ] at 19 33
 \pinlabel {$t_2$} [ ] at 43 33
 \pinlabel {$\dots$} [ ] at 70 33
 \pinlabel {$t_l$} [ ] at 93 33
 \pinlabel {$t_{l+1}$} [ ] at 120 33
 \pinlabel {$s_1$} [ ] at 200 146
 \pinlabel {$s_2$} [ ] at 223 146
 \pinlabel {$\dots$} [ ] at 249 146
 \pinlabel {$s_l$} [ ] at 273 146
 \pinlabel {$y_1$} [ ] at 295 146
 \pinlabel {$y_2$} [ ] at 319 146
 \pinlabel {$t_1$} [ ] at 200 73
 \pinlabel {$t_2$} [ ] at 223 73
 \pinlabel {$\dots$} [ ] at 249 73
 \pinlabel {$t_l$} [ ] at 273 73
 \pinlabel {$t_{l+1}$} [ ] at 299 73
 \pinlabel {$a$} [ ] at 319 66
 \pinlabel {$s_1$} [ ] at 378 146
 \pinlabel {$s_2$} [ ] at 403 146
 \pinlabel {$\dots$} [ ] at 430 146
 \pinlabel {$y_1$} [ ] at 450 146
 \pinlabel {$s_l$} [ ] at 476 146
 \pinlabel {$y_2$} [ ] at 497 146
 \pinlabel {$t_1$} [ ] at 378 73
 \pinlabel {$t_2$} [ ] at 403 73
 \pinlabel {$\dots$} [ ] at 430 73
 \pinlabel {$t_{l}$} [ ] at 450 73
 \pinlabel {$t_{l+1}$} [ ] at 479 73
 \pinlabel {$a$} [ ] at 499 66
\endlabellist
\centering
\ig{.85}{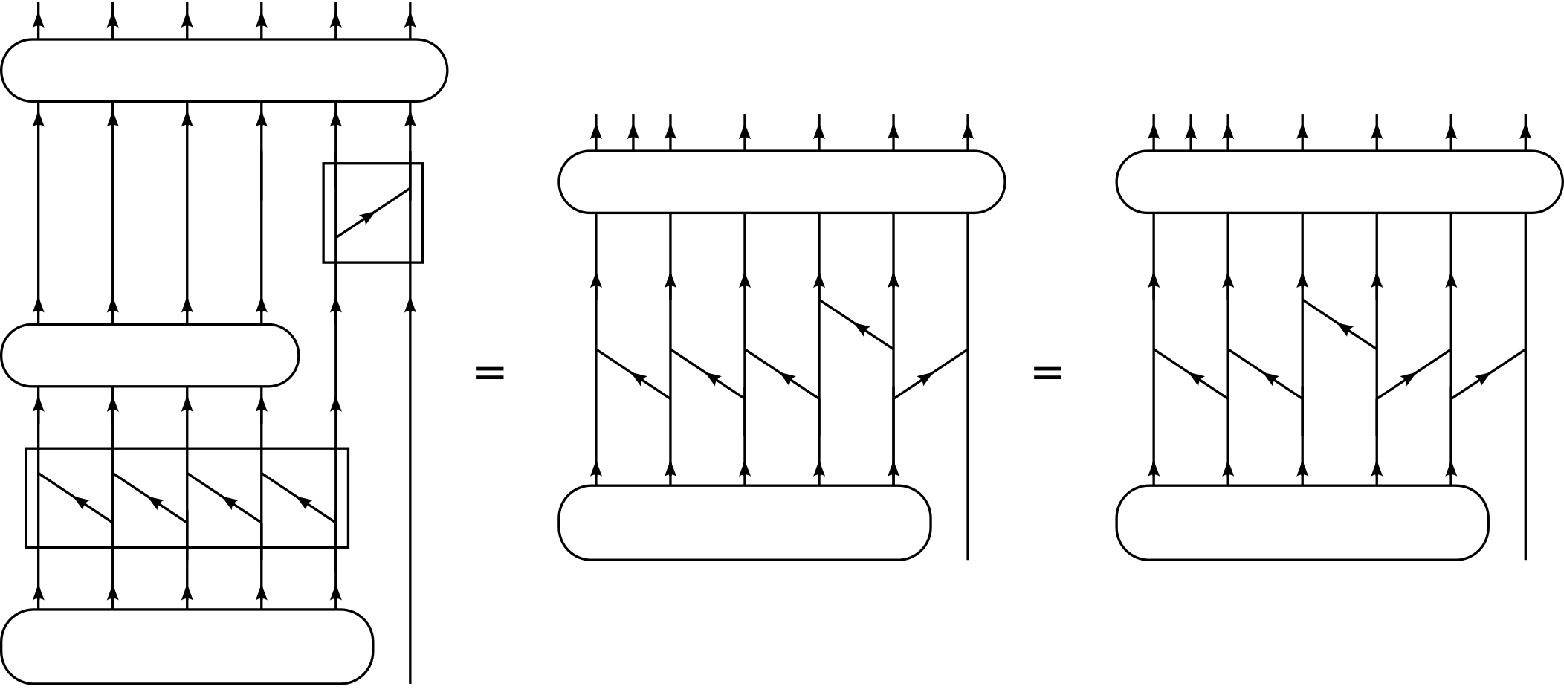}
}\end{equation}
To get the second equality, we seesaw the $(s_l, y_1)$ rung (as in \eqref{seesaw1} but upside-down), and then apply associativity again. Recall that seesawing does not change the morphism modulo lower terms and neutral rungs, so this operation can be applied to any rung adjacent to a clasp.

For each rung that was originally in $\overline{E}_{\nu}$, we pull it to the top using associativity, seesaw it, and continue. The result looks like the RHS of \eqref{gammaterm},
with bottom boundary $(t_1, t_2, \ldots, t_l, t_{l+1}, a)$ and top boundary $(y_1, s_1, s_2, \ldots, s_l, y_2)$. If $y_1 < t_1$ and $t_{l+1} < y_2$ then the result is an elementary light
ladder, and the reader can confirm that it is $E_{\s}$.

Suppose that $y_1 > s_i$ for some $i$. Then at some point, the seesawed rung is actually inward, and is thus orthogonal to the $(\l+\mu)$-clasp. Similarly, if $y_1 > t_1$ then the
leftmost rung of the result is also inward and orthogonal. In either of these cases, the reader can confirm that $\s = \mu + \w_{x_1} - \nu$ is not in $\Om(a)$, so the diagram should vanish for trivial reasons. We have already argued that $\s \in \Om(a)$ implies that $t_{l+1} \le y_2$. Thus we need only see what happens when $y_1 = t_1$ or when $y_2 = t_{l+1}$. In both cases, $\s \in \Om(a)$, so the diagram should not vanish.

If $y_1 = t_1$, the leftmost rung of the result has crossbar labeled $0$, and effectively does not exist. In this case, the remainder is $E_{\s}$ so long as $t_{l+1} < y_2$.

Suppose that $t_{l+1} = y_2$, and the rightmost rung is neutral. We then resolve the diagram along the lines of the proof of Lemma \ref{rewritelemma}. Using \eqref{claspN} to produce neutral ladders on bottom, we need to consider the LHS of the following diagram.
\begin{equation} \label{gammaEasy2} {
\labellist
\small\hair 2pt
 \pinlabel {$\l^- + \nu$} [ ] at 57 21
 \pinlabel {$\l + \mu$} [ ] at 70 139
 \pinlabel {$\l^- + \nu$} [ ] at 257 22
 \pinlabel {$\l + \mu$} [ ] at 269 115
\tiny\hair 2pt
 \pinlabel {$y_2$} [ ] at 31 38
 \pinlabel {$t_1$} [ ] at 55 38
 \pinlabel {$\dots$} [ ] at 92 38
 \pinlabel {$t_l$} [ ] at 113 38
 \pinlabel {$a$} [ ] at 137 21
 \pinlabel {$t_1$} [ ] at 41 85
 \pinlabel {$\dots$} [ ] at 68 85
 \pinlabel {$t_l$} [ ] at 90 85
 \pinlabel {$y_2$} [ ] at 113 85
 \pinlabel {$y_1$} [ ] at 41 123
 \pinlabel {$s_1$} [ ] at 65 123
 \pinlabel {$\dots$} [ ] at 92 123
 \pinlabel {$s_l$} [ ] at 113 123
 \pinlabel {$y_2$} [ ] at 137 123
 \pinlabel {$y_2$} [ ] at 243 40
 \pinlabel {$t_1$} [ ] at 267 40
 \pinlabel {$\dots$} [ ] at 296 40
 \pinlabel {$t_l$} [ ] at 317 40
 \pinlabel {$a$} [ ] at 340 22
 \pinlabel {$y_2$} [ ] at 243 100
 \pinlabel {$y_1$} [ ] at 267 100
 \pinlabel {$s_1$} [ ] at 293 100
 \pinlabel {$\dots$} [ ] at 319 100
 \pinlabel {$s_l$} [ ] at 340 100
\endlabellist
\centering
\ig{1}{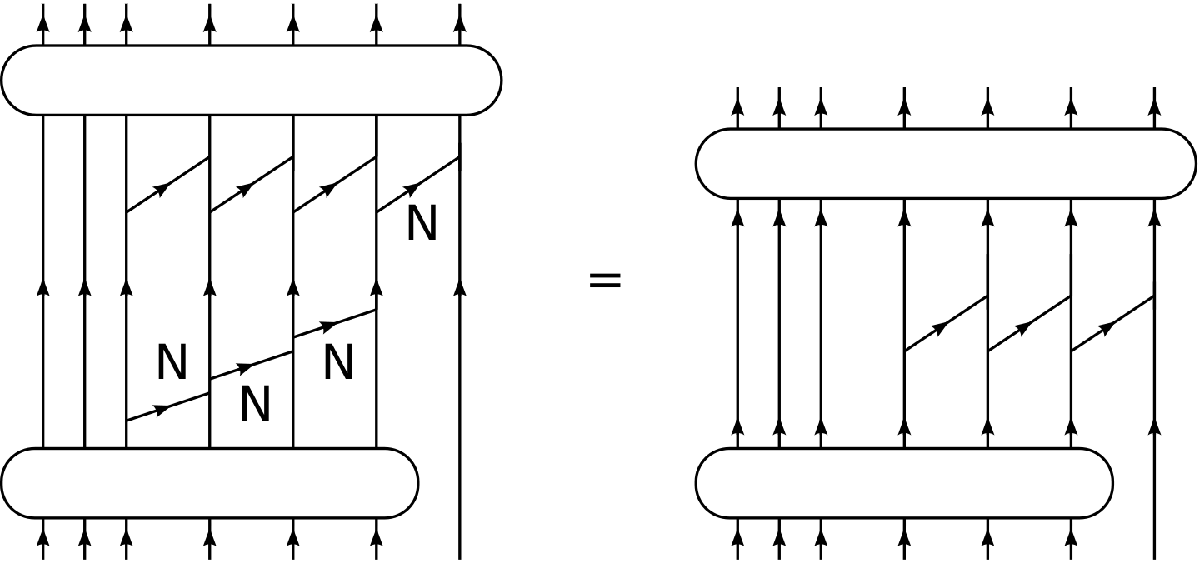}
} \end{equation}
Now we can apply \eqref{R3special} repeatedly, until we reach the RHS of \eqref{gammaEasy2}. In doing so, we ignore terms with an inward ladder on top (since they are orthogonal to the $(\l+\mu)$-clasp), and allow the $(\l+\mu)$-clasp to absorb all neutral ladders on top. As noted below \eqref{R3special}, this operation has coefficient $1$. The reader can once again confirm that the result is $E_{\s}$. \end{proof}

This proposition is sufficient to compute that the coefficients $\g_{\l, \mu, \nu}$ are either $1$ or $0$ for $n \le 3$, and most of the time for $n=4$. Before continuing with the arduous computation of $\g$ for $n \ge 4$, let us unravel the recursions for $n=2,3$.

\begin{subequations} When $n=2$, there is only one weight $\mu \in \Om^-(a)$ for any $a$, namely $(01) \in \Om^-(1)$. The recursion \eqref{recursive1} simplifies to \begin{equation}
\label{recursivesl2} \k_{(b), (01)} = [2] - \frac{1}{\k_{(b-1), (01)}} \end{equation} so long as $b \ge 2$, and \begin{equation} \k_{(1), (01)} = [2]. \end{equation} Note that the sum
over $\nu$ in \eqref{recursive1} has a single term $\nu = (01)$ when $b \ge 2$, and no terms when $b =1$ (and the intersection form $\k_{(b), (01)}$ is not defined for $b=0$, since $\l +
\mu$ is not dominant). Moreover, both $\k_{\l^-, \mu^-}$ and $\k_{\l^- + \nu, \s}$ are just $+1$. It is easy to verify that 
\begin{equation} \k_{(b), (01)} = \frac{[b+1]}{[b]} \end{equation} solves the recursion.
\end{subequations}

When $n=3$, there are four recursive formulas, for $\mu = (010)$ or $(001)$ in $\Om^-(1)$, and for $\mu = (101)$ or $(011)$ in $\Om^-(2)$. It is always the case that $\k_{\l^-, \mu^-} = 1$.

\begin{subequations}
When $\mu=(010)$ then $x_k = 1$ and $\nu \in \Om^-(1)$ can be either $(010)$ or $(001)$. In the former case, $\s = (100)$. In the latter case, $\s$ is not in $\Om(1)$. Thus we have
\begin{equation} \k_{(b,c), (010)} = [2] - \frac{1}{\k_{(b-1, c), (010)}}, \end{equation}
so long as $b \ge 2$. This is the same recursion as in the $n=2$ case, with solution
\begin{equation} \k_{(b,c), (010)} = \frac{[b+1]}{[b]}. \end{equation}	
Similarly, when $\mu = (101)$, $x_k = 2$, and $\nu \in \Om^-(2)$ can be either $(101)$ or $(011)$. In the former case, $\s = \w_2$, while in the latter case, $\s \notin \Om(2)$. Thus we have
\begin{equation} \k_{(b,c), (101)} = [2] - \frac{1}{\k_{(b,c-1), (101)}}, \end{equation}
for $c \ge 2$, with solution
\begin{equation} \k_{(b,c), (101)} = \frac{[c+1]}{[c]}. \end{equation}

When $\mu = (001)$ then $x_k = 2$, and $\nu \in \Om^-(2)$ can be either $(101)$ or $(011)$. In the former case, $\s = (010)$, and in the latter case, $\s = \w_1$. Thus we have
\begin{equation} \k_{(b,c), (001)} = [3] - \frac{\k_{(b+1,c-2), (010)}}{\k_{(b,c-1), (101)}} - \frac{1}{\k_{(b,c-1), (011)}}, \end{equation}
where the first fraction disappears when $c = 1$, and the second disappears when $b = 0$ (and the intersection form is not defined when $c=0$). Plugging in our known formulas from above, we get
\begin{equation} \k_{(b,c), (001)} = [3] - \frac{[b+2]}{[b+1]}\frac{[c-1]}{[c]} - \frac{1}{\k_{(b,c-1), (011)}}. \end{equation}
Similarly, when $\mu = (011)$ we end up with the recursive formula
\begin{equation} \k_{(b,c), (011)} = [3] - \frac{[c+2]}{[c+1]}\frac{[b-1]}{[b]} - \frac{1}{\k_{(b-1,c), (001)}}, \end{equation}
where the first fraction disappears when $b=1$, and the second disappears when $c=0$ (and the intersection form is not defined when $b=0$). These two formulas can be used to solve each other.
We encourage the reader to confirm that
\begin{equation} \k_{(b,c), (001)} = \frac{[c+1]}{[c]}\frac{[b+c+2]}{[b+c+1]} \end{equation}
and
\begin{equation} \k_{(b,c), (011)} = \frac{[b+1]}{[b]}\frac{[b+c+2]}{[b+c+1]} \end{equation}
solves the recursion, which is a (nontrivial) exercise in manipulating quantum numbers. For purposes of sanity, the reader may want to confirm this first when $q=1$, so that quantum numbers are ordinary numbers.
\end{subequations}

Note that the recursive formula for $\k_{(b,c), (010)}$ involved computing other coefficients $\k_{\l, \mu}$ of the same kind, but not coefficients of other kinds, so that this recursion
could be solved before moving on to more difficult recursions. This appears to be a general phenomenon: there appears to be a preorder for which the recursive formula for $\k_{\l, \mu}$
only depends on $\k_{?, \xi}$ for $\xi \le \mu$ in the preorder. Looking in \eqref{recursive1}, the recursion formula for $\k_{\l, \mu}$ only involves $\k_{?, \xi}$ for $\xi = \mu^-$,
$\nu$, or $\s$. The following loosely-worded claim is a first attempt at understanding this preorder.

\begin{claim} One has $\langle \mu, \rho \rangle \le \langle \xi, \rho \rangle$ for $\xi = \mu^-$, $\nu$, or $\s$, so long as $\k_{?, \xi}$ plays an ``interesting" role in
\eqref{recursive1}. \end{claim}
	
\begin{proof} It is easy to see that $\langle \mu, \rho \rangle$ is equal to $-1$, plus $1$ if $y_{k+1}<n$, plus $1$ if $y_1 > 0$. That is, the pairing of a 01 sequence with $\rho$ can
only be decreased by: putting a 0 at the start, when before there was a 1; putting a 1 at the end, when before there was a 0. Neither of these operations can be achieved by $\nu$ or
$\s$, lest $\g$ vanish for trivial reasons. The only time that $\mu^-$ can put a 0 at the start is when $\mu = \w_a$, which is not an interesting case. \end{proof}

Thus, for instance, pairing with $-\rho$ seems like an interesting preorder. However, the problem is that \eqref{recursive1} uses $\g_{\l, \mu, \nu}$, whose computation may secretly
involve other $\k_{?, \xi}$. Thus, it is not known whether pairing with $-\rho$ is a preorder governing the recursive formula. Moreover, computations for $n=4$ show that an even stronger
preorder is playing a role, but we have not put in the effort to make this preorder precise.

When $n=4$, it remains to compute the coefficient $\g_{\l, \mu, \nu}$ when $\mu = (0101)$ and $\nu = (1101)$ or $\nu = (0111)$, as the other coefficients $\g$ are either $1$ or $0$.
Thus, many $\k_{\l, \mu}$ can be computed using recursive formulas that do not involve any of the unknown coefficients $\g$. We list them here.

\begin{subequations}
\begin{equation} \k_{(b,c,d), (0100)} = [2] - \frac{1}{\k_{(b-1, c, d), (0100)}}. \end{equation}
\begin{equation} \k_{(b,c,d), (1010)} = [2] - \frac{1}{\k_{(b, c-1, d), (1010)}}. \end{equation}
\begin{equation} \k_{(b,c,d), (1101)} = [2] - \frac{1}{\k_{(b, c, d-1), (1101)}}. \end{equation}

\begin{equation} \k_{(b,c,d), (0010)} = [3] - \frac{[b+2]}{[b+1]}\frac{[c-1]}{[c]} - \frac{1}{\k_{(b, c-1, d), (0110)}}. \end{equation}
\begin{equation} \k_{(b,c,d), (0110)} = [3] - \frac{[c+2]}{[c+1]}\frac{[b-1]}{[b]} - \frac{1}{\k_{(b-1, c, d), (0010)}}. \end{equation}
\begin{equation} \k_{(b,c,d), (1011)} = [3] - \frac{[d+2]}{[d+1]}\frac{[c-1]}{[c]} - \frac{1}{\k_{(b, c-1, d), (1001)}}. \end{equation}
\begin{equation} \k_{(b,c,d), (1001)} = [3] - \frac{[c+2]}{[c+1]}\frac{[d-1]}{[d]} - \frac{1}{\k_{(b, c, d-1), (1011)}}. \end{equation}	
	
\begin{equation} \k_{(b,c,d), (0001)} = [4] - \frac{\k_{(b,c+1,d-2), (0010)}}{\k_{(b,c,d-1), (1101)}} - \frac{\k_{(b+1,c-1,d-1), (0100)}}{\k_{(b,c,d-1), (1001)}} - \frac{1}{\k_{(b,c,d-1), (0111)}}. \end{equation}
\begin{equation} \k_{(b,c,d), (0111)} = [4] - \frac{\k_{(b-2,c+1,d), (1011)}}{\k_{(b-1,c,d), (0100)}} - \frac{\k_{(b-1,c-1,d+1), (1101)}}{\k_{(b-1,c,d), (0010)}} - \frac{1}{\k_{(b-1,c,d), (0001)}}. \end{equation}
\end{subequations}

We claim that these recursions are solved by \eqref{n=4LIF}. This is a long and difficult computation involving manipulations with quantum numbers; the author verified it by computer.

Except for the last two, these recursions look exactly like the recursions from $\sl_3$. In fact, it seems that placing a $1$ at the start of each $\mu$, or a $0$ at the end of each
$\mu$, will not change the recursion equation. This is a consequence of Conjecture \ref{MainConj}.

\begin{prop} Assume Conjecture \ref{MainConj}. Let $\l = (b_1, b_2, \ldots, b_{n-1})$ be a dominant $\sl_n$ weight, and let $\mu$ be a 01-sequence of length $n$. Let $\mu 0$ and $1\mu$
be the sequences obtained from $\mu$ by adding a zero at the end or a 1 at the start. Let $\l c = (b_1, \ldots, b_{n-1}, c)$ and $c \l = (c, b_1, \ldots, b_{n-1})$ be dominant
$\sl_{n+1}$ weights, for $c \ge 0$. Then \begin{subequations} \begin{equation} \k_{\l, \mu} = \k_{\l c, \mu 0}, \end{equation} \begin{equation} \k_{\l, \mu} = \k_{c \l, 1 \mu}.
\end{equation} \end{subequations} \end{prop}

\begin{proof} As discussed in Claim \ref{minusoneclaim}, the roots in $\Phi(\mu)$ are in bijection with \emph{inversions}, or pairs consisting of a 0 before a 1 in $\mu$. So $\Phi(\mu 0)$
and $\Phi(1 \mu)$ are in a straightforward bijection with $\Phi(\mu)$, since adding a 1 to the start or a 0 to the end will not produce any new inversions, though it may reindex them. For
each root in $\Phi(\mu)$, it is easy to see that pairing it against $\l$ (resp. $\l + \mu$, $\rho$) will give the same result as pairing the corresponding root in $\Phi(\mu 0)$ against
$\l c$ (resp. $\l c + \mu 0$, $\rho$). The same is true for $1 \mu$ and $c \l$. \end{proof}

\subsection{Deriving a recursive formula: part II}
\label{subsec-recursive2}

Let us compute some of the remaining coefficients $\g$, to illustrate the subtleties involved in a general computation. First, a lemma.

\begin{lemma} \label{sillylemma} Suppose that $\mu \in \Om^-(a)$ has more than one 1-string, and consider $\mu^- \in \Om(\a)$, obtained by removing the last 1-string from $\mu$. If $\s^-
\in \Om(\a)$ is bigger than $\mu^-$ in the dominance order, then one can add back in the 1-string to obtain $\s \in \Om(a)$ which is bigger than $\mu$, and whose last 1-string is the
same size as that of $\mu$. \end{lemma}

\begin{proof} Let $\{x_i\}$ and $\{y_i\}$ and $k$ be determined from $\mu$ as usual. Then $\mu^-$ has a final 0-string from positions $y_k+1$ through $n$. Since $\s^- > \mu^-$, it too must have zeroes in positions $y_k + 1$ through $n$. Thus one can reinsert ones in positions $x_k + 1$ through $y_{k+1}$, and obtain a weight $\s \in \Om(a)$ as desired. \end{proof}

The purpose of this silly lemma is as follows: consider an elementary light ladder $E_{\mu}$, and ignore the last rung to obtain $E_{\mu^-}$. If $\mu^-$ is then involved in a diagram like
the LHS of \eqref{gammaterm}, and is replaced with $E_{\s^-}$ like the RHS of \eqref{gammaterm}, then one can reinsert the last rung to obtain an elementary light ladder $E_{\s}$. This
will become clear in practice below.

Now we compute $\g_{\l, (0101), (0111)}$. In this case, $\s = \w_a = (1110)$. We use $\t$ for $\l - \w_3 - \w_1$.
\begin{equation} \label{computation1} {
\labellist
\small\hair 2pt
 \pinlabel {$\l - \w_3 + \nu$} [ ] at 40 29
 \pinlabel {$\l - \w_3$} [ ] at 28 96
 \pinlabel {$\l+\mu$} [ ] at 52 156
 \pinlabel {$\l - \w_3 + \nu$} [ ] at 182 30
 \pinlabel {$\t$} [ ] at 157 96
 \pinlabel {$\l + \mu$} [ ] at 192 156
 \pinlabel {$\l - \w_3 + \nu$} [ ] at 352 12
 \pinlabel {$\t$} [ ] at 326 46
 \pinlabel {$\t + \eta$} [ ] at 339 91
 \pinlabel {$\t$} [ ] at 326 135
 \pinlabel {$\l + \mu$} [ ] at 361 169
\tiny\hair 2pt
 \pinlabel {$0$} [ ] at 54 47
 \pinlabel {$4$} [ ] at 78 47
 \pinlabel {$1$} [ ] at 54 80
 \pinlabel {$3$} [ ] at 78 80
 \pinlabel {$2$} [ ] at 102 80
 \pinlabel {$1$} [ ] at 54 112
 \pinlabel {$3$} [ ] at 78 112
 \pinlabel {$0$} [ ] at 54 141
 \pinlabel {$2$} [ ] at 78 141
 \pinlabel {$4$} [ ] at 102 141
 \pinlabel {$0$} [ ] at 194 47
 \pinlabel {$4$} [ ] at 219 47
 \pinlabel {$1$} [ ] at 194 91
 \pinlabel {$3$} [ ] at 219 91
 \pinlabel {$2$} [ ] at 244 91
 \pinlabel {$0$} [ ] at 194 141
 \pinlabel {$2$} [ ] at 219 141
 \pinlabel {$4$} [ ] at 244 141
 \pinlabel {$0$} [ ] at 363 31
 \pinlabel {$4$} [ ] at 386 31
 \pinlabel {$1$} [ ] at 363 62
 \pinlabel {$3$} [ ] at 386 62
 \pinlabel {$2$} [ ] at 411 62
 \pinlabel {$1$} [ ] at 363 120
 \pinlabel {$1$} [ ] at 386 120
 \pinlabel {$4$} [ ] at 411 152
 \pinlabel {$0$} [ ] at 363 152
 \pinlabel {$2$} [ ] at 386 152
 \pinlabel {\large $\displaystyle - \sum_\eta \frac{1}{\k_{\t, \eta}}$} [ ] at 276 93
\endlabellist
\centering
\ig{1}{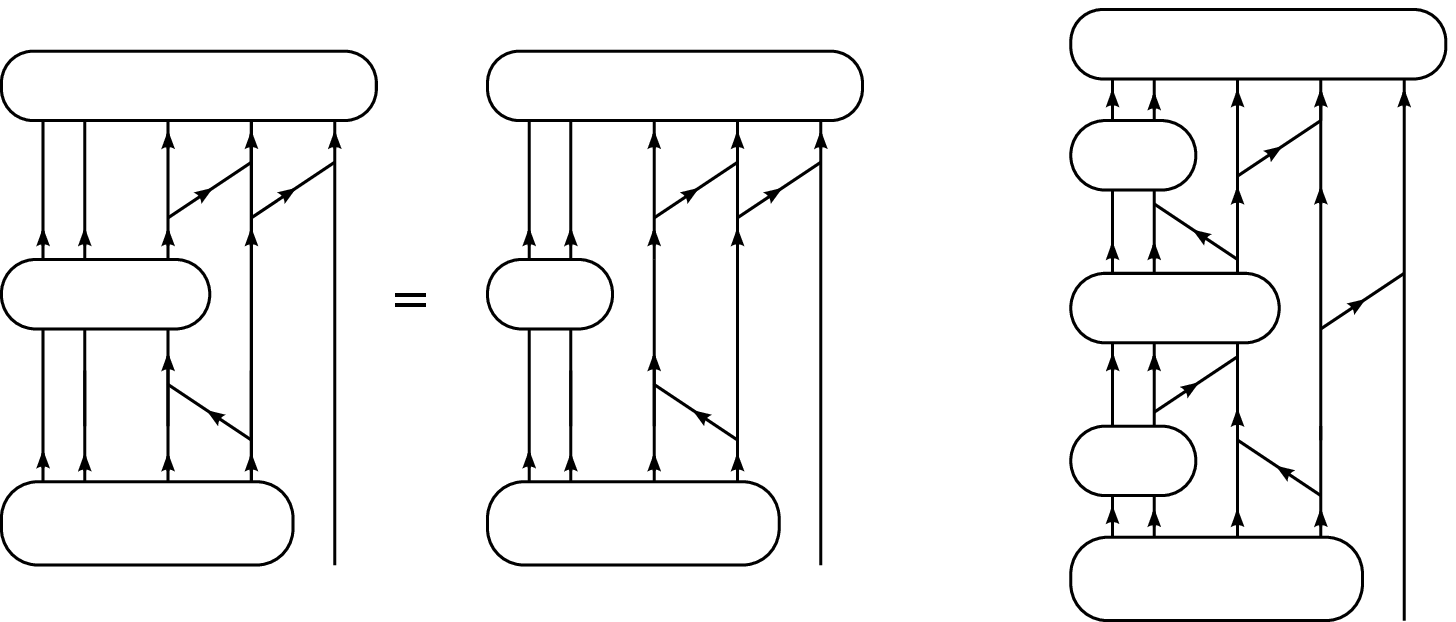}
} \end{equation}
Upon resolving the middle clasp using \eqref{eq:tripleclasp}, we obtain the RHS. The first diagram on the RHS is easy to resolve using associativity and the rung swap: the result is $[2]$ times $E_{\s}$.

Now let us analyze each term in the sum: here $\eta \in \Om^-(1)$ and $\t + \eta$ is dominant. Having shifted the rightmost rung $R$ out of the way, the rest of the diagram splits into
two halves, each of which looks like the LHS of \eqref{gammaterm} (possibly upside-down). Thus the top half yields coefficient $\g_{\l - \w_3, \mu^-, \eta}$ times $E_{\s_1}$ for some
$\s_1$, while the bottom half yields coefficient $\g_{\l - \w_3, \nu, \eta} E_{\s_2}$. Lucky, both of these $\g$ coefficients are $1$ or trivially $0$ by Proposition \ref{gammaEasyprop}. 

Moreover, placing the rung $R$ below $E_{\s_1}$ gives another elementary light ladder $E_{\s_1^+}$, thanks to Lemma \ref{sillylemma}. Now, the entire diagram looks like the LHS of \eqref{gammaterm} again, yielding a factor of $\g_{\l - \w_1, \s_1^+, \s_2}$. The resulting diagram is, finally, $E_{\s}$.

Thus we obtain the recursive formula below.
\begin{equation} \label{grecursionimprecise} \g_{\l, (0101), (0111)} = [2] - \sum_{\eta} \frac{1}{\k_{\t, \eta}} \g_{\l-\w_1, \s_1^+, \s_2}. \end{equation}
To make this more precise: $\eta$ ranges over $\Om^-(1)$ such that $\t + \eta$ is dominant. If $\eta = (0100)$ then $\s_2 = (0111) + \w_1 - \eta = (1011)$, and $\s_1 = (0100) + \w_1 - \eta = \w_1$, so that $\s_1^+ = (1001)$. Then $\g_{\l - \w_1, (1001), (1011)} = 1$ by Proposition \ref{gammaEasyprop}, so long as things are all dominant. If $\eta = (0010)$ or $(0001)$ then $\s_1$ is not in $\Om(1)$, so the term does not appear. Thus we have
\begin{equation} \label{grecursionprecise} \g_{(b,c,d), (0101), (0111)} = [2] - \frac{1}{\k_{(b-1, c, d-1), (0100)}} = \frac{[b+1]}{[b]}. \end{equation}
Here, we have used the formula for $\k_{\t, (0100)}$ from \eqref{n=4LIF}, which is already known.

Now we compute $\g_{\l, (0101), (1101)}$, where $\l^- = \l - \w_3$ and $\s = (0110)$. Here a real complication arises. One wants to take the same approach as before: resolve an interior clasp and reduce to previously computed coefficients.  However, the label on the rightmost strand in the $\l^-$-clasp on top, $x_{k-1}$ from $\mu$, differs from the label on bottom, $s_l$ from $\nu$. Following Remark \ref{notsametopbottom}, we must apply a neutral ladder (and involve a new strand) in order to resolve the $\l^-$-clasp.

\begin{equation} \label{computation2} {
\labellist
\small\hair 2pt
 \pinlabel {$\l^- + \nu$} [ ] at 74 253
 \pinlabel {$\l^-$} [ ] at 63 319
 \pinlabel {$\l + \mu$} [ ] at 86 379
 \pinlabel {$\l^- + \nu$} [ ] at 220 238
 \pinlabel {$\l^-$} [ ] at 213 313
 \pinlabel {$\l + \mu$} [ ] at 235 392
 \pinlabel {$\l^- + \nu$} [ ] at 71 67
 \pinlabel {$\l^{--}$} [ ] at 47 124
 \pinlabel {$\l+\mu$} [ ] at 82 192
 \pinlabel {$\l^-+\nu$} [ ] at 262 12
 \pinlabel {$\l^{--}$} [ ] at 236 48
 \pinlabel {$\l^{--}+\eta$} [ ] at 249 90
 \pinlabel {$\l^{--}$} [ ] at 237 134
 \pinlabel {$\l + \mu$} [ ] at 270 192
 \pinlabel {$- \displaystyle \sum_{\eta} \frac{1}{\k_{\l^{--}, \eta}}$} [ ] at 172 111
 \tiny
 \pinlabel {$2$} [ ] at 136 254
 \pinlabel {$2$} [ ] at 284 254
 \pinlabel {$2$} [ ] at 131 80
 \pinlabel {$2$} [ ] at 321 28
 \pinlabel {$0$} [ ] at 89 364
 \pinlabel {$2$} [ ] at 112 364
 \pinlabel {$4$} [ ] at 138 364
 \pinlabel {$1$} [ ] at 88 336
 \pinlabel {$3$} [ ] at 112 336
 \pinlabel {$3$} [ ] at 88 303
 \pinlabel {$2$} [ ] at 88 272
 \pinlabel {$4$} [ ] at 112 272
 
 \pinlabel {$3$} [ ] at 212 376
 \pinlabel {$0$} [ ] at 235 376
 \pinlabel {$2$} [ ] at 262 376
 \pinlabel {$4$} [ ] at 283 376
 \pinlabel {$1$} [ ] at 235 348
 \pinlabel {$3$} [ ] at 262 348
 \pinlabel {$1$} [ ] at 211 331
 \pinlabel {$3$} [ ] at 235 331
 \pinlabel {$1$} [ ] at 212 299
 \pinlabel {$3$} [ ] at 235 299
 \pinlabel {$2$} [ ] at 235 256
 \pinlabel {$4$} [ ] at 262 256
 
 \pinlabel {$3$} [ ] at 60 176
 \pinlabel {$0$} [ ] at 83 176
 \pinlabel {$2$} [ ] at 109 176
 \pinlabel {$4$} [ ] at 132 176
 \pinlabel {$1$} [ ] at 83 147
 \pinlabel {$3$} [ ] at 109 147
 \pinlabel {$1$} [ ] at 60 138
 \pinlabel {$3$} [ ] at 83 121
 \pinlabel {$2$} [ ] at 73 84
 \pinlabel {$4$} [ ] at 109 84
 
 \pinlabel {$3$} [ ] at 249 176
 \pinlabel {$0$} [ ] at 272 176
 \pinlabel {$2$} [ ] at 296 176
 \pinlabel {$4$} [ ] at 322 176
 \pinlabel {$1$} [ ] at 265 150
 \pinlabel {$1$} [ ] at 293 150
 \pinlabel {$3$} [ ] at 272 131
 \pinlabel {$3$} [ ] at 285 75
 \pinlabel {$3$} [ ] at 272 65
 \pinlabel {$2$} [ ] at 260 33
 \pinlabel {$4$} [ ] at 296 33
 
\endlabellist
\centering
\ig{1}{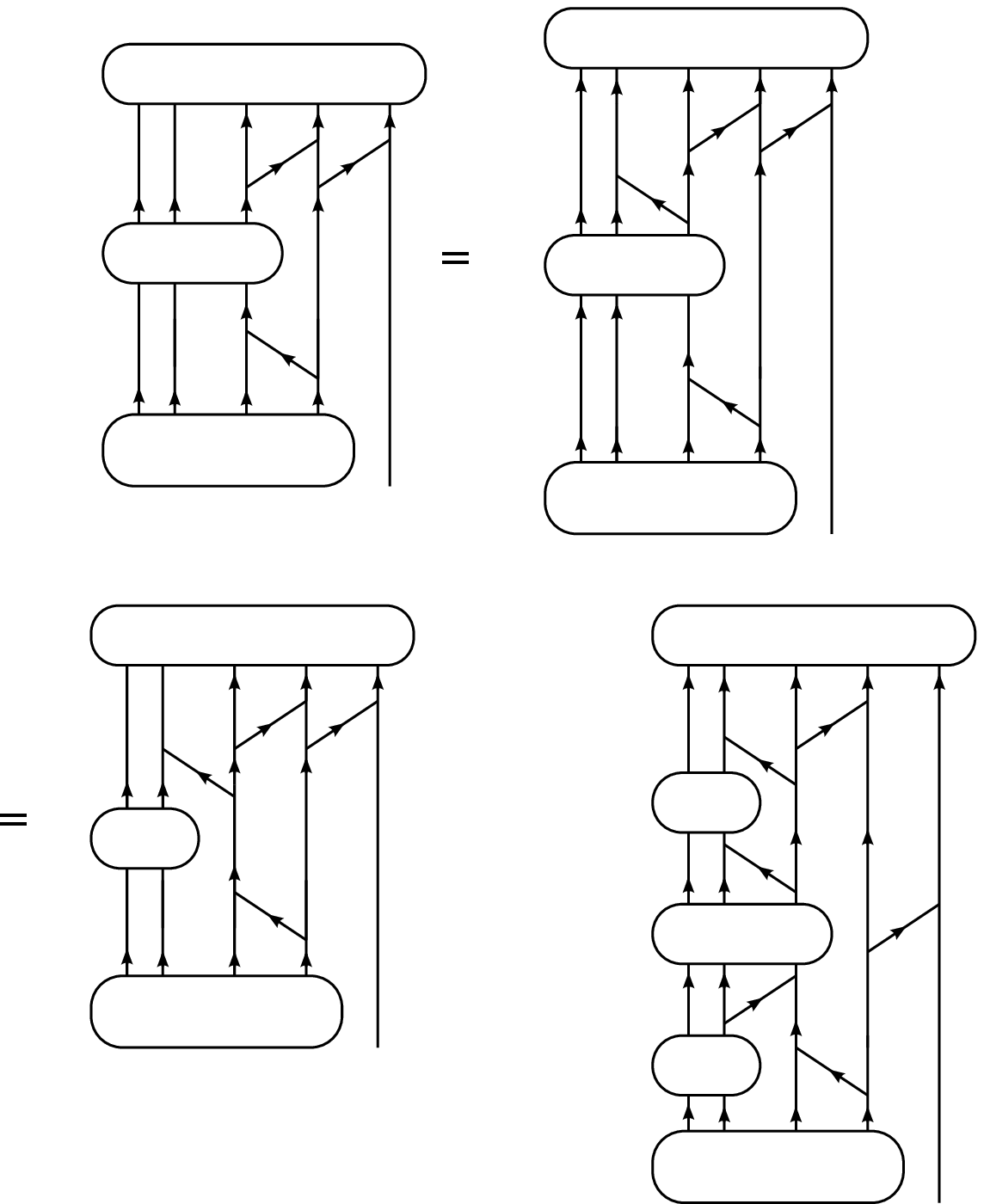}
} \end{equation}
Here, $\l^{--} = \l - 2 \w_3$.

Once again, the first diagram on the bottom row can be resolved directly, to obtain $E_{\s}$ with coefficient $1$. We leave this to the reader. In the second diagram, the final rung $R$
was shunted out of the way, and the remainder splits into two halves. As before, the bottom half is the LHS of \eqref{gammaterm} upside-down, and can be resolved with a coefficient of
$\g_{\l^-, \nu, \eta}$, and this coefficient is $1$ or $0$. In fact, $\g_{\l^-, \nu, \eta} = 0$ for $\eta \in \Om^-(3)$ unless $\eta = (1101) = \nu$, so that there is only one term in the sum. However, the top half is not of the form \eqref{gammaterm}, because the extra neutral ladder stuffs things up!

In a general computation of $\g_{\l, \mu, \nu}$, such neutral ladders will inevitably appear. Thus we should define the following generalization of the coefficients $\g$, for any $c \in \{1, \ldots, n-1\}$.
\begin{equation} \label{gammaterm2} {
 \labellist
 \tiny\hair 2pt
 \pinlabel {\small $\l - \w_b + \nu$} [ ] at 59 10
 \pinlabel {\small $\l - \w_c$} [ ] at 47 152
 \pinlabel {\small $\l+\mu$} [ ] at 68 287
 \pinlabel {\small $\l - \w_b+\nu$} [ ] at 307 80
 \pinlabel {\small $\l+\mu$} [ ] at 318 172
 \pinlabel {\small $\overline{E}_{\nu}$} [ ] at 149 61
 \pinlabel {\small $E_{\mu}$} [ ] at 149 243
 \pinlabel {\small $E_{\s}$} [ ] at 407 126
 \pinlabel {$a$} [ ] at 137 150
 \pinlabel {$c$} [ ] at 113 150
 \pinlabel {$x_k$} [ ] at 114 215
 \pinlabel {$b$} [ ] at 114 88
 \pinlabel {\small $\g_{\l, \mu, \nu, c}$} [ ] at 219 127
 \pinlabel {$c$} [ ] at 28 217
 \pinlabel {$c$} [ ] at 41 83
 \endlabellist
 \centering
 \ig{.8}{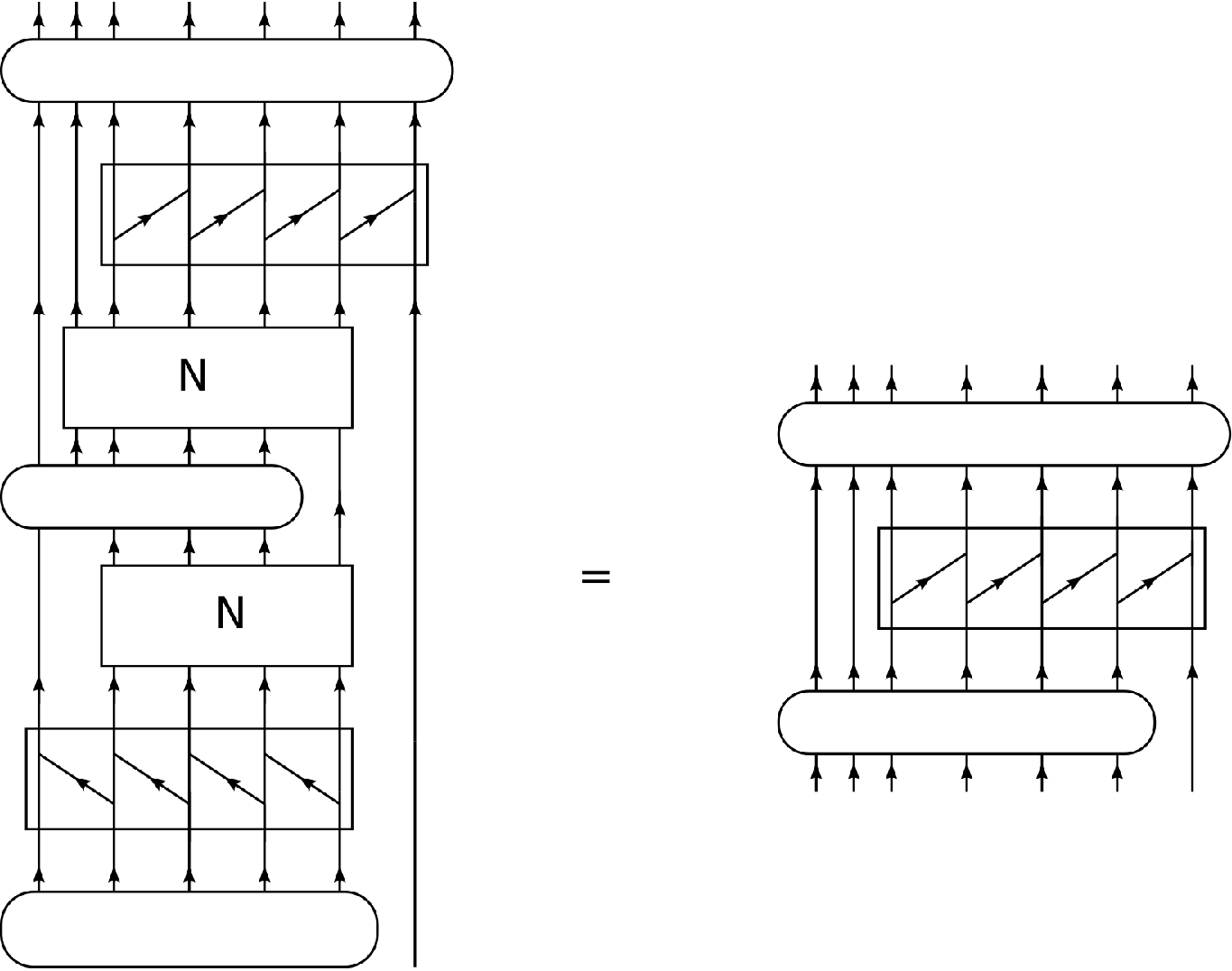}
} \end{equation}

Let us explain this diagram. Here, $\mu \in V_{\w_a}$ and $\nu \in V_{\w_b}$ are two weights in arbitrary fundamental representations (unlike in \eqref{gammaterm} where $b$ was required
to equal $x_k$, coming from $\mu$). Meanwhile, $c$ is a third arbitrary label which must occur in $\l$, and $\l^- = \l - \w_c$. The $N$ boxes represent neutral ladders which bring $c$ to
its appropriate location. Because $c$ appears inside $\l$, it is either one of the inputs involved in the elementary light ladder $E_{\l, \mu}$, or it is one of the irrelevant uprights
to the left, and similarly for $E_{\l^-, \nu}$. In the exemplary diagram above, one of each occurs: it appears that $c$ is one of the irrelevant uprights for $E_{\l, \mu}$, but is one of
the relevant (upside-down) inputs to $E_{\l^-, \nu}$. Finally, $\s = \mu + \w_b - \nu$. When $\s$ is not a weight in $V_{\w_a}$, we say that $\g_{\l, \mu, \nu, c} = 0$ for trivial
reasons.

When $c = b = x_k$, the neutral ladders must be just identity maps, and $\g_{\l, \mu, \nu, c} = \g_{\l, \mu, \nu}$. When the neutral ladders are not identity maps, it is much more
difficult to determine a recursive formula for $\g_{\l, \mu, \nu, c}$. Even worse, we can no longer apply Lemma \ref{sillylemma} in the same way, to assume that the rung $R$ will
recombine to form an elementary light ladder, as we did in the earlier analogous computation. Overall, this requires a much more sophisticated computation.

Let us return to the computation of $\g_{\l, (0101), (1101)}$. We need to resolve the diagram on the bottom right of \eqref{computation2} when $\eta = \nu = (1101)$. Thus we need to compute the following.
\begin{equation} \label{computation3} {
\labellist
\small\hair 2pt
 \pinlabel {$\l^- + \nu$} [ ] at 41 9
 \pinlabel {$\l^{--}$} [ ] at 18 67
 \pinlabel {$\l+\mu$} [ ] at 51 135
 \tiny
 \pinlabel {$3$} [ ] at 29 118
 \pinlabel {$0$} [ ] at 54 118
 \pinlabel {$2$} [ ] at 78 118
 \pinlabel {$4$} [ ] at 103 118
 \pinlabel {$1$} [ ] at 54 91
 \pinlabel {$3$} [ ] at 76 91
 \pinlabel {$2$} [ ] at 103 11
 \pinlabel {$1$} [ ] at 29 83
 \pinlabel {$3$} [ ] at 29 52
 \pinlabel {$3$} [ ] at 54 52
 \pinlabel {$2$} [ ] at 29 28
 \pinlabel {$4$} [ ] at 54 28
\endlabellist
\centering
\ig{1}{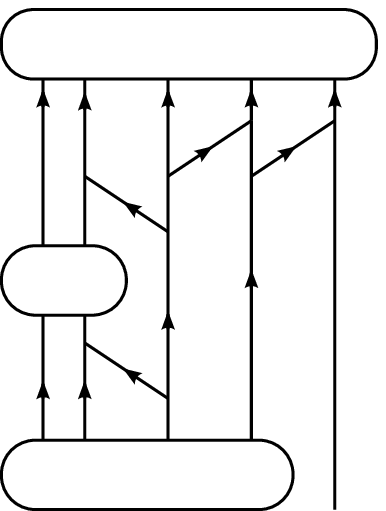}
} \end{equation} 
Now we resolve the $\l^{--}$-clasp, which requires another neutral ladder to be placed, and we get
\begin{equation} \label{computation4} {
\labellist
\small\hair 2pt
 \pinlabel {$\displaystyle - \sum_{\eta} \frac{1}{\k_{\l - 3 \w_3, \eta}}$} [ ] at 153 92
 \pinlabel {$\l^- + \nu$} [ ] at 50 38
 \pinlabel {$\l^{---}$} [ ] at 18 96
 \pinlabel {$\l+\mu$} [ ] at 64 163
 \pinlabel {$\l^-+\nu$} [ ] at 246 12
 \pinlabel {$\l^{---}$} [ ] at 211 47
 \pinlabel {$\l^{---}+\eta$} [ ] at 221 90
 \pinlabel {$\l^{---}$} [ ] at 209 135
 \pinlabel {$\l+\mu$} [ ] at 256 191
 \tiny
  \pinlabel {$4$} [ ] at 127 147
  \pinlabel {$2$} [ ] at 103 147
  \pinlabel {$0$} [ ] at 78 147
  \pinlabel {$3$} [ ] at 53 147
  \pinlabel {$3$} [ ] at 30 147
  \pinlabel {$2$} [ ] at 127 38
  \pinlabel {$3$} [ ] at 103 93
  \pinlabel {$3$} [ ] at 78 93
  \pinlabel {$3$} [ ] at 53 93
  \pinlabel {$1$} [ ] at 30 111
  \pinlabel {$4$} [ ] at 78 55
  \pinlabel {$2$} [ ] at 53 55
 \pinlabel {$4$} [ ] at 319 175
 \pinlabel {$2$} [ ] at 297 175
 \pinlabel {$0$} [ ] at 272 175
 \pinlabel {$3$} [ ] at 247 175
 \pinlabel {$3$} [ ] at 221 175
 \pinlabel {$2$} [ ] at 319 12
 \pinlabel {$3$} [ ] at 297 125
 \pinlabel {$3$} [ ] at 272 125
 \pinlabel {$3$} [ ] at 247 125
 \pinlabel {$1$} [ ] at 221 149
 \pinlabel {$3$} [ ] at 247 61
 \pinlabel {$4$} [ ] at 272 31
 \pinlabel {$2$} [ ] at 247 31
\endlabellist
\centering
\ig{1}{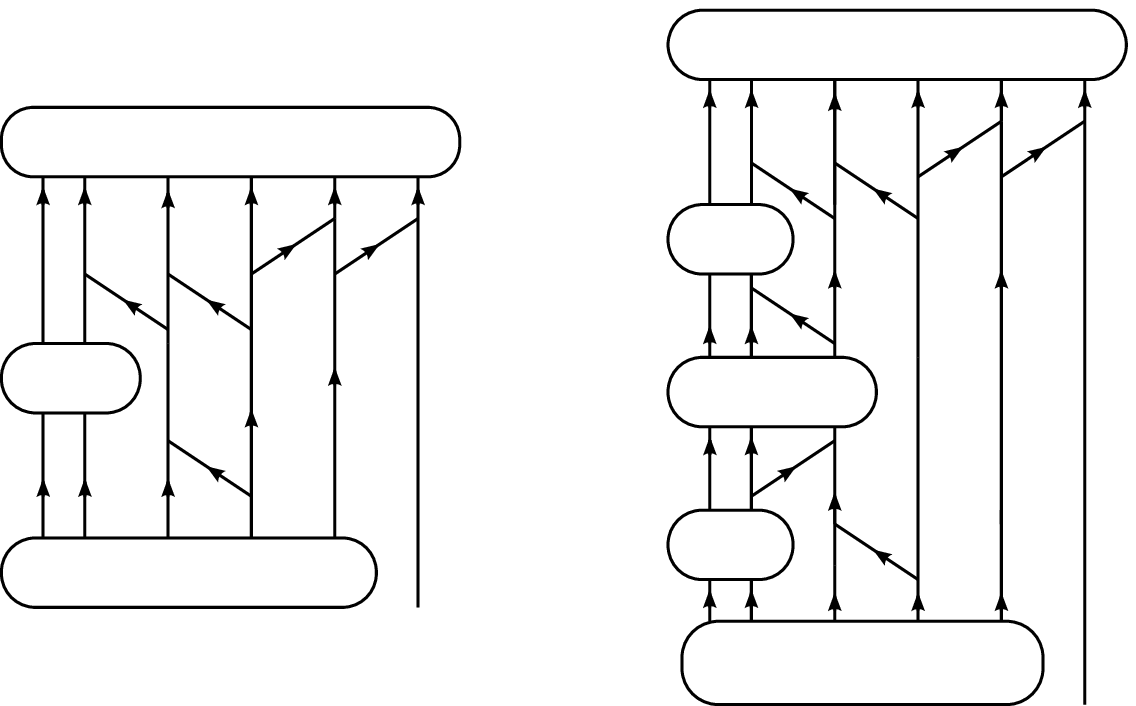}
}. \end{equation}
This time the first diagram is zero. We can apply \eqref{R3special} and move the non-neutral rung to the top of the diagram, where it vanishes by \eqref{claspout}. In the sum over $\eta \in \Om^-(3)$, once again, the lower half is easy to resolve using \eqref{gammaterm}, and the only non-vanishing case is when $\eta = \nu = (1101)$. Thus it remains to compute
\begin{equation} \label{computation5} {
\labellist
\small\hair 2pt
 \pinlabel {$\l + \mu$} [ ] at 64 135
 \pinlabel {$\l^-+\nu$} [ ] at 52 11
 \pinlabel {$\l^{---}$} [ ] at 19 68
\tiny
 \pinlabel {$4$} [ ] at 127 119
 \pinlabel {$2$} [ ] at 102 119
 \pinlabel {$0$} [ ] at 78 119
 \pinlabel {$3$} [ ] at 55 119
 \pinlabel {$3$} [ ] at 29 119
 \pinlabel {$2$} [ ] at 127 8
 \pinlabel {$3$} [ ] at 102 67
 \pinlabel {$3$} [ ] at 78 67
 \pinlabel {$3$} [ ] at 55 67
 \pinlabel {$1$} [ ] at 29 83
 \pinlabel {$4$} [ ] at 55 26
 \pinlabel {$2$} [ ] at 29 26
\endlabellist
\centering
\ig{1}{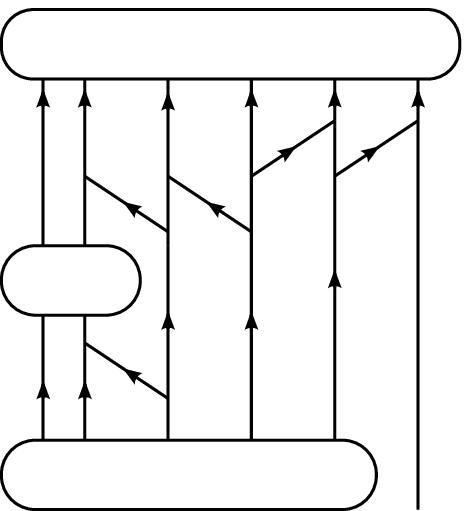}
}. \end{equation}
Resolve the $\l^{---}$-clasp and the pattern repeats. Now, the first term will always vanish by applying \eqref{R3special} and \eqref{claspout}, while a sum over lower terms will have only one term, iterating this procedure. Eventually, $\l - m \w_3$ is no longer dominant, and the process terminates.

We have shown that the bottom right diagram in \eqref{computation2} is actually zero, a fact which was not a priori obvious from examining the weights. Hence, $\g_{\l, (0101), (1101)} =
1$. Interesting patterns appeared in this computation, giving hope that something may be adapted to the general case.

We can now derive the final recursion relations for $n=4$. They are:

\begin{equation} \k_{(b,c,d), (0101)} = [3] \k_{(b,c,d-1), (0100)} - \frac{\k_{(b,c+1,d-2), (0110)}}{\k_{(b,c,d-1), (1101)}} - \frac{\k^2_{(b,c,d-1), (0100)}}{\k_{(b,c,d-1), (0111)}}.
\end{equation}
	
\begin{align} \k_{(b,c,d), (0011)} = \qbinom{4}{2} - \frac{\k_{(b+1, c-2, d+1), (0101)}}{\k_{(b,c-1,d), (1010)}}
	- \frac{\k_{(b-1, c-1, d+1), (1001)}}{\k_{(b,c-1,d), (0110)}}  - \frac{\k_{(b+1, c-1, d-1), (0110)}}{\k_{(b,c-1,d), (1001)}} \\
\nonumber \qquad \qquad	- \frac{\k_{(b-1, c, d-1), (1010)}}{\k_{(b,c-1,d), (0101)}} - \frac{1}{\k_{(b,c-1,d), (0011)}}. \end{align}

\noindent That \eqref{n=4LIF} will solve these recursive formulas is again a complicated equality amongst quantum numbers, which was verified by computer. More intuition will need to be
gained before the general form of these recursive formulas becomes obvious, but if the conjecture holds true, one can hope that they too admit a combinatorial description involving roots,
weights, and the Weyl group.

\bibliographystyle{plain}
\bibliography{mastercopy}

\begin{thebibliography}{10}

\bibitem{AST}
Henning~Haahr Andersen, Daniel Tubbenhauer, and Catharina Stroppel.
\newblock Cellular structures using $\textbf{U}_q$-tilting modules.
\newblock Preprint, 2015.
\newblock arXiv 1503.00224.

\bibitem{CautisClasp}
Sabin Cautis.
\newblock Clasp technology to knot homology via the affine {G}rassmannian.
\newblock {\em Math. Ann.}, 363(3-4):1053--1115, 2015.

\bibitem{CKM}
Sabin Cautis, Joel Kamnitzer, and Scott Morrison.
\newblock Webs and quantum skew {H}owe duality.
\newblock {\em Math. Ann.}, 360(1-2):351--390, 2014.

\bibitem{dCaMig05}
Mark Andrea~A. de~Cataldo and Luca Migliorini.
\newblock The {H}odge theory of algebraic maps.
\newblock {\em Ann. Sci. \'Ecole Norm. Sup. (4)}, 38(5):693--750, 2005.

\bibitem{EQuantumI}
Ben Elias.
\newblock Quantum {S}atake in type {A}: part {I}.
\newblock Preprint.
\newblock arXiv:1403.5570.

\bibitem{ELauda}
Ben Elias and Aaron~D. Lauda.
\newblock Trace decategorification of the {H}ecke category.
\newblock Preprint, 2015.
\newblock arXiv 1504.05267.

\bibitem{ELib}
Ben Elias and Nicolas Libedinsky.
\newblock Soergel bimodules for universal {C}oxeter groups.
\newblock {\em Transactions of the American Mathematical Society}, to appear.
\newblock arXiv:1401.2467.

\bibitem{EWGr4sb}
Ben Elias and Geordie Williamson.
\newblock Soergel calculus.
\newblock Preprint.
\newblock arXiv:1309.0865.

\bibitem{EWShadows}
Ben Elias and Geordie Williamson.
\newblock Kazhdan-{L}usztig conjectures and shadows of {H}odge theory.
\newblock Mathematische Arbeitstagung, 2013.
\newblock arXiv:1403.1650.

\bibitem{EWHodge}
Ben Elias and Geordie Williamson.
\newblock The {H}odge theory of {S}oergel bimodules.
\newblock {\em Ann. of Math. (2)}, 180(3):1089--1136, 2014.

\bibitem{Font12}
Bruce Fontaine.
\newblock Generating basis webs for {${\rm SL_n}$}.
\newblock {\em Adv. Math.}, 229(5):2792--2817, 2012.

\bibitem{GraLeh}
J.~J. Graham and G.~I. Lehrer.
\newblock Cellular algebras.
\newblock {\em Invent. Math.}, 123(1):1--34, 1996.

\bibitem{Jon3}
V.~F.~R. Jones.
\newblock Braid groups, {H}ecke algebras and type {${\rm II}_1$} factors.
\newblock In {\em Geometric methods in operator algebras ({K}yoto, 1983)},
  volume 123 of {\em Pitman Res. Notes Math. Ser.}, pages 242--273. Longman
  Sci. Tech., Harlow, 1986.

\bibitem{Kauf}
Louis~H. Kauffman.
\newblock State models and the {J}ones polynomial.
\newblock {\em Topology}, 26(3):395--407, 1987.

\bibitem{KhoKup}
Mikhail Khovanov and Greg Kuperberg.
\newblock Web bases for {${\rm sl}(3)$} are not dual canonical.
\newblock {\em Pacific J. Math.}, 188(1):129--153, 1999.

\bibitem{KhoLau10}
Mikhail Khovanov and Aaron~D. Lauda.
\newblock A categorification of quantum {${\rm sl}(n)$}.
\newblock {\em Quantum Topol.}, 1(1):1--92, 2010.

\bibitem{Kim03}
Dongseok Kim.
\newblock {\em Graphical calculus on representations of quantum {L}ie
  algebras}.
\newblock ProQuest LLC, Ann Arbor, MI, 2003.
\newblock Thesis (Ph.D.)--University of California, Davis.

\bibitem{Kim07}
Dongseok Kim.
\newblock Jones-{W}enzl idempotents for rank 2 simple {L}ie algebras.
\newblock {\em Osaka J. Math.}, 44(3):691--722, 2007.

\bibitem{Kupe}
Greg Kuperberg.
\newblock Spiders for rank {$2$} {L}ie algebras.
\newblock {\em Comm. Math. Phys.}, 180(1):109--151, 1996.

\bibitem{LauSL2}
Aaron~D. Lauda.
\newblock A categorification of quantum {${\rm sl}(2)$}.
\newblock {\em Adv. Math.}, 225(6):3327--3424, 2010.

\bibitem{LibLL}
Nicolas Libedinsky.
\newblock Sur la cat\'egorie des bimodules de {S}oergel.
\newblock {\em J. Algebra}, 320(7):2675--2694, 2008.

\bibitem{Morr02}
Scott Morrison.
\newblock A formula for the {J}ones-{W}enzl projections.
\newblock Preprint, 2002.
\newblock arXiv:1503.00384.

\bibitem{Morr07}
Scott Morrison.
\newblock A diagrammatic category for the representation theory of
  {$U_q(\mathfrak{sl}_n)$}.
\newblock Preprint, 2007.
\newblock arXiv:0704.1503.

\bibitem{TemLie}
H.~N.~V. Temperley and E.~H. Lieb.
\newblock Relations between the ``percolation'' and ``colouring'' problem and
  other graph-theoretical problems associated with regular planar lattices:
  some exact results for the ``percolation'' problem.
\newblock {\em Proc. Roy. Soc. London Ser. A}, 322(1549):251--280, 1971.

\bibitem{Wenzl}
Hans Wenzl.
\newblock On sequences of projections.
\newblock {\em C. R. Math. Rep. Acad. Sci. Canada}, 9(1):5--9, 1987.

\bibitem{WestG2}
Bruce~W. Westbury.
\newblock Enumeration of non-positive planar trivalent graphs.
\newblock {\em J. Algebraic Combin.}, 25(4):357--373, 2007.

\bibitem{WestCC}
Bruce~W. Westbury.
\newblock Invariant tensors and cellular categories.
\newblock {\em J. Algebra}, 321(11):3563--3567, 2009.

\bibitem{WestWebs}
Bruce~W. Westbury.
\newblock Web bases for the general linear groups.
\newblock {\em J. Algebraic Combin.}, 35(1):93--107, 2012.

\bibitem{WillSingular}
Geordie Williamson.
\newblock Singular {S}oergel bimodules.
\newblock {\em Int. Math. Res. Not. IMRN}, (20):4555--4632, 2011.

\end{thebibliography}

\end{document}